\documentclass[12pt,a4paper]{amsart}
\usepackage[english]{babel}

\newcommand{\numberset}{\mathbb}	
\newcommand{\R}{\numberset{R}}	
\newcommand{\Poly}{\numberset{P}}	
\renewcommand{\vec}{\mathbf}		
\newcommand{\lambdaB}{\boldsymbol{\lambda}}

\newcommand{\vth}{V_{th}}
\newcommand{\pot}{\varphi}
\newcommand{\ecpot}{\varphi^{ec}}
\newcommand{\Vbi}{\pot_{bi}}

\newcommand{\divt}{\text{div}\,}

\newcommand{\dx}[1]{\frac{\partial #1}{\partial x}}

\newcommand{\kb}{k_B}

\newcommand{\doping}{P}
\newcommand{\cref}{c_{ref}}
\newcommand{\Vapp}{V_{app}}

\newcommand{\vsat}{v_{sat}}
\newcommand{\po}{\partial \Omega}

\newcommand{\pk}{\partial K}
\newcommand{\Vtf}{w}
\newcommand{\Qtf}{q}
\newcommand{\Ltf}{\xi}
\newcommand{\Qbf}{\psi}

\newcommand{\Th}{\mathcal{T}_h}
\newcommand{\pec}{\mathcal{P}e}
\newcommand{\DX}{\,dx}
\newcommand{\DS}{\,d\sigma}
\newcommand{\dtk}{D} 
\newcommand{\vtk}{v} 
\newcommand{\ctk}{c} 
\newcommand{\gtk}{g} 

\newcommand{\numel}{N_{el}}

\newtheorem{remark}{Remark}
\newtheorem{theorem}{Theorem}
\newtheorem{definition}{Definition}
\newtheorem{proposition}{Proposition}

\providecommand{\norm}[1]{\lVert#1\rVert}
\providecommand{\abs}[1]{\lvert#1\rvert}
\newcommand*\widebar[1]{%
  \hbox{%
    \vbox{%
      \hrule height 0.7pt 
      \kern0.5ex
      \hbox{%
        \kern-0.1em
        \ensuremath{#1}%
        \kern-0.1em
      }%
    }%
  }%
} 

\newcommand{\unit}[1]{\mathrm{#1}}

\usepackage{amssymb}

\usepackage{amsmath}
\usepackage{amsfonts}
\usepackage{amssymb}
\usepackage{empheq} 

\usepackage{subfigure}
\usepackage[margin=3cm]{geometry} 

\usepackage{url}

\usepackage{graphicx}
\usepackage{color} 
\usepackage{transparent} 

\usepackage[T1]{fontenc} 
\usepackage{microtype}
\usepackage{lmodern}

\usepackage{siunitx}

\usepackage{rotating}
\usepackage{longtable}
\usepackage{multirow}
\usepackage{booktabs}
\usepackage{caption}
\usepackage{tabularx}
\newcolumntype{x}[1]{>{\arraybackslash\hspace{0pt}}p{#1}}

\newcommand{\rs}[1]{\textcolor{black}{#1}}

\graphicspath{{figures/}}

\newcommand{\Frac}{\displaystyle\frac}

\begin{document}

\title[Thermo-Fluid-Electrochemically Driven Ion Flow]
{Modeling and Simulation of Thermo-Fluid-Electrochemical 
Ion Flow in Biological Channels}

\author{Riccardo Sacco$^{1}$ \and Fabio Manganini$^2$ 
\and Joseph W. Jerome$^3$}

\address{$^{1}$ Dipartimento di Matematica,
               Politecnico di Milano \\
	       Piazza L. da Vinci 32, 20133 Milano, Italy \\
$^{2}$ Istituto di Matematica Applicata e 
Tecnologie Informatiche, \\
CNR, Via E. Bassini, 15, 20133 Milano, Italy \\
$^{3}$ Northwestern University, Mathematics Department\\
	2033 Sheridan Road Evanston, IL 60208-2730, USA}

\date{\today}

\begin{abstract}
In this article we address the study of ion charge transport
in the biological channels separating the intra and extracellular
regions of a cell. The focus of the investigation is devoted to
including thermal driving forces in the well-known 
velocity-extended Poisson-Nernst-Planck (vPNP) electrodiffusion model.
Two extensions of the vPNP system are proposed: the velocity-extended 
Thermo-Hydrodynamic model (vTHD) and the velocity-extended Electro-Thermal
model (vET). Both formulations are based on the principles of conservation
of mass, momentum and energy, and collapse into the vPNP model under 
thermodynamical equilibrium conditions. Upon introducing a suitable 
one-dimensional geometrical representation of the channel, we discuss
appropriate boundary conditions that depend 
only on effectively accessible measurable quantities. Then, we 
describe the novel models, the solution map used to iteratively solve them, 
and the mixed-hybrid flux-conservative stabilized finite element
scheme used to discretize the linearized equations.
Finally, we successfully 
apply our computational algorithms to the simulation
of two different realistic biological channels: 1) the Gramicidin-A channel considered in~\cite{JeromeBPJ}; and 
2) the bipolar nanofluidic diode considered in~\cite{Siwy7}.
%
%
\end{abstract}

\maketitle

\section{Introduction and Motivation}\label{sec:introduction}

Ion channels are ubiquitous in the cells 
of the human body. They allow communication 
between the intra and extra-cellular sites
and are responsible for the signaling pathways that 
regulate the activity of each part of the body system.
A complete presentation of ionic channels is far beyond the scope
of this article. For a treatment of the general properties of ionic
channels and of their characterization in cellular biology
and neuroscience we refer to the books~\cite{Hille2001} 
and~\cite{ErmentroutBook}. In broad terms, ion channel 
are \lq\lq biological gates\rq\rq{} that can be activated (opened) 
or deactivated (closed) by the application of electrical, 
chemical, mechanical or thermal stimuli.

In the present article we mainly focus on the interaction among 
electrical, chemical and thermal driving forces, starting from the
evidence that the sensory system in mammals is capable of 
discriminating thermal stimuli ranging from noxious cold 
(< 8$^\circ\,$C) to noxious heat (> 52$^\circ\,$C).
Of particular relevance are the biological temperature 
sensors denoted thermo-transient receptor potential (thermo-TRP) 
cation channels.
Thermo-TRP channels constitute a \lq\lq superfamily\rq\rq{} of receptors in sensory neurons that can be activated by noxious stimuli, resulting in the transmission to the spinal cord and brain
of the information of local pain perception~\cite{Fields1987}. 
The fact that the thermal thresholds of many thermo-TRP channels 
can be modulated by extracellular mediators has recently led to the development of antagonists for noxious channel blockage in therapeutic uses as novel analgesics~\cite{vay2012}.
This is the case, for example, of heat-sensitive TRP channels exposed to capsaicin, a natural ingredient of spicy foods such as 
red hot chilli peppers 
and of cold-sensitive TRP channels exposed to 
menthol.
In~\cite{caterina1997} and~\cite{voets2004}, it is shown how the chemical agonists capsaicin and menthol function for both types of TRP channels 
as gating modifiers, shifting activation curves towards physiological membrane potentials.
 
\medskip

Despite a significant amount of experimental data, the mechanisms 
underlying the marked temperature sensitivity of channel gating are still largely unknown. For this reason, the use of mathematical tools
is increasingly becoming popular to support biophysical conjectures and suggest novel theories for data interpretation.

\medskip

In this context, the most widely adopted approach to ion channel
modeling and simulation is represented by compartmental models 
described by a system of ordinary differential equations (ODEs)
based on the solution of Kirchhoff current law (KCL) written at the 
cellular membrane level~\cite{KeenerSneyd1998}.
The KCL equation is supplemented by phenomenological expressions characterizing
the input-output functional response of each protein channel to changes in
ion concentrations and/or electric potential inside and outside the cell~\cite{Hille2001}.
Temperature in these models is usually assumed to be a given parameter.

More sophisticated approaches involve the solution of a system of partial 
differential equations (PDEs) expressing balance of mass of each single ion species flowing across the channel and the Gauss law for the electric field.
The system is supplemented by a transport relation, known as the Nernst-Planck
(NP) equation, that describes ion motion under an electrochemical 
gradient~\cite{Rubinstein1990,Jerome2002,hess2005}. Also in this kind of modeling, well-known as the Poisson-Nernst-Planck (PNP) system and as the Drift-Diffusion (DD) system~\cite{markowich1986stationary,JeromeBook}, 
temperature is a given parameter. 

\medskip

A significant step forward to account for temperature as a dependent variable 
was taken in~\cite{JeromeBPJ}. In this reference
a hydrodynamic (HD) formulation including convective and thermal energy
in the electro-chemical motion of a single cation is proposed and numerically investigated. A remarkable feature of~\cite{JeromeBPJ} is
that model and analysis are inspired and guided by the 
analogy between a biological ion channel and the channel of a 
semiconductor device in which electrons and holes, 
instead of charged ion particles, flow to transport electrical
current between device terminals. This similarity between 
biology and solid-state electronics has been 
thoroughly addressed in the overview paper~\cite{eisenberg_living_transistor},
in the numerical simulations of~\cite{Hess2000,hess2005,Pandey2007} 
and in the MSc thesis~\cite{manganini2013}. 
In this latter reference, the mathematical view 
of~\cite{JeromeBPJ} is generalized by the introduction 
of a hierarchical modeling perspective
to represent ion transport in a biological channel
in which the interstitial electrolyte fluid is assimilated 
to the semiconductor device medium where two monovalent species
(anion and cation) are flowing under the effect of electric, chemical
and thermal forces. The hierarchy proposed in~\cite{manganini2013} 
is based on the ideas discussed in the semiconductor modeling
reference books~\cite{selberherr},~\cite{Markowich},~\cite{JeromeBook} and~\cite{Juengel2009}, and includes
four members: the basic DD formulation, the electro-thermal (ET) 
formulation, the hydrodynamic (HD) formulation and the thermo-hydrodynamic
(THD) formulation. For each member of the hierarchy, a set of conservation
laws for mass, momentum and energy is written. In the case of DD, ET and
HD models, the energy exchange between particles and medium are neglected
while in the case of the THD model a supplementary conservation law is
added to account for the dynamical thermo-electrochemical balance among
the three interacting subsystems. The hierarchy is extensively investigated 
in a series of numerical computations performed in a simplified one-dimensional
channel geometry. Externally applied data on which a sensitivity analysis
is carried out are the values of bulk ion concentrations in the intra and
extra-cellular sites and the applied potential drop across the channel.
Simulations indicate that: {\em i}) channel heating 
is, in general, relatively small compared to ion heating; 
and {\em ii}) for certain ranges of
model parameters, the high-order effects introduced by the THD
picture can significantly affect the input-output transfer characteristics of the \lq\lq biological transistor\rq\rq{}. 

\medskip

Based on the experience and results described above,
in the present article we propose the following mathematical
structure for the modeling of ion charge transport in biological
channels:
\begin{description}
\item[(a)] in Section~\ref{sec:PNP_model} we review 
the classic DD (PNP) model for ion electrodiffusion;
\item[(b)]
in Section~\ref{sec:geometry} we use the channel conformation model proposed 
in~\cite{Flux_coupling} to supply the DD (PNP) formulation with 
a set of \lq\lq reduced order\rq\rq{} boundary conditions that allow one to 
account for the experimentally accessible values of applied bias, 
ion concentrations and bathing temperatures;
\item[(c)] in Section~\ref{sec:extensionPNP} 
we review two members of the hierarchy of models for ion charge transport introduced in~\cite{manganini2013}, the Thermo-Hydrodynamic
and Electro-Thermal systems, and include electroosmotic effects by adding a 
linear advective term into the momentum balance equations
proportional to the given electrolyte fluid velocity. 
The resulting modified 
hierarchy is a family of (approximated) velocity-extended models, 
as the self-consistent study of electrolyte fluid dynamics 
through Navier-Stokes equations is neglected unlike 
in~\cite{Rubinstein1990,Jerome2002,Sacco2009,Schmuck}; 
\item[(d)] in Section~\ref{sec:solution_map} we illustrate the 
solution map used to solve iteratively the models 
in (a), (b) and (c) in steady-state conditions;
\item[(e)] in Section~\ref{sec:mod_prob} we study 
a linear model advective-diffusive and reactive 
boundary value problem (BVP) that represents each of the 
subproblems arising in the solution map, focusing
on the continuous maximum principle;
\item[(f)] in Section~\ref{sec:FE_discretization} 
we address the numerical study of the model BVP in (e)
based on the dual-mixed hybridized (DMH) finite element method
proposed and analyzed in~\cite{ArnoldBrezzi,BrezziFortin1991,RobertsThomas1991} 
for elliptic equations and extended to more general differential
operators in~\cite{ArbogastChen1995};
\item[(g)] in Section~\ref{sec:DMP} we propose a stabilization 
of the DMH scheme to deal with the case where advection
or reaction dominate over diffusion in the model BVP;
\item[(h)] in Section~\ref{sec:estimates} we conduct a series of numerical
tests to verify the convergence and stability properties of the DMH method;
\item[(i)] in Section~\ref{sec:numerical_results} 
we apply external temperature gradients, besides the usual 
electrochemical gradients, to the
biological transistor to replicate {\em in vitro} the biophysical situations
occurring {\em in vivo} in thermo-TRP channels;
\item[(j)] in the concluding Section~\ref{sec:conclusions} we summarize
the main contributions of the present research and we address 
a list of forthcoming activities for model improvement.
\end{description}

\section{The Poisson-Nernst-Planck model for ion electrodiffusion}
\label{sec:PNP_model}

The Poisson-Nernst-Planck system (PNP) for ion electrodiffusion 
reads~\cite{Rubinstein1990}:
\begin{subequations}\label{eq:PNP_flux}
\begin{align}
& \dfrac{\partial c_{i}}{\partial t}+\text{div}\,\mathbf{f}_{i}
\left(c_{i},\varphi\right) = 0 && i=1,\ldots,M
\label{eq:continuity-1-1}
\\
& \mathbf{f}_{i}\left(c_{i},\varphi\right)
= -D_{i}\nabla c_{i} + \mu_{i} \Frac{z_{i}}{|z_i|} c_{i}\mathbf{E}
&& i = 1,\ldots,M
\label{eq:nernstplanck-1}
\\
& D_{i} = \dfrac{\mu_{i}V_{th}}{\left|z_{i}\right|} && i=1,\ldots,M
\label{eq:einstein_rel-1}
\\
& \text{div}\,\mathbf{E} = \dfrac{\rho}{\varepsilon}
\label{eq:poisson-1-1}
\\
& \mathbf{E} =  -\nabla\varphi \label{eq:cel_2}
\\
& \rho = q \sum_{i=1}^M z_{i}c_{i} + q \doping.
\label{eq:rho}
\end{align}
In the PNP system,~\eqref{eq:continuity-1-1} is the continuity equation
describing mass conservation for each ion whose concentration
is denoted by $c_i$ (\si{\metre^{-3}}), $i=1, \ldots, M$, $M \geq 1$ being the number
of ions flowing in the electrolyte fluid.
Each ion flux density $\mathbf{f}_{i}$ (\si{\metre^{-2} s^{-1}}) 
is defined by the Nernst-Planck relation~\eqref{eq:nernstplanck-1}
in which it is possible to recognize a chemical contribution and an electric contribution, in such a way that the model be regarded as an extension of 
Fick's law of diffusion to the case where the diffusing particles are 
also moved by electrostatic forces with respect to the fluid. 
The quantity $z_i$ is the valence of the $i$-th ion, while 
$\mu_{i}$ and $D_{i}$ are the mobility and diffusivity of the chemical 
species, respectively, related by the Einstein relation \eqref{eq:einstein_rel-1} where $V_{th}= k_{B} T_{sys}/q$ (\si{V}) is the thermal potential, our having denoted by $k_{B}$, $T_{sys}$ and $q$, the Boltzmann constant, the absolute temperature of the system 
and the elementary charge, respectively.
We remark that in the PNP modeling approach, $T_{sys}$ is a parameter 
representing the constant temperature of the system, without distinction 
between the ionic species and the electrolyte fluid. This assumption 
is going to be relaxed in Section~\ref{sec:extensionPNP}.
The function $\rho$ is the space charge density (\si{C m^{-3}})
in the electrolyte and is given by the sum of two contributions, 
the mobile ion charge $q \sum_{i=1}^M z_{i}c_{i}$ and the fixed 
charge $q \doping$.
The electric field $\mathbf{E}$ (\si{V m^{-1}}) 
due to space charge distribution $\rho$ in the electrolyte is determined 
by the Poisson equation~\eqref{eq:poisson-1-1} which represents Gauss' 
law in differential form, $\varepsilon$ being the dielectric permittivity of 
the electrolyte fluid medium.
For further development, it is useful to introduce 
the electrical current density 
$\mathbf{j}_{i}$ (\si{A \metre^{-2}}), equal to the 
number of ion charges flowing through a given surface area per unit time
and defined as
\begin{equation}
\mathbf{j}_{i}:=q z_i \mathbf{f}_{i} = 
- q z_i D_i \nabla c_i + q \mu_i |z_i| c_i \mathbf{E} \qquad 
i=1, \ldots, M.
\label{eq:current_density}
\end{equation}
\end{subequations}

\begin{remark}
The PNP system~\eqref{eq:PNP_flux} has the same format and structure
as the Drift-Diffusion (DD) equations for 
semiconductors (see, e.g.,~\cite{JeromeBook}),
but it is applied to a different medium (water instead of a semiconductor
crystal lattice) and includes, in general, 
more charge carriers than just holes and electrons, as in the case
of semiconductor device theory.
\end{remark}

\section{Geometry, biophysical assumptions and boundary
conditions}\label{sec:geometry}

\begin{figure}[h!]
\begin{center}
\includegraphics[width=0.65\linewidth]{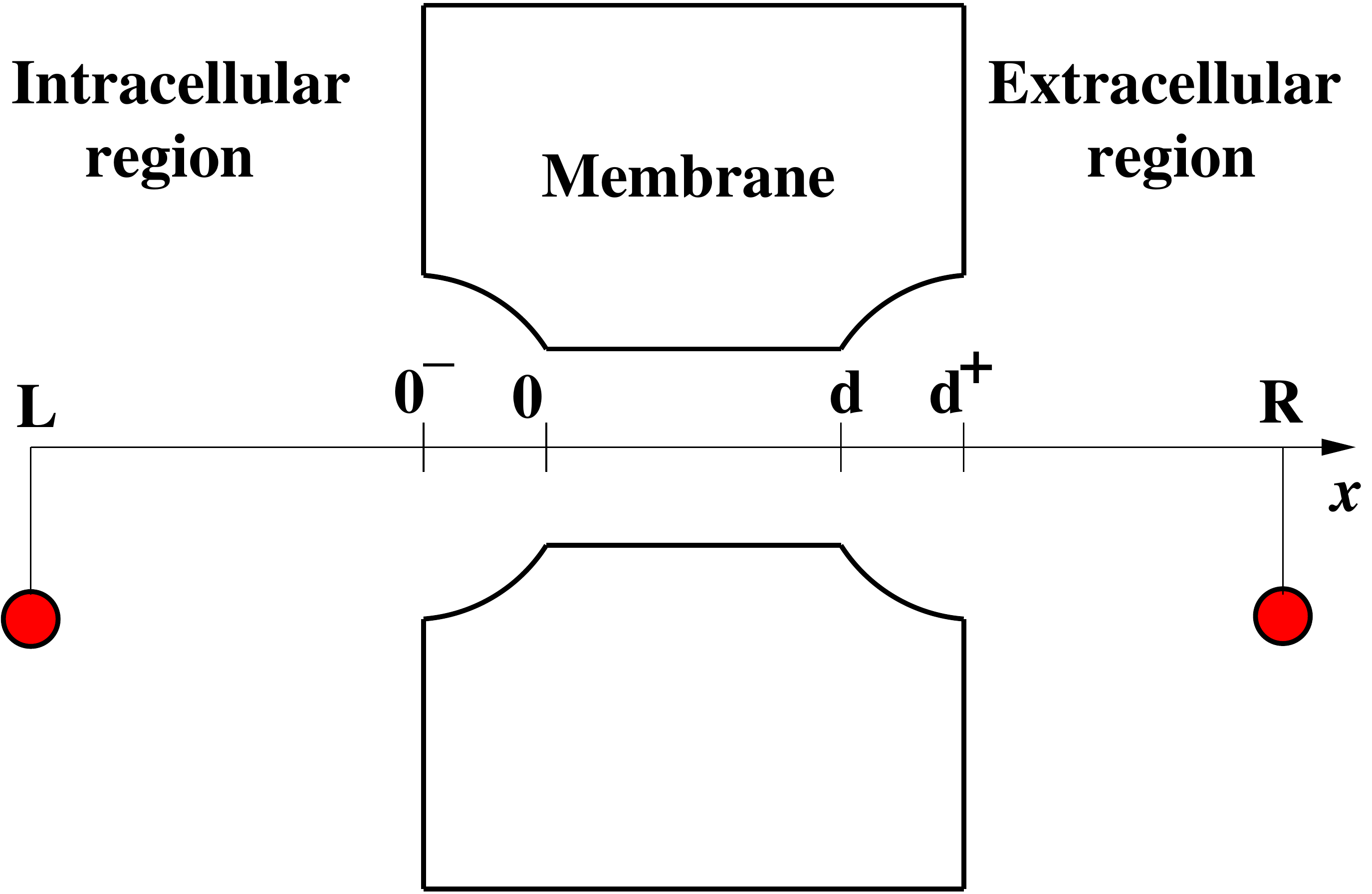}
\end{center}
\caption[channel: schematic view]{Schematic view of the biophysical
problem. Five regions can be distinguished. From left to right: 
intracellular bathing solution, antichamber, channel region, 
antichamber, bathing solution in the extracellular side. At the 
endpoints, $x=0$ and $x=L$, terminal contacts are marked in red color.}
\label{fig:channel}
\end{figure}

Following the presentation of~\cite{Flux_coupling}, 
we illustrate in Fig.~\ref{fig:channel} the schematic representation
of a cross-section of a biological channel, assuming rotational invariance
around the channel axis ($x$ axis).

Five regions can be distinguished, ordered from left to right:
\begin{enumerate}
\item Region 1, $L \leq x \leq 0^-$: this is the bathing solution 
in the intracellular side.
\item Region 2, $0^- < x \leq 0$: this is the channel antichamber 
(or access region) from the intracellular side. 
\item Region 3, $0 \leq x < d$: this is the channel region.
\item Region 4, $d \leq x < d^+$: this is the channel antichamber
from the extracellular side.
\item Region 5, $d^+ \leq x \leq R$: this is the bathing solution 
in the extracellular side.
\end{enumerate} 
The above geometrical representation leads naturally
to a multi-domain formulation for ion channel
simulation. The biophysical treatment of such a problem is fully
carried out in~\cite{Flux_coupling} whilst its mathematical and
numerical treatment is the object of the present article, where 
the PNP formulation~\eqref{eq:PNP_flux} is extended to include thermo-hydrodynamical phenomena to ion electrodiffusion in the channel.

\subsection{Assumptions}\label{sec:assumptions}
The endpoints of the domain, $x=L$ and $x=R$, are located 
sufficiently far from the antichamber and channel regions, in such a way
that appropriate equilibrium conditions can be applied.
More importantly, the endpoints are assumed to be the physical 
place where the solution intra and extra-cellular electrochemical conditions 
are accessible to experimental measurements.
Because of this, following~\cite{one_conformation,Flux_coupling}, 
we assume that at $x=L$ and $x=R$:
\begin{subequations}\label{eq:BCs_endpoints}
\begin{description}
\item[(\textbf{A1})] 
the electric potential $\pot$ is a known given quantity, so that:
\begin{align}
& \pot(L) = \pot_L \label{eq:phi_L} \\
& \pot(R) = \pot_R; \label{eq:phi_R}
\end{align}
\item[(\textbf{A2})] 
the ion concentrations $c_i$ are known given quantities, so that:
\begin{align}
& c_i(L) = c_{i,L}, \qquad i=1, \ldots, M \label{eq:c_i_L} \\
& c_i(R) = c_{i,R}, \qquad i=1, \ldots, M \label{eq:c_i_R}
\end{align}
where $M \geq 1$ is the number of ions flowing in the cellular solution.
The boundary values for the ion concentrations satisfy the electroneutrality
constraint:
\begin{align}
& \sum_{i=1}^M z_i c_{i,L} = 0 \label{eq:electroneutrality_L} \\
& \sum_{i=1}^M z_i c_{i,R}= 0 \label{eq:electroneutrality_R}
\end{align}
where $z_i$ is the charge number associated with each ion species
$c_i$ ($z_i>0$ for cations, $z_i<0$ for anions and $z_i=0$ for neutral
species).
\end{description}
We close the characterization of the bathing regions by assuming that:
\begin{description}
\item[(A3)] the ion flux densities vanish inside the baths 
\begin{align}
& J_i(x) = 0	\qquad L \leq x \leq 0^-, \qquad 
i=1, \ldots, M \label{eq:J_i_L} \\
& J_i(x) = 0	\qquad d^+ \leq x \leq R, \qquad 
i=1, \ldots, M. \label{eq:J_i_R}
\end{align}
\end{description}

Continuing to move from the periphery of the domain towards the channel
region, we encounter the antichamber openings, at $x=0^-$ and
$x=d^+$, respectively. Located on the membrane lipid bilayer, a 
space-dependent fixed charge density $\doping=\doping(x)$ is 
distributed. We make the following assumptions within the antichamber
regions:
\begin{description}
\item[(A4)] 
electroneutrality holds at the channel mouth entrances:
\begin{align}
&q \doping (0^-) + q\sum z_j c_j(0^-) = 0 \label{eq:neutrality_BathL}\\
&q \doping (d^+) + q\sum z_j c_j(d^+) = 0; \label{eq:neutrality_BathR}
\end{align}
\item[(A5)] the electric potential is constant in the antichamber regions:
\begin{align}
& \varphi(x) = \varphi(0^-) \qquad 0^- \leq x \leq 0 
\label{eq:pot_constant_antich_L}\\
& \varphi(x) = \varphi(d^+) \qquad d^- \leq x \leq d^+; 
\label{eq:pot_constant_antich_R}
\end{align}
\item[(A6)] the ion concentrations are constant 
in the antichamber regions:
\begin{align}
& c_i(x) = c_i(0^-) \qquad 0^- \leq x \leq 0 \qquad
i=1, \ldots, M
\label{eq:c_i_constant_antich_L}\\
& c_i(x) = c_i(d^+) \qquad d \leq x \leq d^+ \qquad
i=1, \ldots, M.
\label{eq:c_i_constant_antich_R}
\end{align}
\end{description}

\end{subequations}

\subsection{Boundary Conditions 
at Channel Openings}\label{sec:BCs}

In this section we use the geometrical multi-domain representation
of the problem of Sect.~\ref{sec:geometry} and the
assumptions made in Sect.~\ref{sec:assumptions} to derive 
the boundary conditions to be supplied to the PNP Model 
at $x=0$ and $x=d$ (channel openings). Using Einstein's 
relation~\eqref{eq:einstein_rel-1} in~\eqref{eq:current_density}
we can write the PNP current density as
\begin{subequations}\label{eq:BCs}
\begin{align}
& J_i = -q\abs{z_i}\mu_i c_i \dx{\ecpot_i} \qquad
i=1, \ldots, M \label{eq:pnp_flux}
\end{align}
where
\begin{align}\label{eq:ecpot}
& \ecpot_i := \pot + \frac{1}{z_i} \vth \ln \left( \frac{c_i}{\cref}\right)
\qquad i=1, \ldots, M 
\end{align}
is the electrochemical potential and 
$\cref := \displaystyle 
\max_{i=1, \ldots, M}\{c_{i,L}; c_{i,R}\}$ is a reference concentration.
Inverting~\eqref{eq:ecpot} we obtain the well-known 
Maxwell-Boltzmann (MB) statistics for the ion densities
\begin{align}\label{eq:MB_relations}
& c_i = \cref \exp\left( z_i \Frac{\ecpot_i - \pot}{\vth} \right)
\qquad i=1, \ldots, M.
\end{align}
Using (\textbf{A3}) we get 
\begin{align}
\ecpot_i(x) = const =: \overline{\ecpot_{i,L}} \qquad L \leq x \leq 0^- 
\qquad i=1, \ldots, M \label{eq:ec_const_L}\\
\ecpot_i(x) = const =: \overline{\ecpot_{i,R}} \qquad d^+ \leq x \leq R
\qquad i=1, \ldots, M \label{eq:ec_const_R}
\end{align}
where $\overline{\ecpot_{i,L}}$ and $\overline{\ecpot_{i,R}}$, 
$i=1, \ldots, M$, are constants yet to be determined.
From~\eqref{eq:MB_relations} we get
\begin{align}\label{eq:conc}
& c_i(x) = \cref \exp{ \left(z_i \frac{(\ecpot_i(x) - 
\pot(x))}{\vth}\right)} \qquad L \leq x \leq R \qquad i=1, \ldots, M.
\end{align}
Thus, using~\eqref{eq:ec_const_L} and~\eqref{eq:ec_const_R} and 
assumptions (\textbf{A5}) and (\textbf{A6}), we find:
\begin{align}
& c_i(x) = \cref \exp{\left( z_i \frac{\left( \overline{\ecpot_{i,L}} - \pot(x)\right)}{\vth}\right) } \qquad L \leq x \leq 0 
\qquad i=1, \ldots, M \nonumber \\
& c_i(x) = \cref \exp{\left( z_i \frac{\left( \overline{\ecpot_{i,R}} - \pot(x)\right)}{\vth}\right) } \qquad d \leq x \leq R 
\qquad i=1, \ldots, M. \nonumber
\end{align}
Using (\textbf{A2}) 
in the equations above at $z=L$ and $z=R$, respectively, we obtain
the boundary values for the electrochemical potential:
\begin{align}
\overline{\ecpot_{i,L}} = \pot_L + \frac{\vth}{z_i}\ln\left( \frac{c_{i,L}}{\cref} \right) \qquad i=1, \ldots, M 
\label{eq:ec_bc_L} \\[3mm]
\overline{\ecpot_{i,R}} = \pot_R + \frac{\vth}{z_i}\ln\left( \frac{c_{i,R}}{\cref} \right) \qquad i=1, \ldots, M.
\label{eq:ec_bc_R}
\end{align}
\end{subequations}

\begin{remark}
It is important to note that relations~\eqref{eq:ec_bc_L} 
and~\eqref{eq:ec_bc_R} allow one to express the electrochemical potentials at the contacts as a function of the sole accessible quantities. 
This is in contrast with~\cite{manganini2013}, page 42, 
relations (3.8a)-(3.8b), where the values of the 
electrochemical potentials were imposed and the 
corresponding values of the electrical potential 
were computed using (3.14a)-(3.14b). 
\end{remark}

From the previous discussion, we see that electroneutrality holds 
at the external contacts ($x=L$, $x=R$) and at the channel mouth
entrances ($x=0^-$, $x=d^+$), but not, in general, elsewhere.
Therefore, it makes sense to introduce the voltage drops occurring 
in the bathing solution regions, the so-called built-in potentials:
\begin{subequations}\label{eq:builtin}
\begin{align}
& \Vbi(0^-) := \pot(L) - \pot(0^-) \label{eq:vbi_L} \\
& \Vbi(d^+) := \pot(d^+) - \pot(R).  \label{eq:vbi_R}
\end{align}

Using the charge neutrality conditions~\eqref{eq:neutrality_BathL} and~\eqref{eq:neutrality_BathR}, the MB relations~\eqref{eq:MB_relations}
and the definitions~\eqref{eq:builtin}, the built-in potentials 
can be determined by solving the two following nonlinear algebraic
equations:
\begin{align}
& \doping (0^-) + \sum_{i=1}^M z_i c_{i,L} 
\exp{\left( z_i \frac{ \Vbi(0^-)}{\vth}\right) } = 0 
\label{eq:BI_L} \\
& \doping (d^+) + \sum_{i=1}^M z_i c_{i,R} 
\exp{\left( -z_i \frac{ \Vbi(d^+)}{\vth}\right) } = 0.
\label{eq:BI_R}
\end{align}

\begin{remark}
In the particular case of the $KCl$ solution ($z_1 = +1,\,z_2=-1$) 
considered in~\cite{one_conformation}, 
equations~\eqref{eq:BI_L} and~\eqref{eq:BI_R} 
can be explicitly solved, to yield:
\begin{align}
&\Vbi(0^-) = \vth \ln \left( \frac{-\doping(0^-) + \sqrt{(\doping(0^-))^2 + 4(c_{1,L})^2}}{2c_{1,L}} \right) \nonumber \\
&\Vbi(d^+) = \vth \ln \left( \frac{\doping(d^+) + \sqrt{(\doping(d^+))^2 + 4(c_{1,R})^2}}{2c_{1,R}} \right) \nonumber
\end{align}
where we notice that $c_{1,L}=c_{2,L}$ and $c_{1,R}=c_{2,R}$
because of the electroneutrality
constraints~\eqref{eq:electroneutrality_L} and~\eqref{eq:electroneutrality_R}.
\end{remark}
\end{subequations}

In order to complete the characterization of the boundary values for the dependent variables of the problem we need the value of the 
electric potential and of the ion concentrations at the two endpoints 
of the channel. The potential is determined 
using~\eqref{eq:builtin} and (\textbf{A5}), which yield:
\begin{subequations}\label{eq:channel_BCS}
\begin{align}
& \pot(0) = \pot(0^-) = \pot_L - \Vbi(0^-) \label{eq:POT_0} \\
& \pot(d) = \pot(d^+) = \pot_R + \Vbi(d^+). \label{eq:POT_d}
\end{align}
To determine the ion concentrations we use~\eqref{eq:conc},~\eqref{eq:ec_bc_L}
and~\eqref{eq:vbi_L} at $x=0$ to obtain
\begin{align}
c_i(0) & = \cref \exp{ \left(z_i \frac{(\ecpot_i(0) - \pot(0))}
{\vth}\right)} = c_{i,L} 
\exp{\left( z_i \frac{\Vbi(0^-)}{\vth}\right)} 
\qquad i=1, \ldots, M
\label{eq:CONC_0} 
\end{align}
and~\eqref{eq:conc},~\eqref{eq:ec_bc_R}
and~\eqref{eq:vbi_R} at $x=d$ to obtain 
\begin{align}
c_i(d) & = \cref \exp{ \left(z_i \frac{(\ecpot_i(d) - \pot(d))}{\vth}\right)} = c_{i,R} 
\exp{\left( -z_i \frac{\Vbi(d^+)}{\vth}\right)} 
\qquad i=1, \ldots, M. \label{eq:CONC_d} 
\end{align}
\end{subequations}

\null
\noindent
Summarizing, the Dirichlet boundary conditions 
for the PNP system in the reduced channel 
domain $\Omega = (0, d)$ are at $x=0$:
\begin{subequations}\label{eq:bc_summary}
\begin{align}
& \doping (0^-) + \sum_{i=1}^M z_i c_{i,L} \exp{\left( z_i \frac{\Vbi(0^-)}{\vth}\right)} = 0 & \label{eq:bi_val_0} \\
&\pot(0) = \pot_L - \Vbi(0^-) & \label{eq:bc_pote_0}\\
&\ecpot_i(0) = \pot_L +\frac{1}{z_i}\vth 
\ln\left( \frac{c_{i,L}}{\cref} \right) & i=1, \ldots, M 
\label{eq:bc_ec_0} \\
&c_i(0) = c_{i,L}\exp{\left( z_i \frac{\Vbi(0^-)}{\vth} \right)}
& i=1, \ldots, M
\label{eq:bc_conc_0}
\end{align}
and, at $x=d$:
\begin{align}
&\doping (d^+) + \sum_{i=1}^M z_i c_{i,R} \exp{\left( z_i \frac{\Vbi(d^+)}{\vth}\right)} = 0 & \label{eq:bi_val_d} \\
&\pot(d) = \pot_R + \Vbi(d^+) &\label{eq:bc_pote_d}\\
&\ecpot_i(d) = \pot_R + 
\frac{1}{z_i}\vth \ln\left( \frac{c_{i,R}}{\cref} \right) 
& i=1, \ldots, M
\label{eq:bc_ec_d} \\
&c_i(d) = c_{i,R} \exp{\left( -z_i \frac{\Vbi(d^+)}{\vth} \right)}
& i=1, \ldots, M.
\label{eq:bc_conc_d}
\end{align}
\end{subequations}

\section{Extensions of the PNP model}\label{sec:extensionPNP}
In this section we present two extensions 
of the PNP equation system illustrated
in Sections~\ref{sec:PNP_model} and~\ref{sec:geometry}.
From now on we restrict our attention to the study of a channel
in stationary conditions and to simplify the exposition, 
in close analogy with semiconductor device physics, 
we consider a binary mixture of \emph{monovalent} anions and cations, 
i.e., such that their chemical valence is equal to $\pm1$
as in the case of the $KCl$ solution or the $NaCl$ solution.
We use the symbol $p$ (positive) and $n$ (negative) to refer to the concentrations of cations and anions, respectively.
The symbol $\nu$ is used to indicate either $p$ or $n$. 
When the operators $\mp$ or $\pm$ are used, the upper sign always refers 
to $n$-type particles, the lower sign to $p$-type particles, and the
value for $\nu$ must be chosen accordingly. The quantities associated 
with the electrolyte medium are referred to by the symbol $e$, 
while $x$ denotes henceforth the spatial coordinate.

In our analysis we follow~\cite{manganini2013}, the series of references~\cite{JeromeBook,Baccarani,MomentsRudan}
in the context of semiconductor device modeling, 
and the theory of~\cite{Rubinstein1990,Jerome2002,Sacco2009,Schmuck}.

The first extension is denoted velocity-extended Thermo-Hydrodynamic model (vTHD) and the second extension is denoted velocity-extended Electro-Thermal model (vET). Both extensions include thermal driving forces 
and the translational contribution of electrolyte fluid flow, 
besides the usual electrochemical forces accounted for by the PNP equations.
In the vTHD model the temperature of anions, cations and 
lattice are considered as distinct dependent variables, 
while in the vET model there is only one temperature, generally 
varying with the spatial coordinate $x$, to describe the global system.

Finally, in the present article we assume 
the electrolyte fluid velocity $v_e$ to be a given function of $x$. 
This is a strong simplification
which is going to be removed in a future publication.
Moreover, to self-consistently account for electrodynamical
effects, we adjoin to both vTHD and vET model the solution of the Poisson
equation, written below in conservation form:
\begin{subequations}\label{eq:POISSON}
\begin{align}
& \dx{E} = \Frac{q}{\epsilon}(p-n+\doping) \label{eq:poisson_1D}\\
& E = -\dx{\varphi}. \label{eq:efield_1D} 
\end{align}
\end{subequations}

\subsection{Thermodynamical Equilibrium}\label{sec:thermal_eq}
Let $T_p=T_p(x)$, $T_n=T_n(x)$ and 
$T_e=T_e(x)$ denote the temperatures
of cations, anions and electrolyte, respectively, with 
$x \in [0, d]$. Let also $T_{sys}$ be 
a given value representing the constant temperature 
of the whole environment (ions + fluid).
By definition, thermal equilibrium is the condition 
corresponding to applying no external electrical, mechanical, chemical
and thermal forces to the biophysical system. 
In such a case, system response is represented by the
following values of the dependent variables:
\begin{subequations}\label{eq:thermal_eq}
\begin{align}
& v_\nu(x) =0 & \label{eq:equilibrium_v_nu}\\
& v_e(x)  = 0 & \label{eq:equilibrium_v_e}\\
& T_p(x) = T_n(x) = T_e(x) = T_{sys}. & \label{eq:equilibrium_T}
\end{align}
Eq.~\eqref{eq:equilibrium_v_nu} expresses the fact that
the ion drift velocity is equal to zero, while 
Eq.~\eqref{eq:equilibrium_v_e} expresses the fact that
the fluid translational velocity is equal to zero. 
Eq.~\eqref{eq:equilibrium_T} expresses the fact that 
the system settles in an isothermal condition.
The three above conditions imply that: 
\begin{align}
& J_\nu(x) =0 & \label{eq:equilibrium_J_nu}\\
& S_\nu(x) = S_e(x) = 0 & \label{eq:equilibrium_S} 
\end{align}
where $J_p$ and $J_n$ are
anion and cation current flows, while $S_\nu$ and $S_e$ are
ion and fluid energy fluxes.
Using~\eqref{eq:pnp_flux} in~\eqref{eq:equilibrium_J_nu} 
tells us that at thermal equilibrium each ion
electrochemical potential is constant in $[0, d]$
\begin{align}
& \ecpot_\nu(x) = \ecpot_{\nu, eq} & \label{eq:equilibrium_ec_nu} 
\end{align}
where the constant value $\ecpot_{\nu, eq}$ is different for each ion
and can be computed as detailed in Sect.~\ref{sec:BCs}.
Model consistency requires that under thermal equilibrium 
any extension of the PNP system verifies all 
conditions~\eqref{eq:thermal_eq}. This fact is documented 
in the later sections where the vTHD and vET models are introduced.
\end{subequations}
 
\subsection{Velocity-extended Thermo-Hydrodynamic (vTHD) model}
\label{sec:TFHD}
The stationary velocity-extended thermo-hydrodynamic model in the 
one-dimensional setting consists of the following system of conservation laws~(see~\cite[Chap. 2]{selberherr} and~\cite{JeromeBook,Baccarani,MomentsRudan}):
\begin{subequations}
\label{eq:THD1D}
\begin{align}
& \frac{1}{q} \dx{J_\nu} = 0 \label{eq:THD1D_cont} \\
& J_\nu + \nu \tau_{p_\nu} v_\nu \dx{} 
\left( \dfrac{J_\nu}{\nu}\right) = \pm qD_\nu \dx{\nu} - q\mu_\nu \nu \dx{}\left( \varphi \mp \frac{\kb T_\nu}{q}\right) \mp q \nu v_e
\label{eq:THD_Mom_1D}\\ 
& \dx{S_\nu} = EJ_\nu - \frac{\nu}{\tau_{w_\nu}}(w_\nu - w_\nu^{eq}) 
\label{eq:THD1D_en}\\
& \dx{S_e} = \left[ \dfrac{p}{\tau_{w_p}}(w_p- w_p^{eq}) + \dfrac{n}{\tau_{w_n}}(w_n-w_n^{eq})	\right]
\label{eq:THD1D_en_L}\\
& S_\nu = -\kappa_\nu \dx{T_\nu} \mp \frac{J_\nu}{q}(w_\nu + \kb T_\nu) 
\label{eq:THD_s_1D} \\
& S_e = -\kappa_e \dx{T_e} + N_e v_e w_e
\label{eq:THD_sL_1D} \\ 
& w_\nu = \frac{1}{2} m_\nu |v_\nu|^2 + \frac{3}{2} \kb T_\nu
\label{eq:THD_W_1D}\\
& w_{e} = \frac{1}{2} m_e |v_e|^2 + \frac{3}{2} \kb T_e
\label{eq:THD_Weq1D} \\
& w_\nu^{eq} = \frac{3}{2} \kb T_e. \label{eq:THD_Weq_1D} 
\end{align}
\end{subequations}
In system~\eqref{eq:THD1D}, $\nu$ and $T_{\nu}$ represent 
ion density and temperature, $N_e$ and $T_e$ are 
electrolyte number density and 
temperature, while $v_\nu$,
$J_\nu$, $S_\nu$ and $S_e$ represent
the drift (or translational) velocity of each ion, the ion current density, the ion energy density flux
and the electrolyte energy density flux, respectively.
The drift velocities and the current densities are related
through the phenomenological law
\begin{subequations}
\begin{align}
& v_\nu = \mp \frac{J_\nu}{q \, \nu} \label{eq:THD_v_1D}
\end{align}
\end{subequations}

Eq.~\eqref{eq:THD1D_cont} expresses conservation of ion density, 
eq.~\eqref{eq:THD_Mom_1D} expresses conservation of ion momentum 
density, eqns.~\eqref{eq:THD1D_en} and~\eqref{eq:THD1D_en_L} 
express conservation of the energy densities for ions and electrolyte 
fluid, 
eqns.~\eqref{eq:THD_s_1D} and~\eqref{eq:THD_sL_1D} are constitutive laws for ion and electrolyte energy density fluxes $S_\nu$ and $S_e$,
respectively.

The quantities $\kappa_\nu$ and $\kappa_e$ are
the thermoconductivity coefficients of ions and electrolyte, respectively,
while $\tau_{p_\nu}$ and $\tau_{w_\nu}$ are 
the momentum and energy relaxation times, respectively.
Phenomenological expressions for these coefficients are: 
\begin{subequations}
\begin{align}
& \tau_{p_\nu} = \dfrac{m_\nu \mu_{0\nu}T_e}{qT_\nu}, & \label{eq:tau_p} \\
& \tau_{w_\nu} = \frac{3}{2} \frac{\mu_{0\nu} k_B T_\nu T_e}
{q \vsat^2(T_\nu+T_e)} + \frac{\tau_{p_\nu}}{2},& \label{eq:tau_w_vsat}\\
& \kappa_\nu = \frac{3}{2} \frac{ \mu_{0_\nu} k_B^2 T_e  } {q} \nu,
& \label{eq:kappa}
\end{align}
\end{subequations}
where $\mu_{0_\nu}$ is the low-field mobility and $\vsat$ is the saturation velocity. For a physical interpretation of the parameter $\vsat$ and its influence on system behavior, we refer to~\cite{JeromeBPJ, manganini2013}. 

The ion mobility
\begin{subequations}
\begin{align}
& \mu_\nu = \dfrac{q}{m_\nu}\tau_{p_\nu} \label{eq:mobility}
\end{align}
is related to the diffusion coefficient $D_\nu$ by the generalized Einstein relation
\begin{align}
& D_\nu = \mu_\nu \dfrac{\kb T_\nu}{q} = 
\dfrac{\kb T_\nu}{m_\nu}\tau_{p_\nu} \label{eq:Einstein}
\end{align}
where $m_\nu$ is the ionic mass.

Finally, $w_\nu$ and $w_e$ are the ion and electrolyte energies,
in which we can distinguish a kinetic contribution and a thermal 
contribution. The quantities $w_p^{eq}$ and $w_n^{eq}$ are the
equilibrium energies of cations and anions.
\end{subequations}

\begin{proposition}[Consistency with respect to thermal equilibrium]
The vTHD satisfies relations~\eqref{eq:thermal_eq}.
\end{proposition}
\begin{proof}
Using~\eqref{eq:equilibrium_v_nu},~\eqref{eq:equilibrium_v_e} and~\eqref{eq:equilibrium_T} in~\eqref{eq:THD_Mom_1D} 
and~\eqref{eq:Einstein} we get
$$
J_\nu = \pm q \mu_\nu \frac{\kb T_{sys}}{q} \dx{\nu} - 
q\mu_\nu \nu \dx{\varphi} = 0
$$
or, equivalently,
$$
J_\nu = - q \mu_\nu \nu \dx{}\left( \mp \vth 
\ln \left( \nu/\nu_{ref} \right) + \varphi \right) = 0.
$$
The above equation represents the particle flux equilibrium 
condition~\eqref{eq:equilibrium_J_nu} and can be written in the form~\eqref{eq:pnp_flux} upon defining the
equilibrium electrochemical potential
$$
\ecpot_\nu := \pot \mp \vth \ln \left( \frac{\nu}{\nu_{ref}}\right)
$$
which corresponds to~\eqref{eq:ecpot} having set $z_n=-1$ and $z_p=+1$.
Having proved~\eqref{eq:equilibrium_J_nu} 
(and, equivalently,~\eqref{eq:equilibrium_ec_nu}) and 
using~\eqref{eq:equilibrium_T} in~\eqref{eq:THD_s_1D} 
and~\eqref{eq:THD_sL_1D}, we immediately obtain 
$$
S_n = S_p = S_e = 0
$$
that is the energy flux equilibrium 
condition~\eqref{eq:equilibrium_S}. We notice that at thermal equilibrium 
the ion 
and fluid energies coincide with their rest energy
$$
w_n = w_p = w_e = \frac{3}{2} \kb T_{sys}.
$$
\end{proof}

\subsection{Velocity-extended Electro-Thermal (vET) model}\label{subsec:ET}
This model differs quite significantly from the 
vTHD because of the following simplifying assumptions:
\begin{description}
\item[(H1)]
cations, anions and electrolyte are in local equilibrium, i.e.
\begin{subequations}\label{eq:assumptions_ET}
\begin{align}
&	T_n = T_p = T_e =: T(x).
\label{eq:ET_T_assumption}
\end{align}
This situation occurs when the system is assumed to have had enough time to evolve and to reach steady state;
\item[(H2)] the kinetic energy is small compared to the thermal energy, i.e. 
\begin{align}
&	\frac{1}{2}m_\nu |v_\nu|^2 \ll \frac{3}{2}\kb T.
\label{eq:ET_kinetic_assumption}
\end{align} 
Thus we can write
\begin{align}
& w_\nu \simeq \frac{3}{2}\kb T
\label{eq:ET_energy}
\end{align} 
\end{subequations}
and neglect the nonlinear inertial term
in the momentum balance equation~\eqref{eq:THD_Mom_1D}.
\end{description}

Using (H1), 
and summing up the three resulting energy balance equations
we obtain the following equations of the velocity-extended 
ET system:
\begin{subequations}\label{subeq:ET}
\begin{align}
& \frac{1}{q}\divt J_\nu = 0
\label{eq:ET_Cont}\\
& J_\nu = q \mu_\nu \nu E \pm \mu_\nu \dx{}(\nu \kb T) \mp q \nu v_e
\label{eq:ET_Mom}\\
& \divt S_{tot} = E \cdot (J_n + J_p) - \divt \left( \frac{5}{2}\dfrac{\kb T}{q} \left( J_p - J_n\right) \right)
\label{eq:ET_En} \\
& S_{tot} =
\left( -\kappa_{tot} \dx{ T } +  \frac{3}{2}\kb T N_e v_e \right)
\label{eq:ET_S_tot} \\
& \kappa_{tot} = \kappa_e + \kappa_n + \kappa_p.
\label{eq:ET_kappa_total}
\end{align}

\begin{remark}
Looking at the (simplified) momentum conservation 
equation~\eqref{eq:ET_Mom}, we see that it now furnishes 
an \emph{explicit} expression for the current density 
as a function of concentration, 
electric potential and temperature. This allows one 
to adopt the basic ideas of the Scharfetter-Gummel discretization 
method~\cite{SG}, which is the most widely used numerical 
scheme in contemporary semiconductor device simulation 
because of its remarkable stability and accuracy.
\end{remark}

\begin{remark}
Simple manipulations of~\eqref{eq:ET_Mom} show that it can be written
in the following equivalent DD form
\begin{align}
& J_\nu = \pm q D_\nu \dx{\nu} + q \mu_\nu \nu E_\nu 
\label{eq:ET_Mom_DD}
\end{align}
where 
\begin{align}
& E_\nu = E \mp \Frac{v_e}{\mu_\nu} 
\pm \dx{} \left(\dfrac{\kb T}{q}\right).
\label{eq:ET_generalized_E}
\end{align}
The effective electric drift $E_\nu$ acting on each ion species
can be interpreted as the linear superposition of the electric, fluid 
and thermal forces that drive ion motion in the electrolyte channel fluid.
Relation~\eqref{eq:ET_Mom_DD} can be interpreted as the formal limit 
of~\eqref{eq:THD_Mom_1D} as the relaxation time $\tau_{p_\nu} \rightarrow 0$.
\end{remark}
\end{subequations}

\begin{proposition}[Consistency with respect to thermal equilibrium]
The vET satisfies relations~\eqref{eq:thermal_eq}.
\end{proposition}
\begin{proof}
The proof is similar to that for the vTHD model and therefore it is omitted.
\end{proof}

\section{Solution map for the vTHD system}\label{sec:solution_map}
In this section we describe the solution map that is used 
to iteratively solve the vTHD and vET models. The description of the
method is discussed in detail in the case of the vTHD system as 
a similar approach is used to treat also the vET model.
The adopted solution map is an extension of 
the decoupled algorithm known as Gummel's map, a functional tool
widely employed in contemporary Drift-Diffusion simulation of semiconductor devices~\cite{JeromeBook,selberherr,Markowich}. 
The Gummel solution map is basically a nonlinear block 
Gauss-Seidel iteration that allows one to subdivide 
the considered PDE model system into blocks of single equations 
to be solved separately and in sequence until convergence is achieved.

To describe the method, we 
follow the idea proposed in~\cite{Jacoboni1993} (cf. Eq. (25) of this
latter reference), and we introduce the generalization of the 
Maxwell-Boltzmann statistics~\eqref{eq:MB_relations} 
to the case of the vTHD system:
\begin{subequations}\label{eq:MB_relations_thd}
\begin{align}
& n = \cref \exp\left( q \Frac{\pot - \varphi_n}{\kb T_n} \right)
& \\ 
& p = \cref \exp\left( q \Frac{\varphi_p - \pot}{\kb T_p} \right)
& 
\end{align}
\end{subequations}
where, for notational brevity, we have omitted the superscript 
$ec$ in the electrochemical potentials for anions and cations.
The main difference between the 
definitions~\eqref{eq:MB_relations_thd}
and the corresponding~\eqref{eq:MB_relations} valid in the case of the
standard PNP model, is that the ion temperatures $T_n$ and $T_p$ are
used instead of the system (constant) temperature $T_{sys}$.
Assuming local thermal equilibrium at the channel entrance and outlet, 
the temperature of the ions and of the fluid 
satisfy the following Dirichlet boundary conditions:
\begin{subequations}\label{eq:BCT_thd}
\begin{align}
& T_n^{0}(0) = T_p^{0}(0) = T_e^{0}(0) = T_L &  \\
& T_n^{0}(d) = T_p^{0}(d) = T_e^{0}(d) = T_R. & 
\end{align}
\end{subequations}
We notice that $T_L$ and $T_R$ are the 
externally accessible temperatures of the intra and extracellular 
baths, not necessarily being equal to the same value.
This feature of the model is important when an external temperature
gradient has to be enforced across the channel to model 
the heat biosensors described in Sect.~\ref{sec:introduction}. 
Conditions~\eqref{eq:BCT_thd} imply that
$\vth = \kb T_L/q$ at $x=0$ and $\vth = \kb T_R/q$ at $x=d$.

The Gummel algorithm for the iterative solution of the 
vTHD model starts with an initial guess 
$\left[\varphi^{0}, \varphi_p^{0}, \varphi_n^{0}, 
T_p^{0}, T_n^{0}, T^{0}_e\right]^T$ that satisfies 
the boundary conditions~\eqref{eq:bc_summary}
and~\eqref{eq:BCT_thd}.
Then, for $k \geq 0$ until convergence, the algorithm consists of the
following steps:
\begin{enumerate}
\item solve with the damped Newton method (see~\cite{selberherr}, 
Chapt. 7) the nonlinear Poisson equation (NLP) resulting
from substituting~\eqref{eq:MB_relations_thd} into~\eqref{eq:POISSON}.
This returns the updated potential $\varphi^{k+1}$;
\item solve the continuity and momentum balance 
equations~\eqref{eq:THD1D_cont}-~\eqref{eq:THD_Mom_1D}. 
This returns the updated concentrations $p^{k+1}$, $n^{k+1}$;
\item solve the energy equations~\eqref{eq:THD1D_en}-~\eqref{eq:THD_s_1D} .
This returns the updated ion temperatures $T_p^{k+1}$, $T_n^{k+1}$;
\item solve the fluid 
energy equation~\eqref{eq:THD1D_en_L}-~\eqref{eq:THD_sL_1D}.
This returns the updated fluid temperature $T_e^{k+1}$;
\item update the electrochemical potentials by 
inverting~\eqref{eq:MB_relations_thd}:
\begin{subequations}
	\begin{align}
	\varphi_n^{k+1} = \varphi^{k+1} - \frac{\kb T_n^{k+1}}{q}\ln\left( \frac{n^{k+1}}{\cref} \right)\\
	\varphi_p^{k+1} = \varphi^{k+1} + \frac{\kb T_p^{k+1}}{q}\ln\left( \frac{p^{k+1}}{\cref} \right);
	\end{align}
\item check the convergence of the iteration 
by controlling whether the maximum absolute difference between two consecutive iterations $k$ and $k+1$ is less than a prescribed tolerance
\begin{equation}\label{eq:toll}
	\max_{\eta \in U} || \eta^{k+1} - \eta^k||_\infty < {\tt toll}
\end{equation}
\end{subequations}
where {\tt toll} is a given tolerance, $U$ is the set of unknowns 
$\left[\varphi, \varphi_p, \varphi_n, T_p,T_n, T_e \right]^T$
and $\| \cdot ||_\infty$ is the $L^\infty$ norm.
If condition \eqref{eq:toll} is satisfied, the algorithm stops.
\end{enumerate}

\null
A deep analysis of the convergence properties and mathematical
foundations of the 
Gummel decoupling algorithm for the stationary DD model of 
semiconductors can be found in~\cite{JeromeBook}. 
Preconditioning strategies for improving 
the convergence rate of the map are investigated in~\cite{Saad}.

\section{The Advection-Diffusion-Reaction Model Problem}\label{sec:mod_prob}
As outlined in items (e)-(h) of the introduction, we are providing a detailed description of the algorithm, beginning in this section, and continuing through Section 9.

All the linearized equations of the vTHD system~\eqref{eq:THD1D} that need 
be solved during the iterative procedure described in 
Sect.~\ref{sec:solution_map} can be cast into the following general 
Boundary Value model Problem (BVP): \\

find $u$ and $J$ such that:
\begin{subequations}
\label{eq:BVP_ch4}
\begin{align}
		&\dx{J} + c u = g & &\text{in } \Omega \label{eq:BVP_cons}\\
		&J = v u - D \dx{u} & &\text{in } \Omega \label{eq:BVP_cost}\\
		&u= \overline{u} & &\text{on } \po.
\end{align}
\end{subequations}
In~\eqref{eq:BVP_ch4}, $\Omega$ is the open interval $(0, d)$, 
$u$ and $J$ denote the primal unknown and the
associated flux, while $c$ and $g$ are the reaction and production terms
such that~\eqref{eq:BVP_cons} can be regarded as a stationary 
conservation law for $u$. In view of the numerical treatment of~\eqref{eq:BVP_ch4} we assume $D$, $v$, $c$ and $g$ to be 
piecewise smooth functions in $\Omega$, 
the diffusion coefficient $D$ being a positive bounded 
function while $c$ and $g$ are nonnegative given functions. 
Dirichlet boundary conditions are expressed by the 
nonnegative function $\overline{u}: \partial \Omega \rightarrow \R^+$.
Let
\begin{subequations}\label{eq:cont_max}
\begin{equation}
\mathcal{L}u := \dx{J} + cu = \dx{} \left( -D \dx{u} + vu \right) + cu : 
\Omega \rightarrow \R.
\end{equation}
The following properties express important biophysical 
features of the solution of the linearized BVP~\eqref{eq:BVP_ch4}.
We refer the reader to~\cite{RoosStynesTobiska1996} for details and examples.
\begin{definition}[Inverse monotonicity]
Let $w \in C^2(\Omega) \cap C^0(\widebar{\Omega})$. 
We say that $\mathcal{L}$ is inverse-monotone if the inequalities
\begin{align}
	& \mathcal{L}w(x) \geq 0 && \forall\,x\in \Omega\\
	& w(x) \geq 0 && \forall\, x \in \po
\end{align}
together imply that
\begin{equation}\label{eq:w_geq_0}
	w(x) \geq 0 \qquad \forall\, x \in \widebar{\Omega}.
\end{equation}
\end{definition}

\begin{remark}
Inverse monotonicity expresses the fact that the dependent variable of the problem, say, a concentration, a temperature or a mass density, cannot take negative values.
\end{remark}

\begin{theorem}[Comparison principle]
Suppose that there exists a function $\phi \in C^2(\Omega) \cap C^0(\widebar{\Omega})$ such that
\begin{align}
	\mathcal{L}u(x) \leq \mathcal{L}\phi(x) &&& \forall\,x\in \Omega\\
	u(x) \leq \phi(x) &&& \forall\,x\in \po.
\end{align}
Then, we have
\begin{equation}\label{eq:u_leq_phi}
	u(x) \leq \phi(x) \qquad \forall\, x \in \Omega
\end{equation}
and we say that $\phi$ is a \emph{barrier function} for $u$.
\end{theorem}
Combining \eqref{eq:w_geq_0} and \eqref{eq:u_leq_phi}, we obtain the following result which is a very useful tool in the approximation process of the BVP \eqref{eq:BVP_ch4}.
\begin{theorem}[Continuous Maximum Principle]
Suppose that $\mathcal{L}$ is inverse-monotone and
that the comparison principle holds for a suitable barrier function $\phi$. Then, setting $M_\phi : = \displaystyle 
\max_{x\in \Omega} \phi(x)$, we have
\begin{equation}\label{eq:0_leq_u_leq_M}
0 \leq u(x) \leq M_\phi \qquad \forall\, x \in \widebar{\Omega}
\end{equation}
and we say that $u$ satisfies a \emph{continuous maximum principle (CMP)}.
\end{theorem}
\end{subequations}

\section{Finite Element Approximation}\label{sec:FE_discretization}
The content of the present methodological section is organized 
as follows. In Sect.~\ref{sec:DM} we introduce 
the dual-mixed weak formulation of~\eqref{eq:BVP_ch4} 
and in Sect.~\ref{sec:FEM} its corresponding hybridized 
finite element approximation, denoted DMH method.
Then, in Sections~\ref{sec:static_cond} and~\ref{sec:DMH_alg_sys}
we illustrate how to reduce the computational complexity of the 
dual-mixed hybridized method
by the use of static condensation, which leads to solving
a linear algebraic system of the same structure as a standard nodal-based
finite element scheme. 

\subsection{Continuous formulation}
\label{sec:DM}

In the dual mixed method, the variables $u$ and $J$ in~\eqref{eq:BVP_ch4} are treated as independent unknowns, each of which belongs to a different 
functional space, namely $u \in V:=L^2(\Omega)$ and $ J \in Q:=H(\text{div},\Omega) \equiv H^1(\Omega)$.
We proceed by formally multiplying equation \eqref{eq:BVP_cons} by a test function $\Vtf \in V$ and equation \eqref{eq:BVP_cost} by a test function $\Qtf \in Q$, and integrate by parts over $\Omega$ to obtain the problem:\\

find $(u,J) \in V \times Q$ such that, for all $(\Vtf,\Qtf) \in V \times Q$:
\begin{subequations}\label{eq:DM_cont}
\begin{empheq}[left=\empheqlbrace]{align}
		&\int_\Omega D^{-1}\left( J -vu \right) \Qtf \DX
		-\int_\Omega u \dx{\Qtf} \DX= - \int_{\po} (\Qtf \cdot \vec{n}) 
		\overline{u} \DS\\
		&-\int_\Omega \dx{J} \Vtf \DX- \int_\Omega cu \Vtf 
		\DX= -\int_\Omega g\Vtf \DX.
\end{empheq}
\end{subequations}
Existence and uniqueness of the solution pair $(u,J)$ of problem \eqref{eq:DM_cont} can be proved under suitable coercivity assumptions, see \cite{DouglasRoberts}.

\subsection{Discrete formulation}
\label{sec:FEM}
We denote by $\Omega$ the open interval $(0,L)$
and by $\po=\{0,L\}$ its boundary on which an outward unit normal $\vec{n}$ is defined,
in such a way that $\vec{n}(0) = -1$ and $\vec{n}(L) = +1$ 
We introduce a partition (triangulation) 
$\Th$ of $\overline{\Omega}$ into a number $\numel	 \geq 2$ of intervals $K_i = [x_i, x_{i+1}], i=1,\ldots,\numel,$ with $x_1 = 0$ and $x_{\numel+1} = L$. 
On each element $K_i$ we introduce a local unit outward normal vector $\vec{n}|_{\pk_i}$, we denote by $h_i := x_{i+1} - x_i$ the length 
of the interval and set 
$h := \max_{\Th} h_i$.  We also set $N := \numel + 1$ and 
denote with $\mathcal{N}_h := \lbrace x_i\rbrace_{i=1}^N$ 
the set of vertices of $\Th$.
We associate with $\Th$ the following function spaces:
\begin{subequations}
\begin{align}
	& V_h := \lbrace \Vtf_h \in L^2(\Omega) \text{ s.t. } \Vtf_h|_K \in \Poly_0(K)\, \forall K \in \Th \rbrace\\
	& Q_h{} := \lbrace \Qtf_h \in L^2(\Omega) \text{ s.t. } \Qtf_h|_K \in \Poly_1(K)\, \forall K \in \Th \rbrace\\
	& \Lambda_{h,\rho} := \lbrace \lambda_h : \mathcal{N}_h \rightarrow \R \text{ s.t. } \lambda_h(a)=\rho_a, \lambda_h(b)=\rho_b \rbrace
\end{align}
\end{subequations}
where $V_h $ is the vector space of piecewise constant polynomials defined over $\Omega$, $Q_h{}$ is the vector space of piecewise linear polynomials discontinuous over $\overline{\Omega}$, and $\Lambda_{h,\rho}$ is the space of functions defined only at the vertices of $\Th$ with given boundary values.
Notice that $\dim(Q_h{})=2 \numel$ because its functions are not 
continuous at interelement interfaces.
The DMH formulation of problem \eqref{eq:DM_cont} is:\\

find $(J_h{},u_h{},\lambda_h) \in Q_h{} \times V_h \times 
\Lambda_{h,\overline{u}}$ 
such that, for all $(\Qtf_h, \Vtf_h, \Ltf_h) \in Q_h{} \times V_h \times \Lambda_{h,0}$:
\begin{subequations}
\label{eq:HYB_sys_FEM}
\begin{align}
		&\sum_{K\in \Th}\left[ \int_K 
		D_h^{-1}\left( J_h{} -vu_h{} \right) \Qtf_h \DX
		-\int_K u_h{} \dx{\Qtf_h} \DX + \int_{\pk} (\Qtf_h \cdot \vec{n}|_{\pk})\lambda_h \DX \right]=0 \label{eq:DMH_1}\\
		&-\sum_{K\in \Th}\left[\int_K \dx{J_h{}} \Vtf_h \DX- \int_K cu_h{} \Vtf_h \DX \right]= -\sum_{K\in \Th}\int_K g\Vtf_h \DX \label{eq:DMH_2}\\[1em]
		& \sum_{K\in \Th}\int_{\pk} (J_h{} \cdot \vec{n}|_{\pk}) \Ltf_h \DX = 0 \label{eq:HY_cont_cond}
\end{align}
\end{subequations}
where $D_h : \Th \rightarrow \R^+$ is a modified diffusion coefficient containing a stabilization term to deal with advective-dominated
problems. 
\begin{remark}
The proof of existence and uniqueness of the solution 
$(J_h{},u_h{},\lambda_h)$ of~\eqref{eq:HYB_sys_FEM} is not an easy 
task because of the advective term.
In what follows, we assume that such a solution always exists and 
we refer to~\cite{BrezziFortin1991,RobertsThomas1991} for the 
general theoretical tools required in the analysis of~\eqref{eq:HYB_sys_FEM}
and to~\cite{ArbogastChen1995} for an example of such analysis in the case of 
mixed methods interpreted as nonconforming discretizations.
\end{remark}

\begin{remark}
Conditions \eqref{eq:HY_cont_cond} imply the continuity of $J_h{}$ across interdomain boundaries.
\end{remark}

\begin{remark}
The \emph{hybrid variable} $\lambda_h$ is the Lagrange multiplier associated with the continuity constraint of $J_h{}\cdot\vec{n}$ and it can be regarded as an approximation of $u$ at the grid nodes:
\begin{equation*}
	\lambda_h \simeq u|_{\mathcal{N}_h}.
\end{equation*}
\end{remark}

\subsection{Static condensation}
\label{sec:static_cond}
Consider a generic element $K \in \Th$. From~\eqref{eq:DMH_1}, 
we can compute $J_h|_K$ by inverting the local flux mass matrix as
\begin{subequations}
\begin{equation}\label{eq:Jk}
\vec{J}_K = -A_K^{-1} \left[ (B_K^T + C_K)u_K + L_K \lambdaB_K\right] 
\in \R^{2 \times 1}.
\end{equation}
Substituting $\vec{J}_K$ into~\eqref{eq:DMH_2}, we get $u_h|K$
\begin{equation}\label{eq:uk}
u_K = H_K^{-1}	\left[ g_K + B_K A_K^{-1} L_K \lambdaB_K\right]
\in \R^{1 \times 1}
\end{equation}
where
\begin{equation*}
H_K := \left[-B_K A^{-1}_K (B_K^T + C_K) + E_K\right] \in \R^{1 \times 1}.
\end{equation*}
Now, using \eqref{eq:uk} in \eqref{eq:Jk}, 
we can express $\vec{J}_K$ in terms of the sole 
hybrid variable $\lambdaB_K$ and of $g_K$ as
\begin{equation}\label{eq:DMH_JK}
	\vec{J}_K = M_K \lambdaB_K + \vec{b}_K
\end{equation}
where
\begin{equation*}
	M_K := -\left[ A_K^{-1}(B_K^T + C_K) H_K^{-1} B_K A_K^{-1} L_K^T + A_K^{-1} L_K \right]
	\in \R^{2 \times 2}
\end{equation*}
and
\begin{equation*}
	\vec{b}_K := \left[-A_K^{-1} (B_K^T + C_K) H_K^{-1} g_K\right]
	\in \R^{2 \times 1}.
\end{equation*}
\end{subequations}
The above procedure is well known as \emph{static condensation} and 
corresponds to Gaussian elimination of the internal variables 
$J_h|_K$ and $u_h|_K$ in favor of $\lambda_h|_{\partial K}$ 
at the level of the element $K\in \Th$ (see~\cite{ArnoldBrezzi}).

\subsection{The algebraic system}\label{sec:DMH_alg_sys}
Enforcing the continuity of $J_h{}$ at each internal node through~\eqref{eq:HY_cont_cond} yields the following $\numel-1$ conditions
\begin{subequations}
\begin{equation}\label{eq:DMH_J_cont}
J_2^{K_{i-1}}(\lambda_{i-1},\lambda_i) = 
J_1^{K_{i}}(\lambda_{i},\lambda_{i+1})
\quad i=1,\ldots,M-1
\end{equation}
which can be written in the form of a linear \emph{tridiagonal} 
algebraic system for $\lambdaB$, i.e., for $i=1,\ldots,M-1$:
\begin{subequations}\label{eq:DMH_tridiagonal}
\begin{empheq}[left=\empheqlbrace]{align}
	& M_{21}^{K_i} \lambda_{i-1} + (M_{22}^{K_i} - M_{11}^{K_{i+1}}) \lambda_{i} - M_{12}^{K_{i+1}} \lambda_{i+1}
	= b_1^{K_{i+1}} - b_2^{K_i}\\
	& \lambda_1=\overline{u_0},\quad \lambda_N=\overline{u_L}.
\end{empheq}
\end{subequations}
Equations~\eqref{eq:DMH_tridiagonal} constitute a linear 
system of the form
\begin{equation}\label{eq:DMH_sys_lambda}
\mathcal{M} \lambdaB = \vec{f}{}
\end{equation}
where the entries of the stiffness 
matrix $\mathcal{M} \in \R^{M-1 \times M-1}$ are defined as:
\begin{subequations}\label{eq:DMH_matrix_M}
\begin{empheq}[left=\empheqlbrace]{align}
	& \mathcal{M}_{i,i-1} = M_{21}^{K_i}\\
	& \mathcal{M}_{i,i}   = M_{22}^{K_i} - M_{11}^{K_{i+1}}\\
	& \mathcal{M}_{i,i+1} = -M_{12}^{K_{i+1}}
\end{empheq}
\end{subequations}
and the entries of the load vector $\vec{f}{} \in \R^{M-1}$ 
are $f{}_i = b_1^{K_{i+1}} - b_2^{K_i}$.
The explicit expressions of these entries will be specified in
Sect.~\ref{sec:matrixentries}.
\end{subequations}

\section{The Stabilized DMH Method and the Discrete Maximum Principle}\label{sec:DMP}

If the solution of a given boundary value problem satisfies a CMP, then a properly designed approximation should behave in the same way. A numerical scheme that does not generate spurious global extrema in the interior of the computational domain is said to satisfy a \emph{discrete maximum principle} (DMP). We assume that the solution $u$ of the BVP \eqref{eq:BVP_ch4} satisfies the a priori estimate \eqref{eq:0_leq_u_leq_M} and we characterize the conditions under which the same estimate is satisfied also by the DMH approximation $\lambda_h$ computed by solving the linear algebraic system \eqref{eq:DMH_sys_lambda} associated with the problem.

\setcounter{equation}{24}  
The following definitions turn out to be useful.
\begin{subequations}\label{eq:monotone_h}
\begin{definition}[Inverse monotone matrix]
An invertible square matrix $A$ of size $n$ is said to be inverse monotone if
\begin{equation}\label{eq:inv_mon_matrix}
A^{-1} \geq 0
\end{equation}
the inequality being understood in the element-wise sense.
\end{definition}

A special class of monotone matrices is that introduced below.
\begin{definition}[M-matrix]
An invertible square matrix $A$ is an M-matrix if:
\begin{itemize}
\item $A_{ij} \leq 0$ for $i \neq j$;
\item A is inverse-monotone.
\end{itemize}
\end{definition}

Given a function $\eta_h \in \Lambda_{h,\rho}$, for any given function
$\rho: \partial \Omega \rightarrow \R$, we denote by
$\rho^\ast_h$ the piecewise linear continuous interpolate of $\eta_h$
over $\Th$.
\begin{theorem}[Sufficient condition for DMP]\label{teo:suff_DMP}
Assume that matrix $\mathcal{M}$ is an M-matrix. Then $\lambda_h^\ast$ 
satisfies the DMP, i.e.
\begin{equation}\label{eq:discrete_MP}
0 \leq \lambda_h^\ast(x) \leq M_\phi \qquad \forall \, x \in \widebar{\Omega}.
\end{equation}
\end{theorem}
The following (necessary and sufficient) condition is useful to verify the property of being an M-matrix.

\begin{theorem}[Discrete comparison principle]
\label{teo:Mmatrix}
Let $A$ be an invertible matrix with non-positive off diagonal entries
($A_{ij} \leq 0 \text{ for } i \neq j$). Then, $A$ is an M-matrix
if and only if there exists a positive vector $\vec{e}$ such that $A\vec{e} \geq 0$ (in the component-wise sense),
with at least one row index $i^*$ such that $(A\vec{e})_{i^*} > 0$.
\end{theorem}
\end{subequations}


\subsection{The stabilized DMH method}
\label{sec:matrixentries}
The solution of the model BVP~\eqref{eq:BVP_ch4} may exhibit sharp
boundary and/or internal layers according to the relative weight of
the advective and reactive coefficients with respect to the diffusion
term (see~\cite{RoosStynesTobiska1996}). Correspondingly, this is 
well known to reflect into numerical instability (spurious unphysical
oscillations) if the local mesh size $h$ is not sufficiently small.
However, a smaller value of $h$ means a larger size of the linear algebraic system and thus an increased computational effort.
To avoid this inconvenience, we introduce in the DMH 
formulation~\eqref{eq:HYB_sys_FEM} two specific approaches 
that ensure numerical stability without necessarily resorting 
to a small value of $h$:
\begin{description}
\item[(S1)] diagonal \emph{lumping} of the local mass flux matrix $A_K$ 
(for dominant reaction);
\item[(S2)] addition of an \emph{artificial diffusion} (for dominant convection).
\end{description}
In particular, in the remainder of the section 
we show that, when (S1) and (S2) are used in conjunction, 
the stiffness matrix $\mathcal{M}$ is an M-matrix irrespective 
of the value of $h$. 
This property guarantees that system~\eqref{eq:DMH_sys_lambda} is uniquely
solvable and that the piecewise linear interpolate of the solution, $\lambda_h^\ast$, satisfies the DMP and,
in particular, the a-priori estimate~\eqref{eq:discrete_MP}.
We refer to~\cite{RoosStynesTobiska1996} 
for a thorough discussion of continuous and 
discrete maximum principles enjoyed by the solution of singularly
perturbed BVPs and by their corresponding finite element and finite
difference approximations. For ease of presentation, 
we assume henceforth that problem coefficients $D,v,c$ and $g$ are
constant positive quantities 
and that the mesh size is uniform and given by $h=1/\numel$.

\subsubsection{(S1): Lumping of the local flux mass matrix}
The stabilization approach to deal with the case where 
the reaction coefficient $c$ dominates over the diffusion coefficient
$D$ in~\eqref{eq:BVP_ch4} consists of replacing the local flux mass matrix
$A_K$ in~\eqref{eq:Jk} with a diagonal matrix $\widehat{A}_K$
such that:
\begin{subequations}\label{eq:lumping}
\begin{align}
& (\widehat{A}_K)_{ii} = \sum_{j=1}^2 (A_K)_{ij} = \dfrac{1}{2} & \qquad
i=1,2 \label{eq:AK_ii_lumped} \\
& (\widehat{A}_K)_{ij} = 0 & \qquad i \neq j. \label{eq:AK_ij_lumped}
\end{align}
The approach~\eqref{eq:lumping} is referred to as \emph{mass lumping} 
and is equivalent to using the trapezoidal quadrature rule to compute the integral
$$
(A_K)_{ij} = \int_K \Qbf_j \Qbf_i \DX \qquad i,j=1,2
$$
where
$$
\Poly_1(K) = {\rm span}\left\{ \Qbf_1, \Qbf_2 \right\} \qquad
\forall K \in \Th.
$$
Using~\eqref{eq:lumping} in~\eqref{eq:DMH_1}, 
when $\Qtf_h = \Qbf_1$, we get
\begin{equation}\label{eq:J_1}
J_1 = \vtk u_K - \dtk_h \dfrac{u_K - \lambda_1}{h/2},
\end{equation}
while, when $\Qtf_h = \Qbf_2$, we similarly get
\begin{equation}\label{eq:J_2}
J_2 = \vtk u_K - \dtk_h \dfrac{\lambda_2 -u_K}{h/2}.
\end{equation}
Equations~\eqref{eq:J_1} and~\eqref{eq:J_1} express each 
local degree of freedom of the flux $J_h|_K$ 
as the sum of an \emph{average} advective flux and a finite difference 
approximation of the Fick diffusive flux flowing between the boundary of
$K$ and its center. The two relations are the DMH analogues
of those of the dual-mixed stabilized method of~\cite{SaccoSaleri97}.
\end{subequations}

\subsubsection{(S2): Artificial diffusion}
The stabilization approach to deal with the case where 
the advection coefficient $v$ dominates over the diffusion coefficient
$D$ in~\eqref{eq:BVP_ch4} consists of replacing the diffusion coefficient $D$ with the following modified expression (cf.~\cite{Brooks1982})
\begin{subequations}
\begin{equation}\label{eq:art_diffusion}
D_h:= 
D\left( 1 + \Phi(\mathcal{P}e_{loc}) \right)
\end{equation}
where $\mathcal{P}e_{loc} := \dfrac{|v|h}{2D}$ is the local P{\'e}clet number and $\Phi : \R^+ \rightarrow \R^+$ is a suitable stabilization function to be chosen in such a way that:
\begin{gather}
	\Phi(t) \geq 0\\
	\lim_{t \rightarrow 0^+} \Phi(t) = 0.
\end{gather}
The approach~\eqref{eq:art_diffusion} is referred to as \emph{artificial diffusion}. Two special choices of $\Phi$ are:
\begin{itemize}
\item Upwind (UP) stabilization function
\begin{equation}\label{eq:UP_peclet}
	\Phi^{UP} (t) := t
\end{equation}
\item Scharfetter-Gummel (SG) stabilization function
\begin{equation}
\label{eq:SG_peclet}
\Phi^{SG} (t) := t-1 + \mathcal{B}(2t), \qquad 
\mathcal{B}(t):= \dfrac{t}{e^t -1 }.
\end{equation}
In the case where~\eqref{eq:SG_peclet} is adopted, the resulting stabilized DMH coincides with the classic 
SG exponentially fitted difference scheme (cf.~\cite{SG}).
\end{itemize}
The accuracy of the two stabilized DMH methods, as $h \rightarrow 0$,
is $\mathcal{O}(h)$ if $\Phi=\Phi^{UP}$ while is $\mathcal{O}(h^2)$ if $\Phi=\Phi^{SG}$. Thus, the SG formulation is far superior than the upwind
method in the limit of a small grid size.
However, for a large value of the P\`eclet number
the two stabilized methods are practically equivalent in terms of the 
amount of added artificial diffusion (see also~\cite{SaccoSaleri97}
and~\cite{manganini2013}).
\end{subequations}

\subsubsection{Matrix form of the stabilized DMH method}

Using superscripts $-$ and $+$ 
to refer to quantities on the interval $K^-_i = [x_{i-1}, x_i]$ and $K^+_i = [x_{i}, x_{i+1}]$, respectively, let us enforce continuity of 
$J_h$ at every interior node $x_i$ $(i=2,\ldots, N-1)$
\begin{subequations}\label{eq:matrix_coeff}
\begin{equation}\label{eq:DMH_cont_stab}
J_2^- = J_1^+. 
\end{equation}
Then, combining equations~\eqref{eq:J_1} and~\eqref{eq:J_2} 
with condition \eqref{eq:DMH_cont_stab}, we get
\begin{multline}\label{eq:tridiagonal}
-\lambda_i\dfrac{2 \dtk^-}{h} + 
\dfrac{(\vtk^- + \frac{2\dtk^-}{h}) 
( \gtk^- h^2 + 2 \lambda_{i-1} \dtk^- 
+ 2 \lambda_i \dtk^- )}{4 \dtk^- + \ctk^- h^2} \\
= \lambda_i\dfrac{2 \dtk^+}{h} + 
\dfrac{(\vtk^+ - \frac{2\dtk^+}{h}) 
( \gtk^+ h^2 + 2 \lambda_i \dtk^+ 
+ 2 \lambda_{i+1} \dtk^+ )}{4 \dtk^+ + \ctk^+ h^2}.
\end{multline}
By inspection on~\eqref{eq:tridiagonal} the entries of the stiffness matrix 
$\mathcal{M}$ and of the load vector $\vec{f}{}$ become:
\begin{subequations}
\label{eq:stab_M_matrix}
\begin{align}
	& \mathcal{M}_{i,i-1} = -\left( \dfrac{v}{2} + 
	\dfrac{D (1 + \Phi)}{h} \right) \frac{1}{1 + \dfrac{c h^2}{4 D}} \\
	& \mathcal{M}_{i,i}   = \dfrac{2 D(1 + \Phi)}{h} 
	\left(\dfrac{1 + \dfrac{c h^2}{2 D}}{1 + \dfrac{c h^2}{4D}}\right)\\[3mm]
	& \mathcal{M}_{i,i+1} = \left( \dfrac{v}{2} - 
	\dfrac{D (1 + \Phi)}{h} \right) \frac{1}{1 + \dfrac{c h^2}{4 D}}\\
	& f_i = \dfrac{g h}{1 + \dfrac{c h^2}{4 D}}.
\end{align}
\end{subequations}

\begin{proposition}[Discrete Maximum Principle in the non stabilized case]
Set $\Phi=0$ and assume that
\begin{align}
& \mathcal{P}e_{loc} < 1. & \label{eq:pecletlessthan1}
\end{align}
Then, the stiffness matrix $\mathcal{M}$ is an M-matrix. 
\end{proposition}
\begin{proof}
Recalling that $v>0$, it turns out that $\mathcal{M}_{i,i-1}$ is $< 0$. Moreover, we have $\mathcal{M}_{i,i} > 0$ and $f_i > 0$.
By inspection on $\mathcal{M}_{i,i+1}$ we see that 
if~\eqref{eq:pecletlessthan1} is satisfied 
then $\mathcal{M}$ is an M-matrix and the DMP holds for the solution
of~\eqref{eq:DMH_sys_lambda}. 
\end{proof}
Condition~\eqref{eq:pecletlessthan1} may be 
too restrictive in the choice of the mesh size so that the following
result may be a convenient remedy.
\begin{proposition}[Discrete Maximum Principle in the stabilized case]
When $\Phi = \Phi^{UP}$ or $\Phi = \Phi^{SG}$, the stabilized DMH
method satisfies the DMP irrespective of the value of $\pec_{loc}$.
\end{proposition}
\begin{proof}
It is easy to verify that
\begin{equation*}
\sum_{j=1, \ldots, N} \mathcal{M}_{i,j} > 0
\end{equation*}
for every interior node ($i=2,\ldots, N-1$), and
\begin{align*}
	& \mathcal{M}_{1,1} + \mathcal{M}_{1,2} > 0 \\
	& \mathcal{M}_{N-1,N-1} + \mathcal{M}_{N,N} > 0
\end{align*}
at boundary nodes. Then, the conclusion immediately follows
by applying Theorem~\ref{teo:Mmatrix}.
\end{proof}
\end{subequations}

\section{Experimental Validation of the DMH method}
\label{sec:estimates}
In this section we numerically verify the convergence 
and stability properties of the DMH method.

\subsection{Convergence Analysis}
We consider the model BVP \eqref{eq:BVP_ch4} with with homogeneous BCs, $u(0)=u(L)=0$, $d=5$, coefficients $D=v=c=1$, and $g(x) = e^{-x}\left[ x^2 - (6+L)x + 3L -2 \right]$ in such a way that the exact solution is the pair
\begin{equation*}
u(x) = xe^{-x}(L-x), \qquad J(x) = -e^{-x}\left[  2x^2 +2x(L-1) +L\right].
\end{equation*}
In the numerical approximation of the problem, the mesh size is uniform and
equal to $h = d / \numel $, with $\numel=[10, 20, 40, 80, 160, 320, 640, 1280, 2560]^T$.
Table \ref{tab:error_u} illustrates the convergence history of the 
DMH method by reporting
various error norms for $u-u_h$ and $u-\lambda_h^*$. 
The discrete maximum norm $\norm{\cdot}_{\infty,h}$ associated with the triangulation $\Th$ is defined for each continuous function $\eta(x):[0,L] \rightarrow \R$ as
\begin{equation*}
\norm{\eta(x)}_{\infty,h} := \max_{x_i \in \mathcal{N}_h} | \eta(x_i)|.
\end{equation*}
Table \ref{tab:error_J} reports the error $J-J_h$ in the $L^2$-norm  and in the $H_{\text{div}} \equiv H^1$-norm.
\begin{table}
\centering
\begin{tabular}{r|c|c|c|c}
	$\numel$ & $\norm{u - u_h}_{L^2}$ & $\norm{\Pi_0 u - u_h}_{L^2}$ & $\norm{u - \lambda_h^*}_{L^2}$ & $\norm{u - \lambda_h^*}_{\infty,h}$ \\
	\toprule
	10 & 4.77887e-01 & 2.21777e-01 & 1.30666e-01 & 8.64963e-02 \\
	\hline
	20 & 1.91869e-01 & 6.25684e-02 & 3.38644e-02 & 2.22424e-02 \\
	\hline
	40 & 8.63963e-02 & 1.61248e-02 & 8.54321e-03 & 5.60052e-03 \\
	\hline
	80 & 4.18183e-02 & 4.06198e-03 & 2.14066e-03 & 1.40542e-03 \\
	\hline
	160 & 2.07297e-02 & 1.01743e-03 & 5.35468e-04 & 3.51510e-04 \\
	\hline
	320 & 1.03422e-02 & 2.54479e-04 & 1.33886e-04 & 8.78873e-05 \\
	\hline
	640 & 5.16826e-03 & 6.36272e-05 & 3.34727e-05 & 2.19724e-05 \\
	\hline
	1280 & 2.58377e-03 & 1.59073e-05 & 8.36820e-06 & 5.49311e-06 \\
	\hline
	2560 & 1.29184e-03 & 3.97683e-06 & 2.09202e-06 & 1.37325e-06
\end{tabular}
\caption[Error estimates for $u_h$ and $\lambda_h^*$ as function of \numel.]
{Various error norms for numerical solutions $u_h$ and $\lambda_h^*$ as a function of the number of intervals $\numel$. No stabilization is adopted.}
\label{tab:error_u}
\end{table}

\begin{table}
\centering
\begin{tabular}{r|c|c}
	$\numel$ & $\norm{J - J_h}_{L^2}$ & $\norm{J - J_h}_{H^1}$ \\
	\toprule
10  &     3.28036e-01  &   2.60239e+00 \\ \hline
20   &    8.75295e-02  &   1.30056e+00 \\ \hline
40    &   2.16140e-02  &   6.47536e-01 \\ \hline
80     &  5.29217e-03  &   3.23222e-01 \\ \hline
160     &  1.30365e-03 &    1.61529e-01 \\ \hline
320      & 3.23124e-04 &    8.07533e-02 \\ \hline
640      & 8.04091e-05 &    4.03752e-02 \\ \hline
1280     &  2.00543e-05&     2.01874e-02 \\ \hline 
2560     &  5.00749e-06&     1.00937e-02 \\
\end{tabular}
\caption[Error estimates oft $J_h$ as function of \numel.]{Error norms $H^1$ and $L^2$ for numerical solution $J_h$ as a function of the number of intervals 
$\numel$. No stabilization is adopted.}
\label{tab:error_J}
\end{table}

The analysis of the asymptotic convergence orders that are predicted by Tables \ref{tab:error_u} and \ref{tab:error_J} show that the computed numerical solutions verify the following error estimates
\begin{itemize}
\item $\norm{u - u_h}_{L^2} \leq \mathcal{C} h$
\item $\norm{\Pi_0 u - u_h}_{L^2} \leq \mathcal{C} h^2$
\item $\norm{u - \lambda_h^*}_{L^2}\leq \mathcal{C} h^2$
\item $\norm{u - \lambda_h^*}_{\infty,h}\leq \mathcal{C} h^2$
\item $\norm{J - J_h}_{H^1} \leq \mathcal{C} h$
\item $\norm{J - J_h}_{L^2} \leq \mathcal{C} h^2$
\end{itemize}
The above results are in excellent agreement with the 
theoretical convergence rates predicted in the elliptic case 
in~\cite{ArnoldBrezzi} and~\cite{BrezziFortin1991,RobertsThomas1991}.

\subsection{Reaction-Dominated and Advective-Dominated Regimes}
\label{sec:num_stab}
In this section we demonstrate the efficacy of the 
stabilization techniques proposed in Sect.~\ref{sec:matrixentries}
in the study of two model problems, special instances of the BVP~\eqref{eq:BVP_ch4}. For ease of presentation, we set $L = 1$ and we subdivide the computational domain into $\numel=10$ uniform intervals of size $h = 1/\numel = 0.1$. We
also assume that the coefficients $D, v, c$ and $g$ are constant, with $g = 1$ and we take homogeneous boundary conditions, $u(0)= u(1)=0$ 
($\overline{u}=0$).

\subsubsection{Diffusion-reaction BVP: mass-lumping}
\begin{figure}[h!]
\centering
\subfigure[No lumping]{
\includegraphics[width=0.4\linewidth]{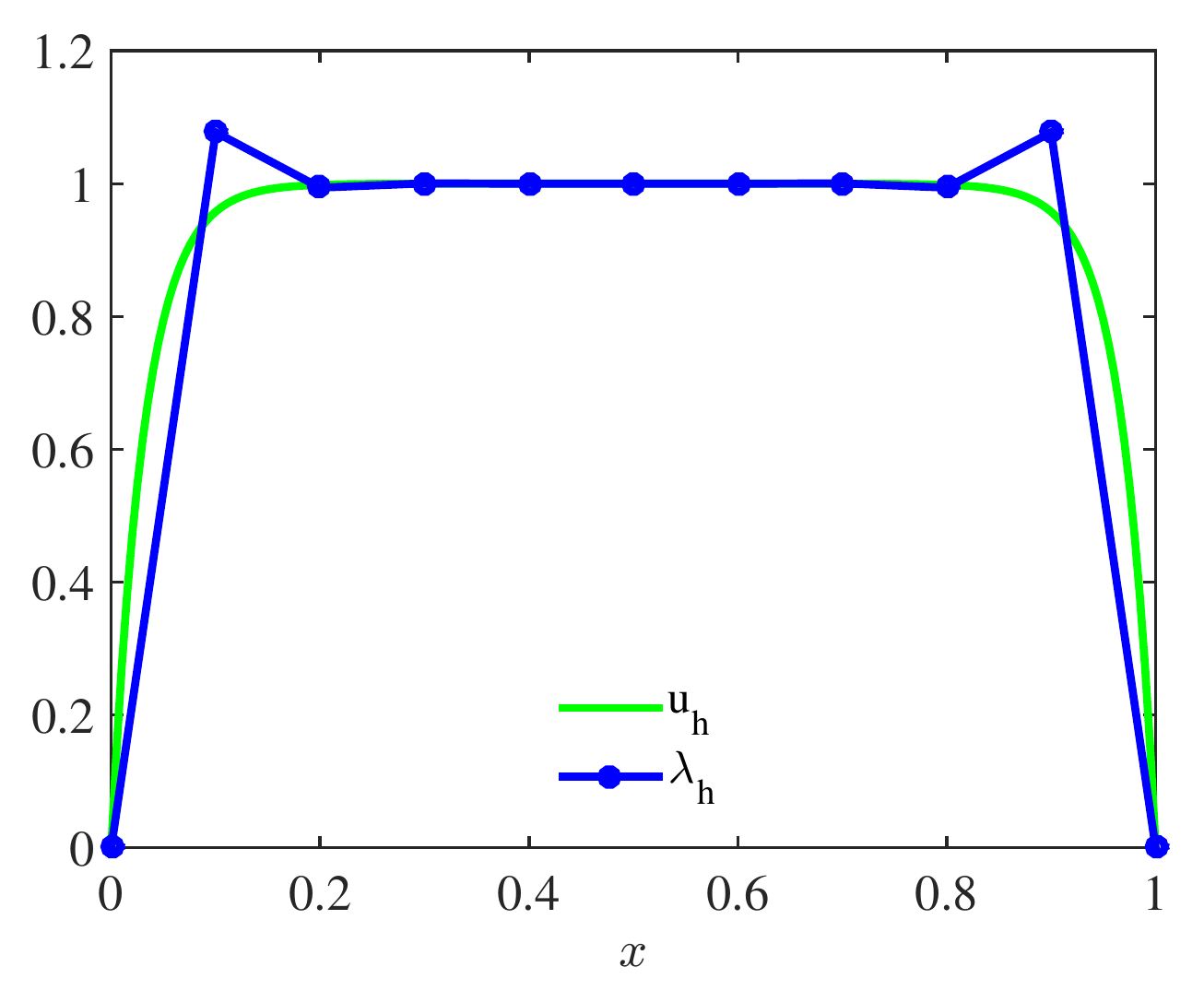}
\label{sfig:no_lump}}
\subfigure[Lumping]{
\includegraphics[width=0.4\linewidth]{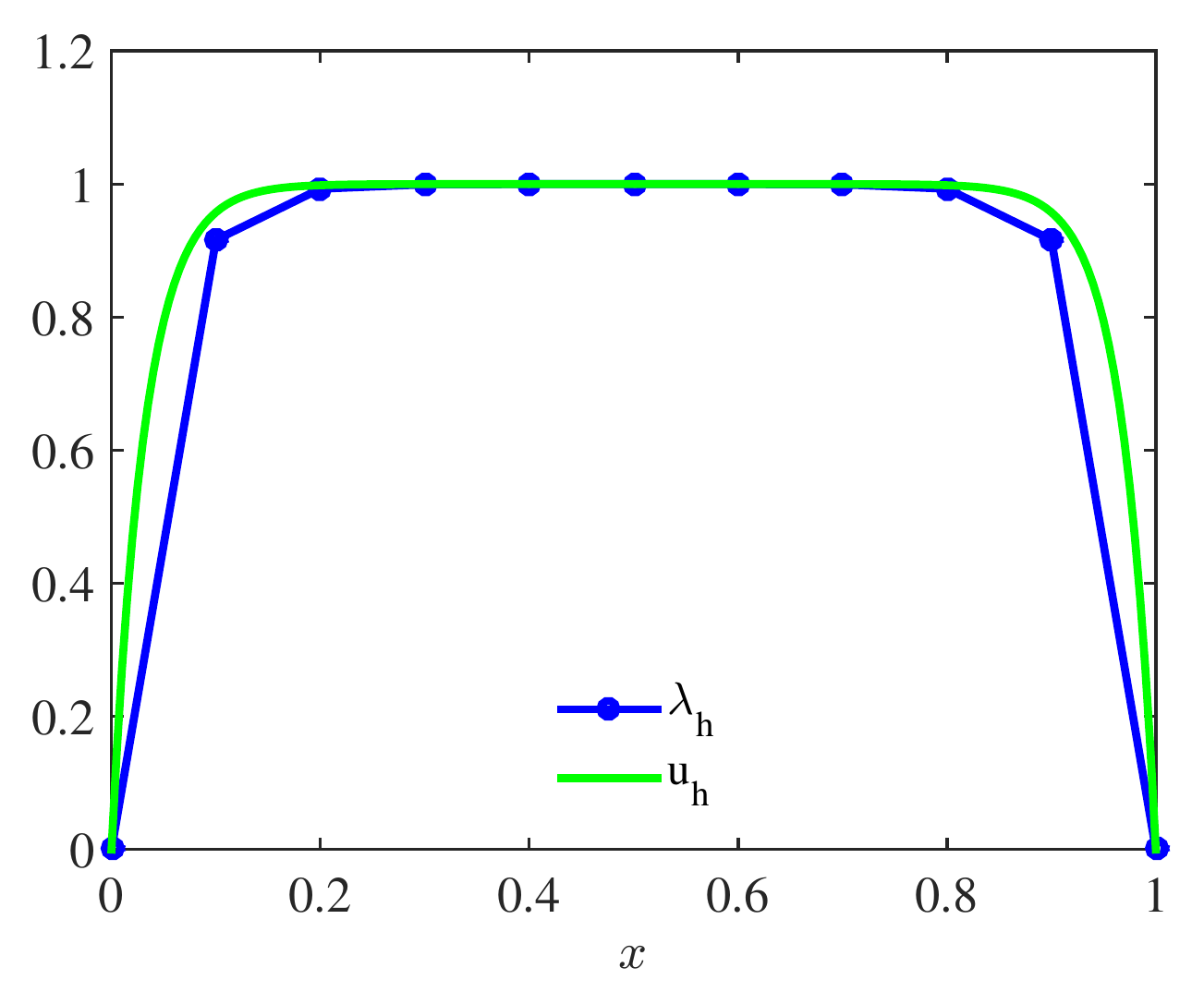}
\label{sfig:lump}}
\caption{Exact and numerical solution of the diffusion-reaction 
BVP. Exact (green, solid line) and 
numerical solutions $\lambda_h{}$ (blue line with bullets).}
\label{fig:diff_reac}
\end{figure}

In this case we have $v=0$. 
A barrier function for $u$ is $\phi(x) = 1/c$.
If the diffusion coefficient is small compared to the reaction term, e.g $D=10^{-3},\, c=1$, spurious (unphysical) oscillations arise near the boundaries, see Figure \ref{sfig:no_lump}. In order to eliminate such oscillations, we adopt the mass-lumping stabilization procedure.
The numerical solution $\lambda_h{}$ obtained with mass-lumping is illustrated in Figure \ref{sfig:lump}.

\subsubsection{Diffusion-advection BVP: artificial diffusion}
\begin{figure}[h!]
\centering
\subfigure[$\Phi = 0$]{
\includegraphics[width=0.3\linewidth]{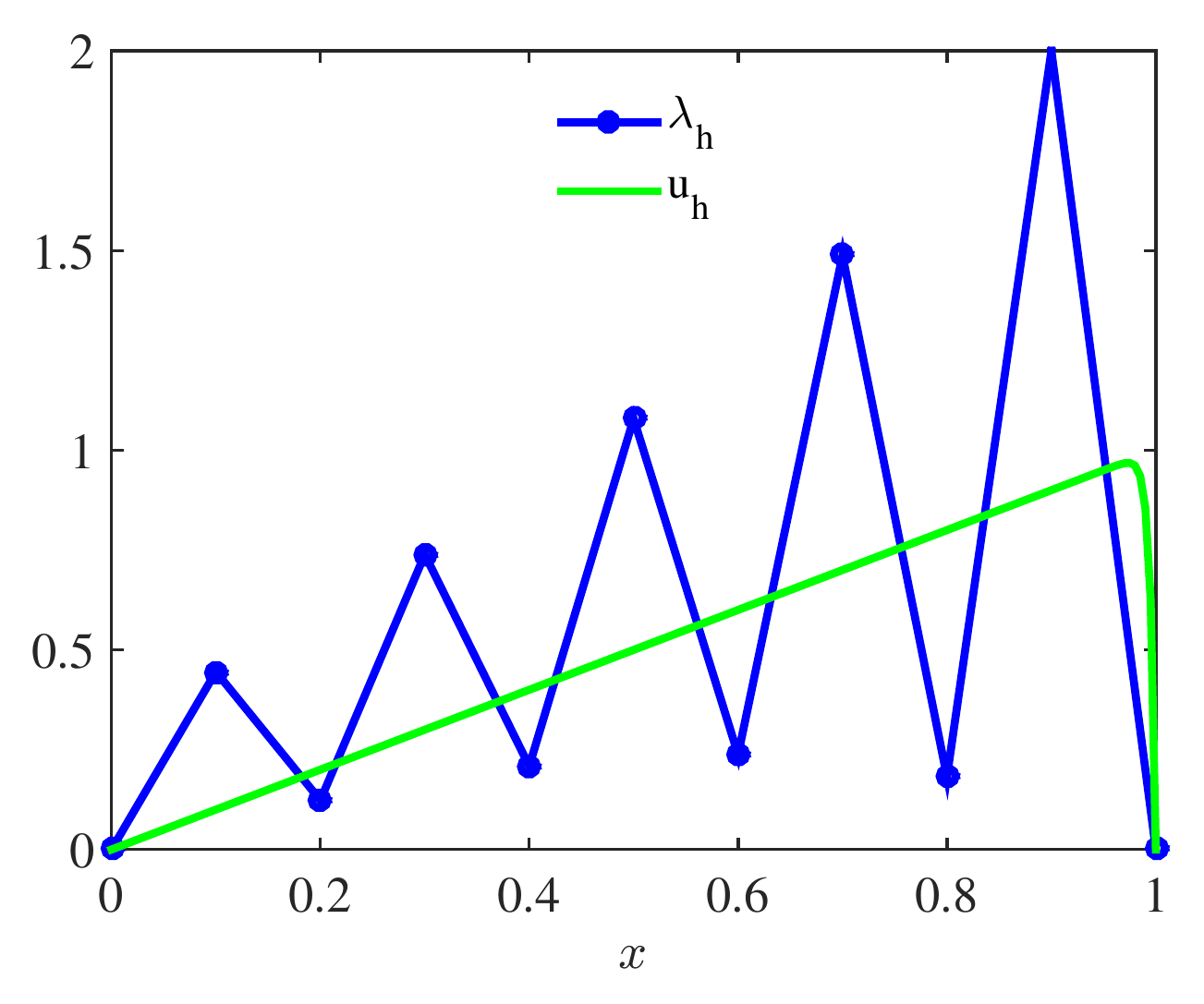}
\label{sfig:no_art_diff}}
\subfigure[$\Phi = \Phi^{UP}$]{
\includegraphics[width=0.3\textwidth]{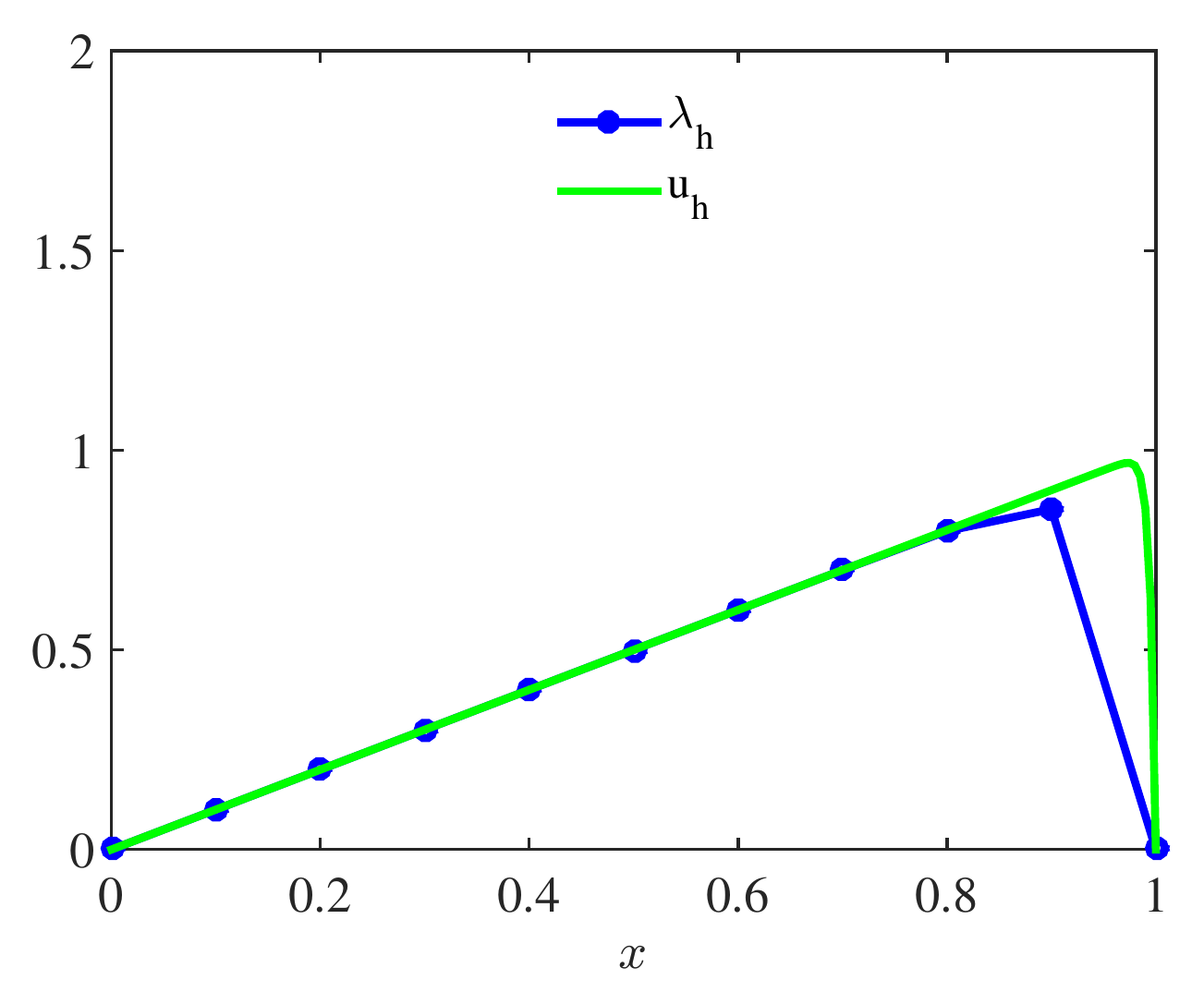}
\label{sfig:UP}}
\subfigure[$\Phi = \Phi^{SG}$]{
\includegraphics[width=0.3\textwidth]{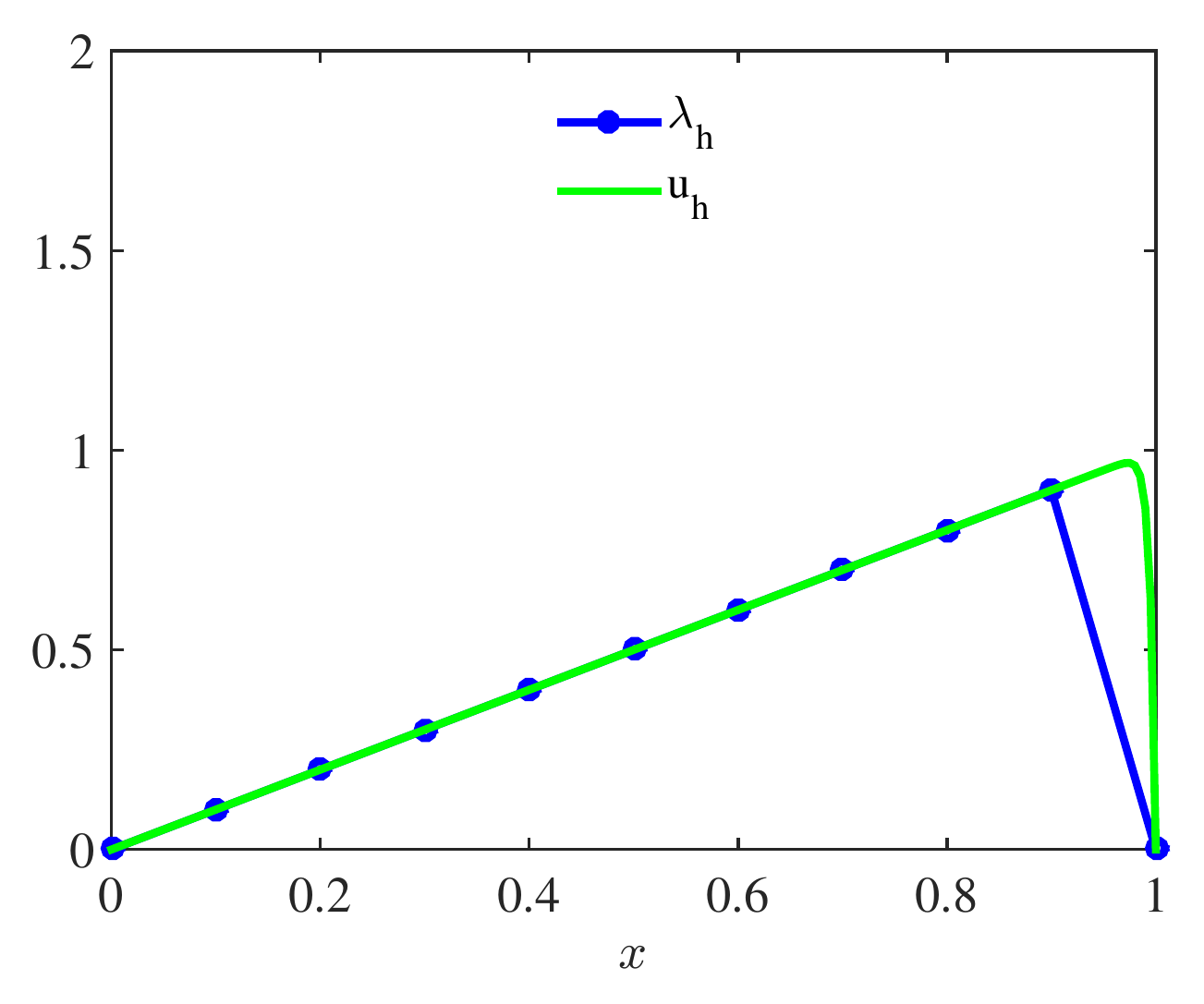}
\label{sfig:SG}}
\caption{Exact and numerical solution of the diffusion-advection 
BVP. Exact (green, solid line) and numerical 
solutions $\lambda_h{}$ (blue line with bullets).}
\label{fig:diff_conv_up_sg}
\end{figure}

In this case we have $c=0$. 
A barrier function for $u$ is $\phi(x) = x/v$.
If the diffusion coefficient is small compared to the advective term, e.g $D=5\cdot 10^{-3},\, v=1$, spurious (unphysical) oscillations arise in the neighborhood of $x = 1$ and propagate throughout the entire 
domain polluting the overall quality of the computed solution, see Figure \ref{sfig:no_art_diff}. 
Numerical solutions $\lambda_h{}$ obtained with UP and SG stabilization functions are illustrated in Figures \ref{sfig:UP} and \ref{sfig:SG}, respectively. In the case of the SG stabilization function, the computed $\lambda_h{}$ is \emph{nodally exact}.

\section{Simulation of biological channels}\label{sec:numerical_results}
In this concluding section we carry out a thorough validation of the 
vET and vTHD models in the study of two different biological channels: 
1) the Gramicidin-A channel considered in~\cite{JeromeBPJ}; and 
2) the bipolar nanofluidic diode considered in~\cite{Siwy7}.
The main focus of the simulations is on the current-voltage (IV) 
characteristics of the channel and on how the
IV curves are affected by the boundary conditions, especially the 
temperature of the two bulk regions, $T(0)=T_L$ and $T(d)=T_R$. 
Also, we aim to investigate the electroosmosis effect, 
which is accounted for by the electrolyte fluid velocity $v_e$.
\rs{In all the reported computations we set $\numel= 200$ and use the SG stabilization function~\eqref{eq:SG_peclet}, with the exception of a specific case discussed in Section~\ref{sec:gramicidin}.}
Moreover, in the case 
of the vTHD model we set the
saturation velocity $\vsat$ equal to $10\,\unit{ms^{-1}}$.
We also assume that the channel has a constant cross-sectional area equal to
$\mathcal{A}$, whose value may vary from channel to channel.
For further information on the simulated channels, we 
suggest to consult~\cite{manganini2013} and the references cited therein.

\subsection{Gramicidin-A channel (ballistic diode)}\label{sec:gramicidin}
This channel is thoroughly analyzed in the work~\cite{JeromeBPJ} that
is here used as a benchmark for the biophysical and numerical 
assessment of models and methods proposed in the present article.
To allow comparison between the results of our models and those
of~\cite{JeromeBPJ} we set $v_e = 0$. Unless otherwise specified, 
the subsequent figures refer to computations with the THD model, 
meaning that the outcome of the ET model would appear indistinguishable.
The channel length $d$ is equal to $2.5 \, \unit{nm}$, $\varphi_R$ is always
set equal to $0 \, \unit{V}$ and $\varphi_L = \Vapp$, $\Vapp$ being
the externally applied voltage.
\subsubsection{Electrochemical variables}
\begin{figure}[h!]
\centering
\textbf{Permanent charge and concentrations}\par\medskip
\includegraphics[width=0.49\textwidth]{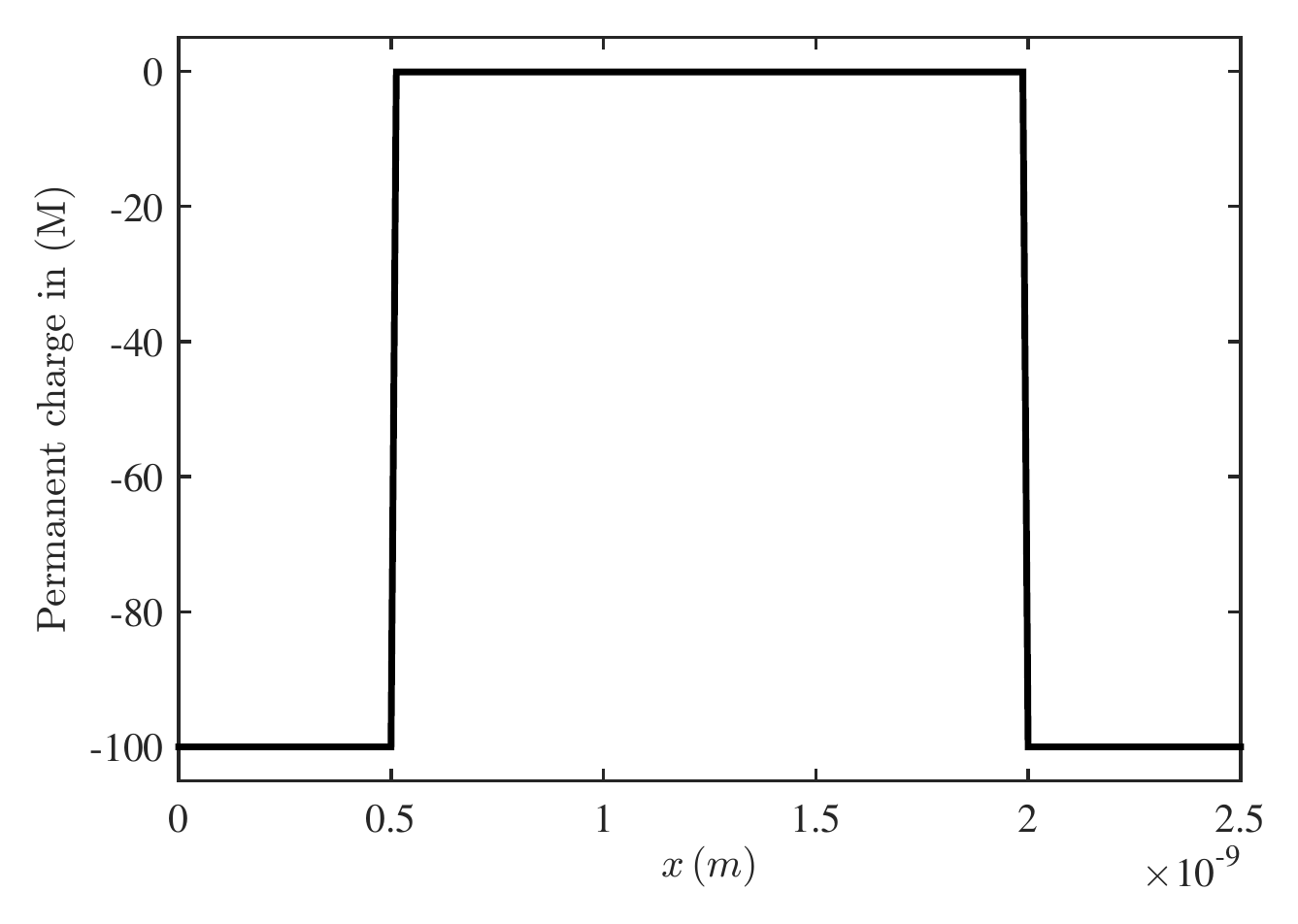}
\includegraphics[width=0.49\textwidth]{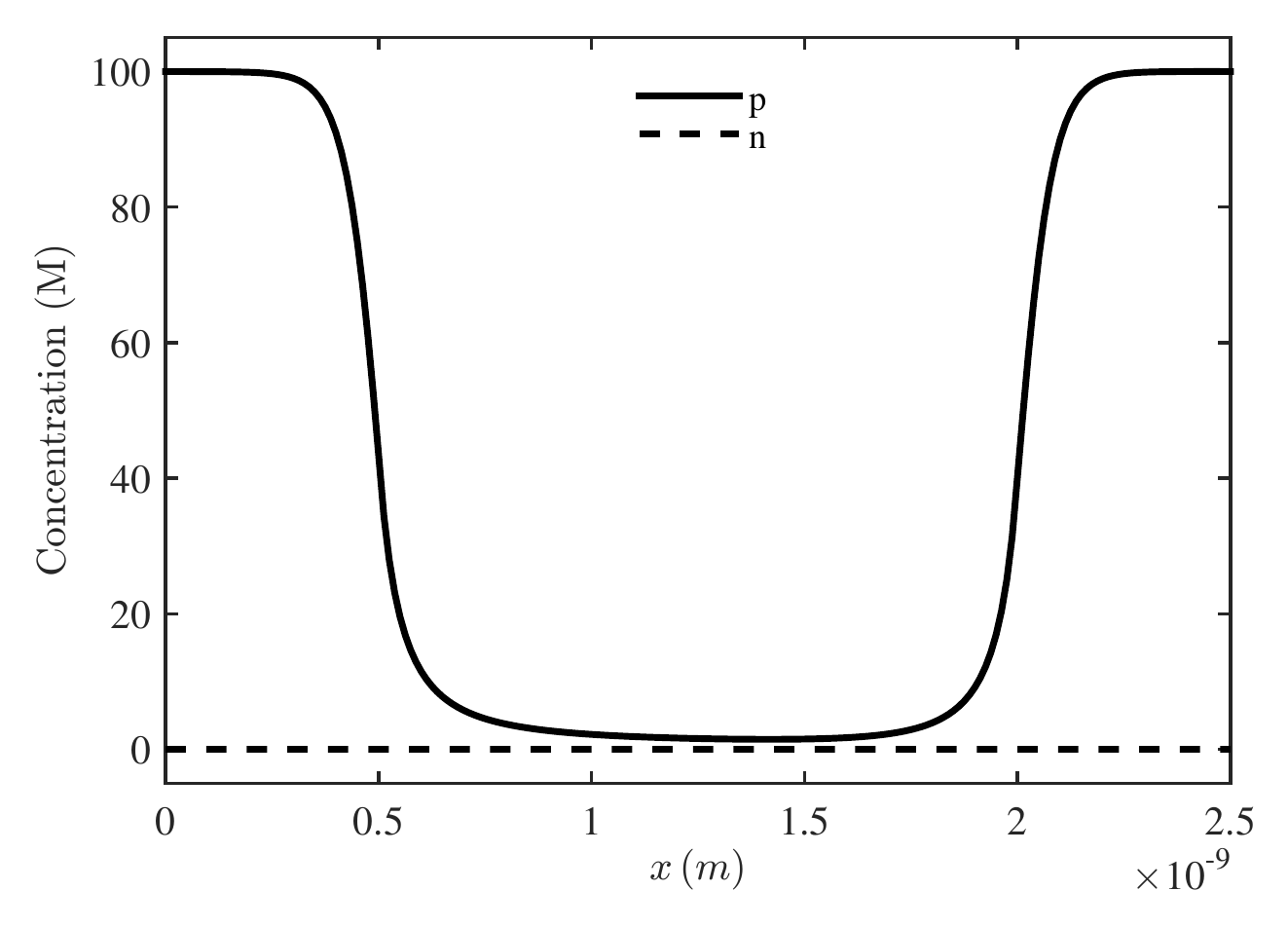}
\caption{Gramicidin-A channel. Left: permanent charge profile.
Right: ion concentrations.}
\label{fig:ballistic_Ndop}
\end{figure}

\begin{figure}[h!]
\centering
\textbf{Concentration profile without and with numerical stabilization}\par\medskip
\includegraphics[width=0.49\textwidth]{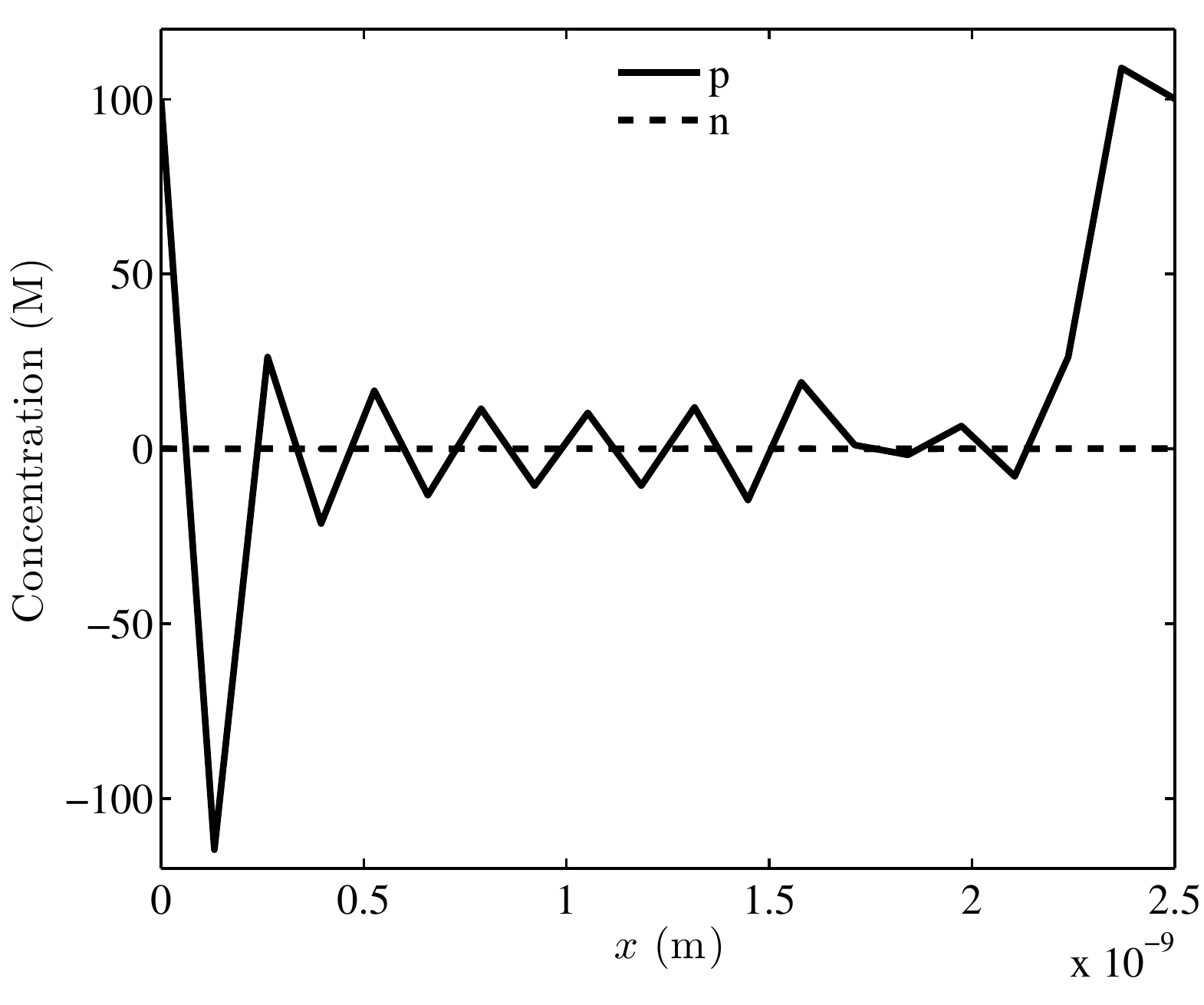}
\includegraphics[width=0.49\textwidth]{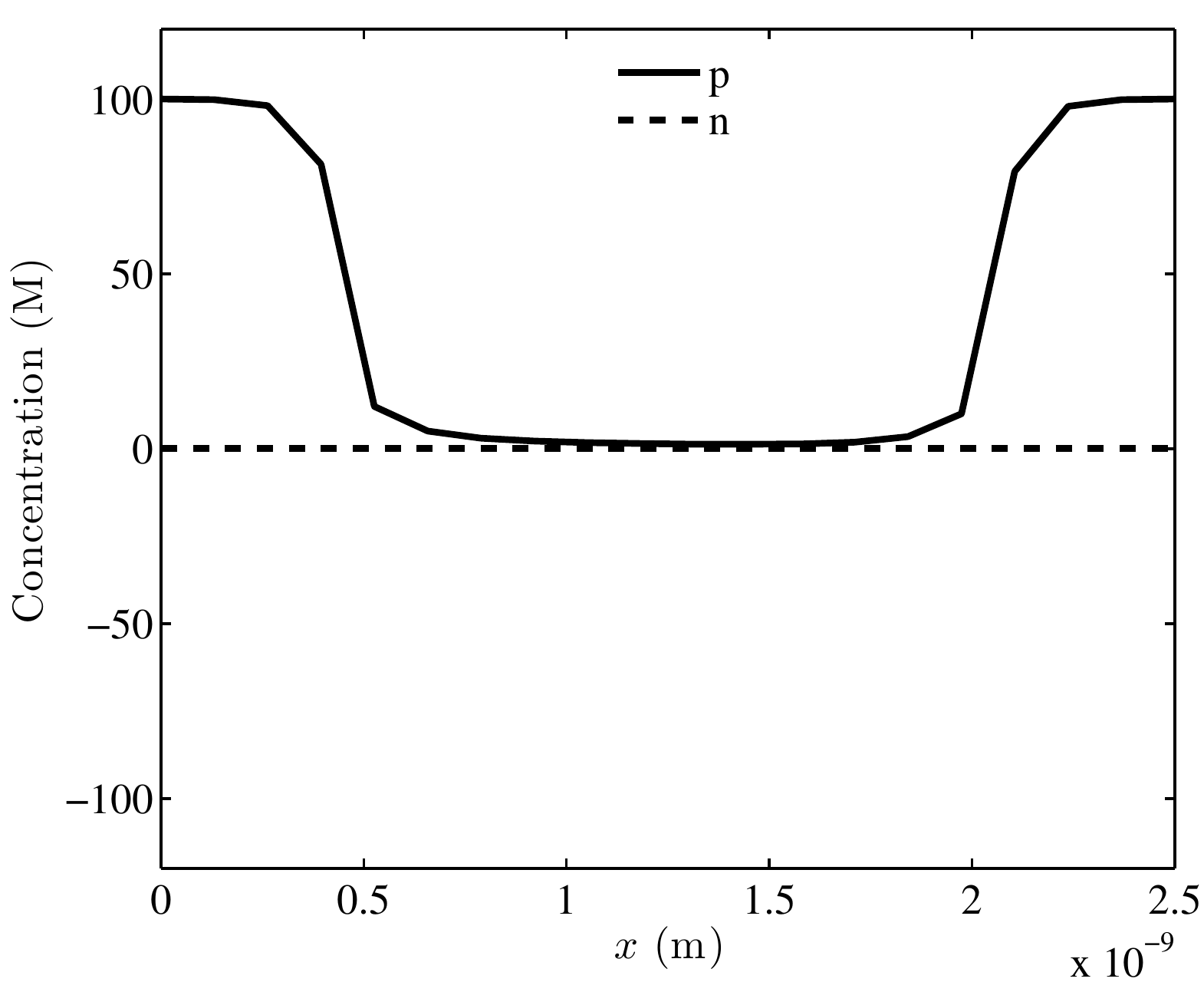}
\caption{\rs{Gramicidin-A channel. Concentration profile with a grid
of 19 elements. When no numerical stabilization is adopted, spurious
oscillations arise (left). The SG stabilization (right) leads to a stable solution.}}
\label{fig:ballistic_conc_no_stab_vs_SG}
\end{figure}

The permanent charge profile of the channel is illustrated in Fig.~\ref{fig:ballistic_Ndop} (left).
Because of the negative fixed charge the channel is highly selective to ion flow and attracts positive $Na^+$ ions while $Cl^-$ ions are mainly repelled.
This behavior is confirmed by the computed ion concentration profiles shown in Fig.~\ref{fig:ballistic_Ndop} (right) in the case where a voltage drop 
$\varphi_L - \varphi_R = \Vapp = 0.1 \, \unit{V}$ 
is applied across the channel.
We first observe that the cation concentration computed by our algorithm is in excellent agreement with that of~\cite{JeromeBPJ}.
We also notice the presence of strong 
internal layers in the $Na^+$ distribution which are effectively captured
by the stabilized DMH discretization without introducing spurious oscillations.
In the mentioned reference~\cite{JeromeBPJ} the hydrodynamic equations 
are solved using a rather different numerical strategy based on the
essentially non oscillatory finite differences with shock capturing
introduced in~\cite{Shu,Shu2}.
\rs{The role of the SG stabilization in the quality of the
computed ion concentrations is well documented by the results shown in Fig.~\ref{fig:ballistic_conc_no_stab_vs_SG}.
In this graph, we see on the left side the cation and anion densities obtained by running the DMH
discretization scheme without stabilization  ($\Phi=0$) on an uniform grid of $N_{el}=19$ elements.
The maximum local P\`eclet number is in this case equal to 68.6288 so that the solution exhibits a
markedly oscillatory behavior and fails to satisfy the DMP.
Conversely, the adoption of the SG stabilization with the same number of elements produces the ion concentrations plotted on the right side of Fig.~\ref{fig:ballistic_conc_no_stab_vs_SG}.
The computed solution is strictly positive and the internal layers
at the interfaces between the channel and highly doped regions are captured with the resolution
allowed by the roughness of the grid size.}
\begin{figure}[h!]
\centering
\textbf{Potential and electric field}\par\medskip
\includegraphics[width=0.49\textwidth]{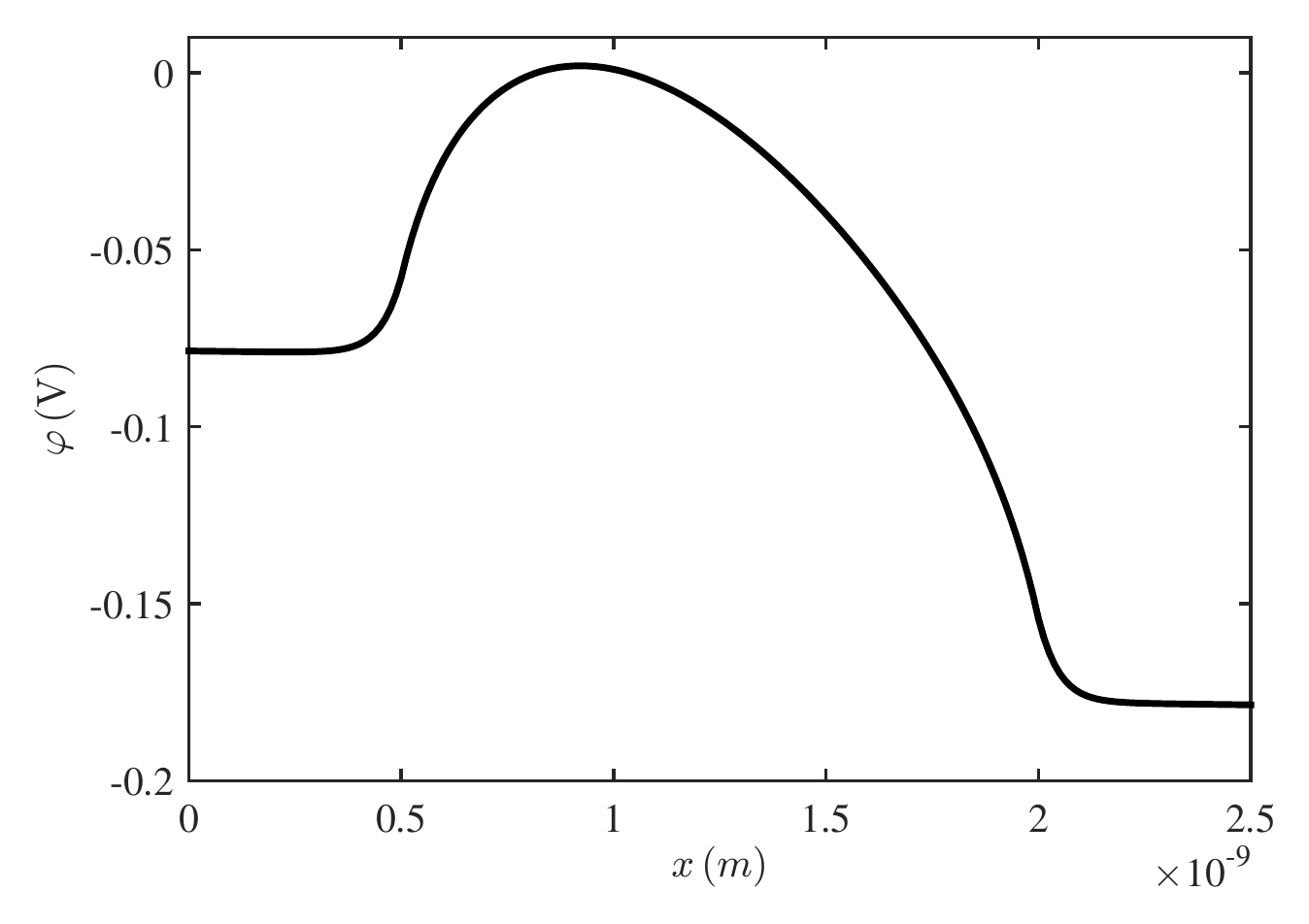}
\includegraphics[width=0.49\textwidth]{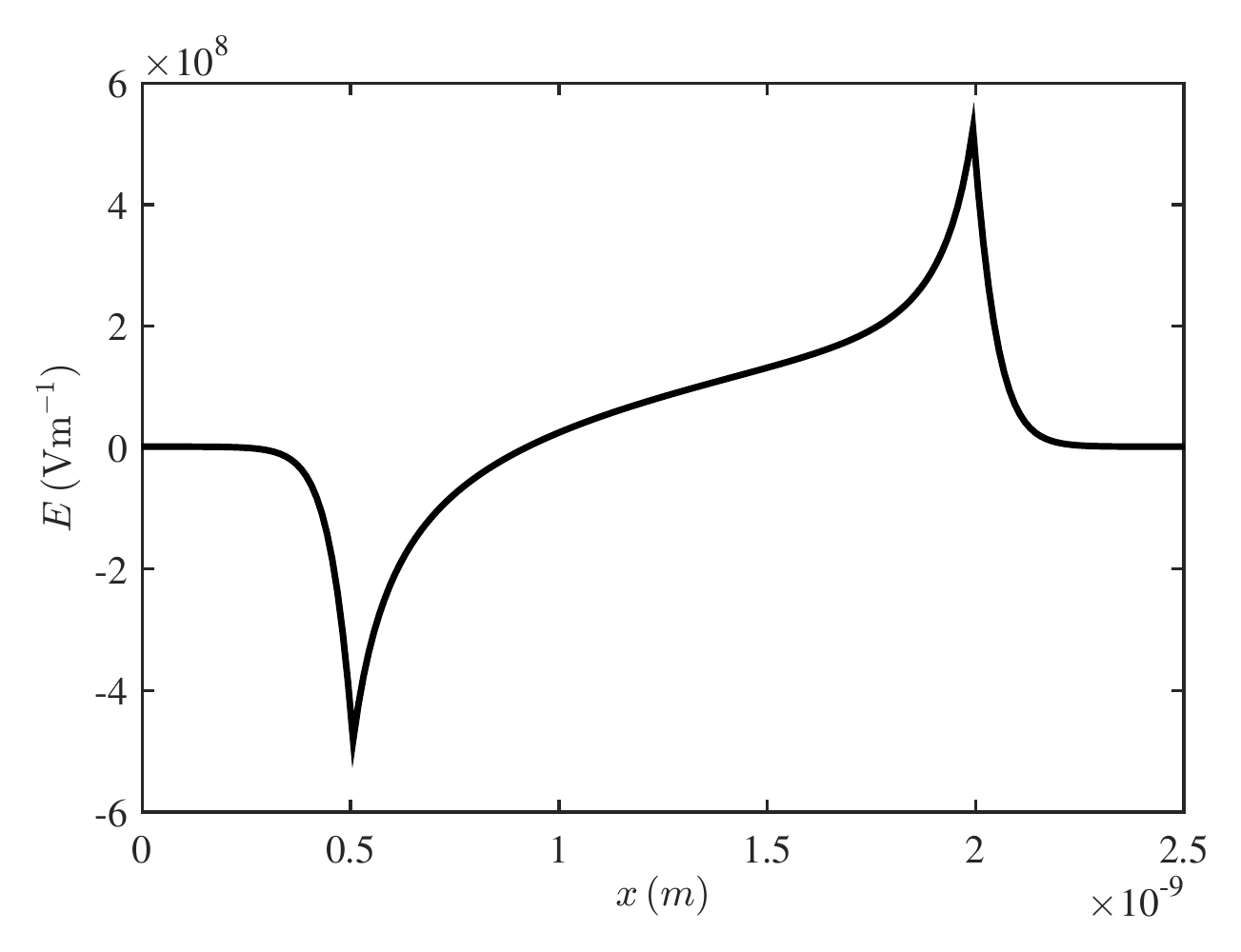}
\caption{Gramicidin-A channel. Left: electric potential.
Right: electric field.}
\label{fig:ballistic_pot}
\end{figure}
\begin{figure}[h!]
\centering
\textbf{IV curves when $T(0)=T(d) \neq 300\,\unit{K}$}\par\medskip
\includegraphics[width=0.49\textwidth]{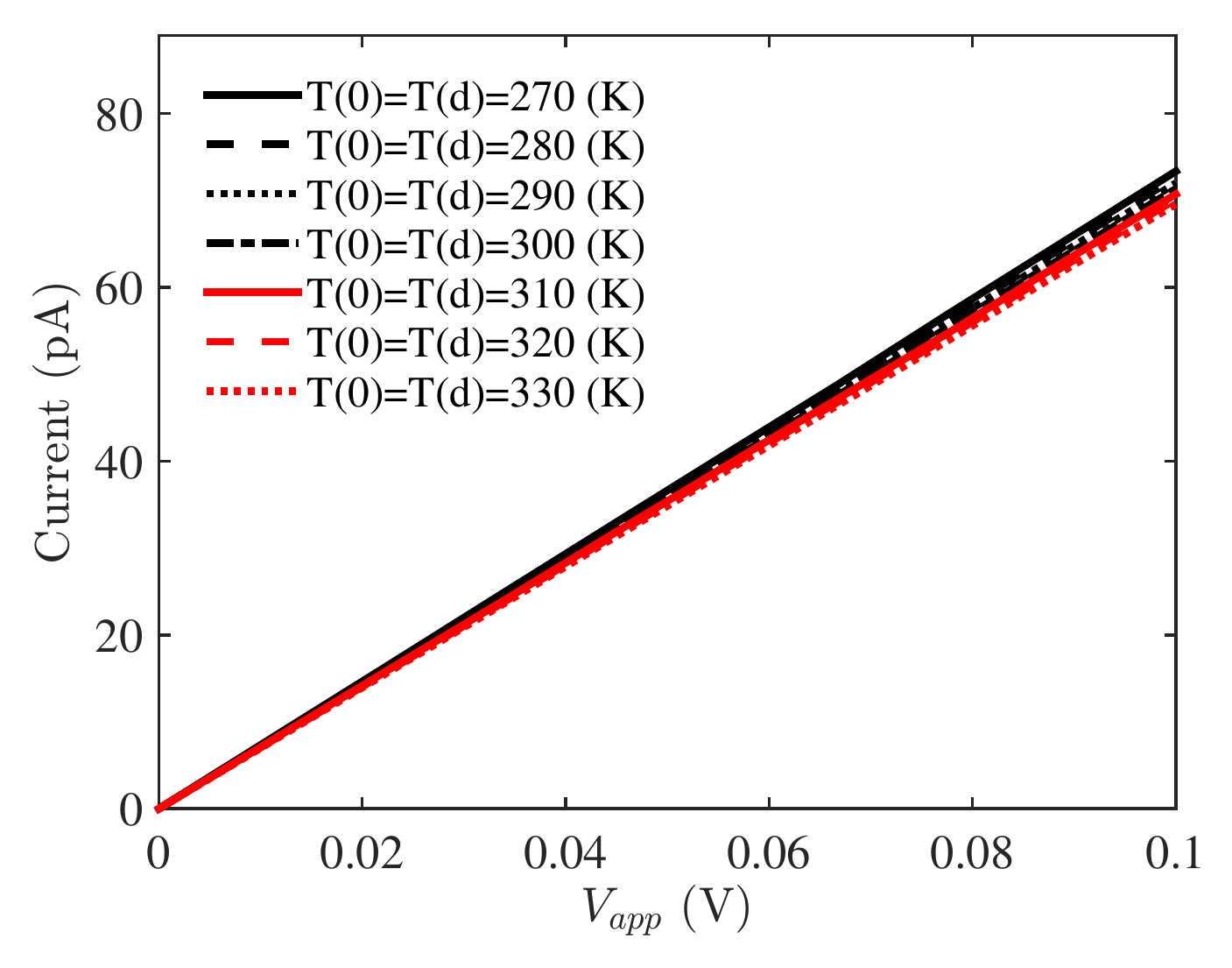}
\includegraphics[width=0.49\textwidth]{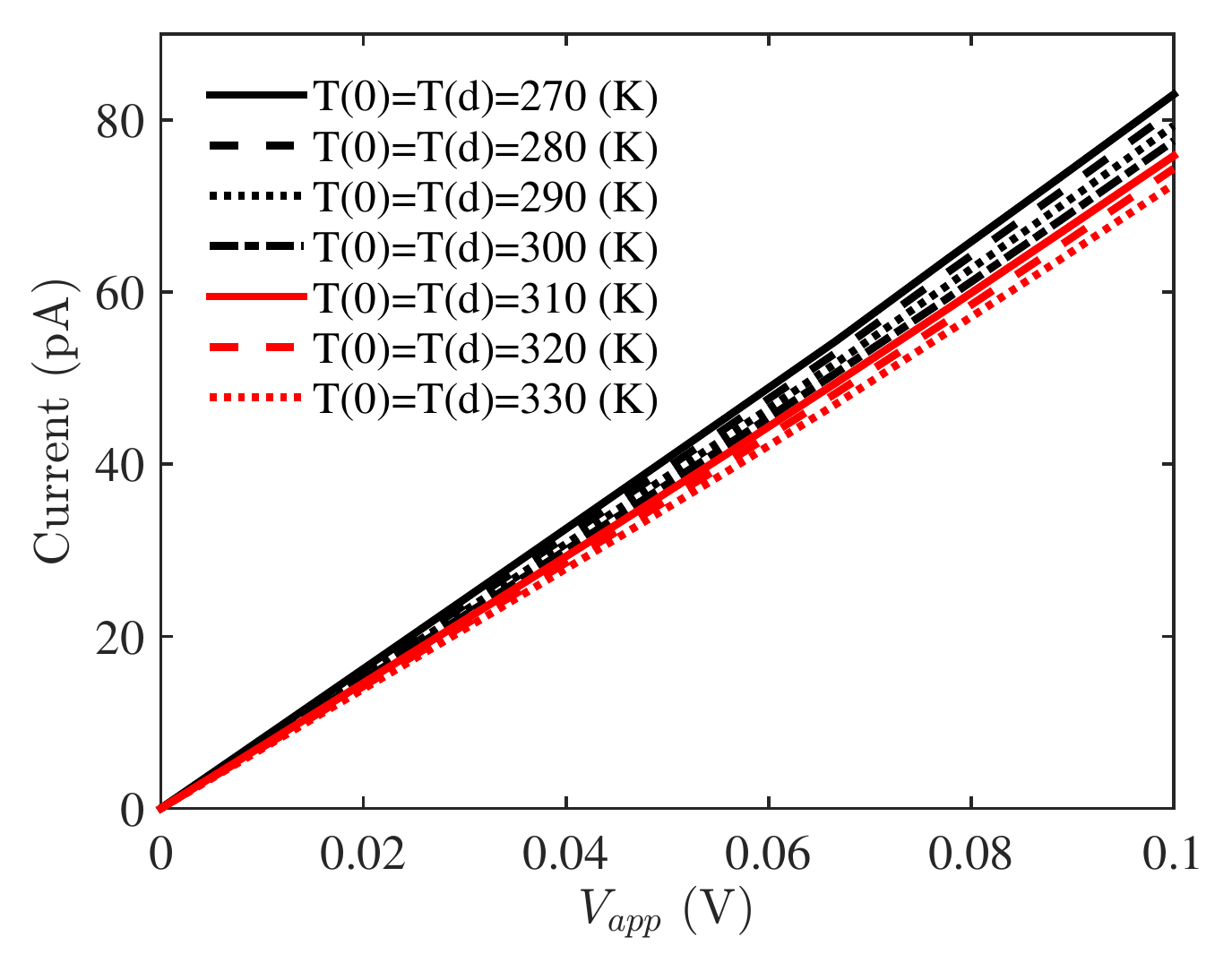}
\caption{Gramicidin-A channel. IV curves when 
$T(0)=T(d) \neq 300 \, K$ (left). THD model (left). 
ET model (right).}
\label{fig:ballistic_IV_THD_T0_Td}
\end{figure}

The permanent charge profile and the consequent distribution of ions allows us to regard the Gramicidin-A channel as the biophysical analogue of an 
electronic "ballistic" diode of type $p^+$-$p$-$p^+$. This analogy is useful
in the interpretation of simulation results and supports the idea 
that "channels are transistors alive" (cf.~\cite{eisenberg_living_transistor}).
In this respect, the region of the channel 
$x \in [0, 0.5] \, \unit{nm}$ corresponds to the Source terminal 
which has the role of emitting cations into the channel, while the
region of the channel $x \in [2, 2.5] \, \unit{nm}$ corresponds to the
Drain terminal which has the role of collecting the cations that 
have travelled throughout the channel under the action of the 
thermo-electro-chemical forces.
Electric potential and electric field profiles are shown in Fig.~\ref{fig:ballistic_pot} (left) and 
Fig.~\ref{fig:ballistic_pot} (right), respectively.
Again, we observe that the electric potential and field computed by our algorithm are in excellent agreement with those of~\cite{JeromeBPJ}.
By inspection of the electric potential distribution, we see that
the positive pole of the bio-electrical system is located at $x=0 , \unit{nm}$
while the negative pole is at $x=2.5 \, \unit{nm}$. Thus, the cations
move from left to right, giving rise to accumulation of (fixed) negative 
charge at the entrance of the channel 
and positive (mobile) charge at the outlet of the channel. 
These accumulations correspond to the two strong peaks visible in 
the electric field distribution.

To investigate the effect of different boundary values 
for the temperature at both entrance and outlet of the channel 
we let $T(0)=T(d)$ range between $270\,\unit{K}$ and $330\,\unit{K}$. 
Fig.~\ref{fig:ballistic_IV_THD_T0_Td} shows
the IV curves for the THD model (left) and the 
ET model (right), respectively. The externally applied voltage 
$\Vapp$ ranges from 0 to 0.1 V. 
Temperature influence is visible in both models: the higher the bath 
temperature the lower the current. This agrees with the 
biophysical insight that an increase of bath temperature 
corresponds to a higher thermal energy for the ions and, consequently,
to a higher energy dissipation because of frictional effects between
particles and fluid. We notice also that the ET model predicts a slightly  higher current than the THD model because this latter model accounts for 
convective energy dissipation within the fluid. 
\subsubsection{Thermal variables}\label{sec:thermal_variables}
\begin{figure}[h!]
\centering
\textbf{Ion temperature when $T(0)$ or $T(d) > 300\,\unit{K}$}\par\medskip
\includegraphics[width=0.49\textwidth]{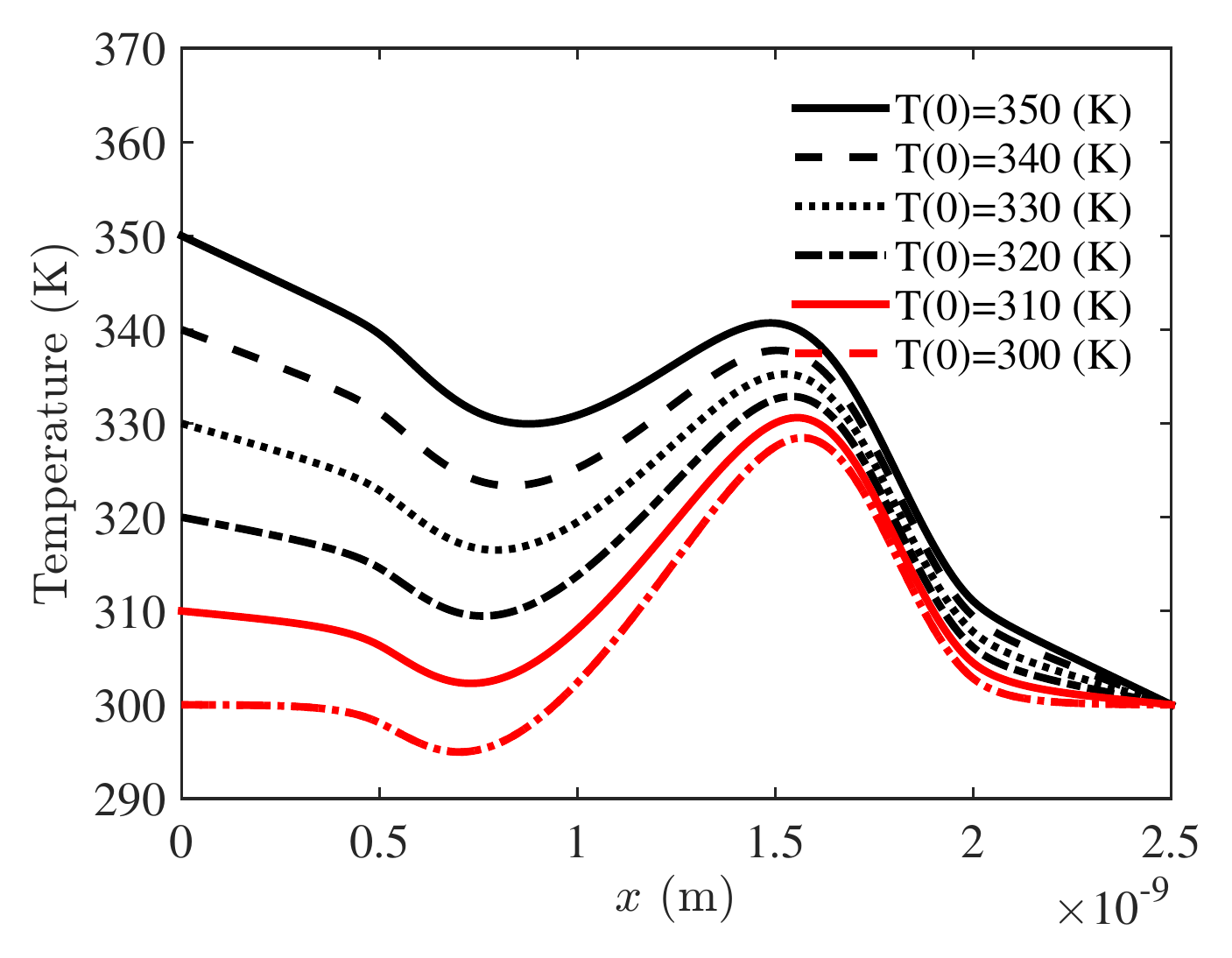}
\includegraphics[width=0.49\textwidth]{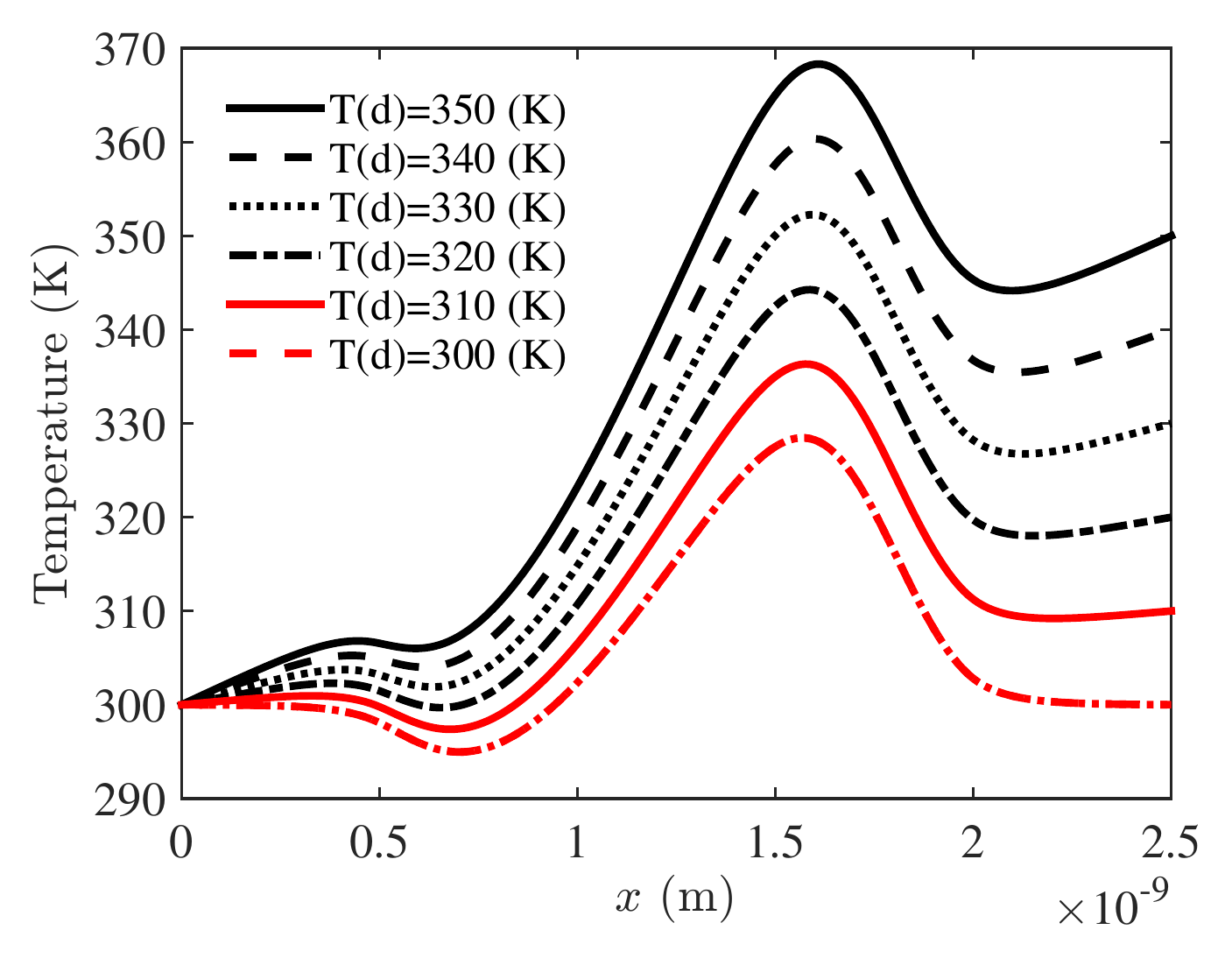}
\caption{Gramicidin-A channel. Cation temperature profiles. 
$T(0)>300 \, K$ (left). $T(d)>300 \, K$ (right).}
\label{fig:ballistic-Tcarr-vs-T0}
\end{figure}

Cation temperature profiles when $T(0)$ and $T(d)$ are 
increased separately, are shown in Fig.~\ref{fig:ballistic-Tcarr-vs-T0}. 

Looking at the two families of distributions, we can distinguish
five distinct subregions in each curve: the Source and Drain regions
(reservoirs), the channel region and the two 
junctions separating Source and Drain from the ion channel. In the 
two reservoirs, the temperature profile is linear, because 
the cation concentration is almost constant and the
electric field is very small. Then, ion temperature increases 
in the channel region according with the fact that particles
are accelerated by the electric field from left to right, 
approximately for  $x \geq 0.8 \, \unit{nm}$ 
(cf.~Fig.~\ref{fig:ballistic_pot}, right). However, the effect of the
electric field at the two junctions is quite different in the
two sets of thermal boundary conditions. 
If $T(0) > T(d)$, electric and thermal fields 
act in the same direction (from left to right) so that
as $T(0)$ increases, the thermal flow at the two channel boundaries 
correspondingly increases. If $T(d) > T(0)$ electric and thermal 
forces act in opposite directions and this contributes to 
diminishing the thermal flow at the two channel boundaries. 
It is interesting to notice that ion heating 
(i.e., the maximum value of ion temperature inside the channel) 
is considerably larger in the case where $T(d) > T(0)$,
because in this condition ion acceleration due to the electric field 
greatly dominates over the thermal field along the channel (moving from
left to right) so that the ion total energy increases and temperature 
gets larger. Then, once the ions are injected across 
the junction between channel and Drain they immediately thermalize to
the local energy of the reservoir distribution and approximately 
cool down to the temperature enforced at $x=d$.

\begin{figure}[h!]
\centering
\textbf{Electrolyte fluid temperature when $T(0)$ or $T(d) > 300\,\unit{K}$}\par\medskip
\includegraphics[width=0.49\textwidth]{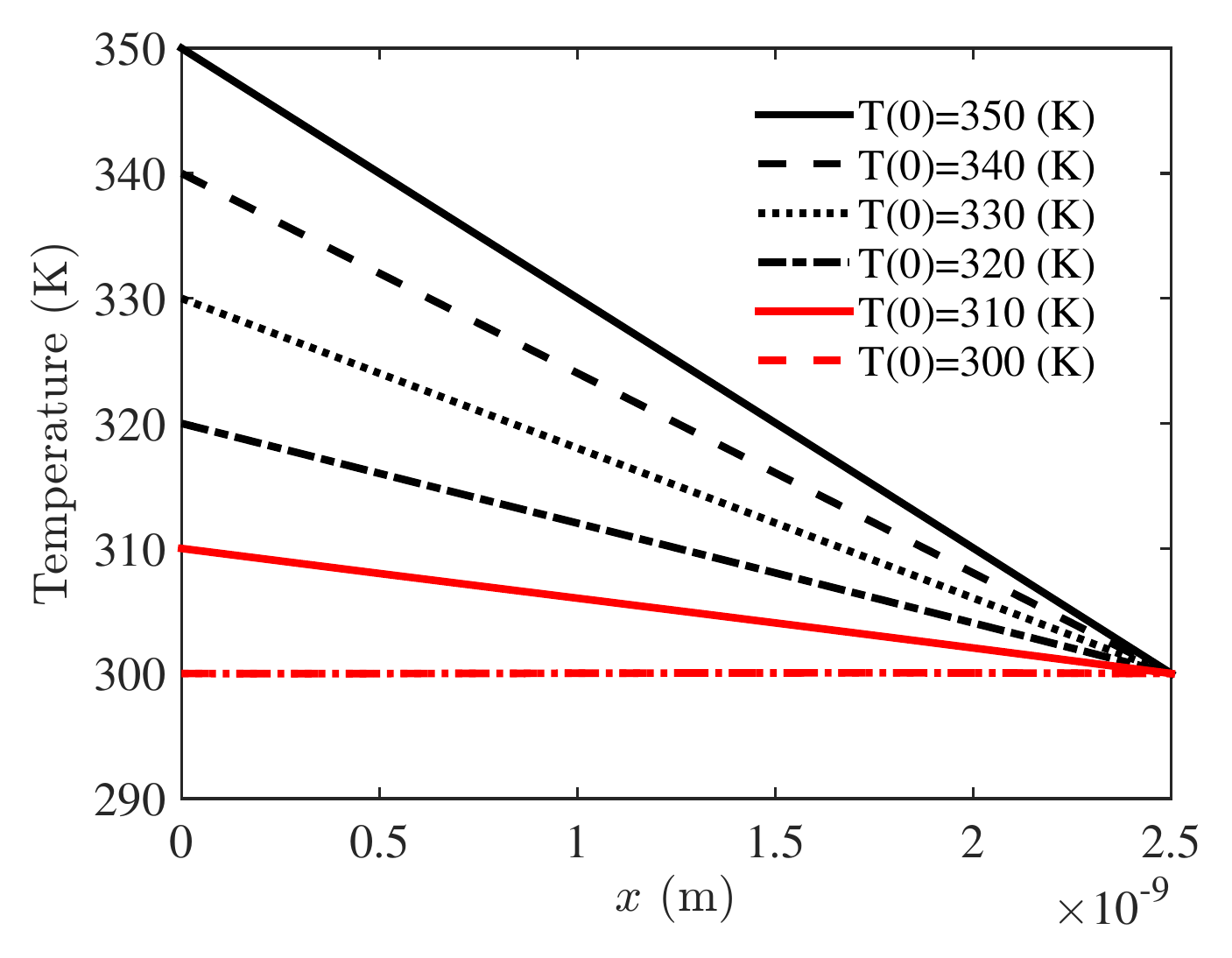}
\includegraphics[width=0.49\textwidth]{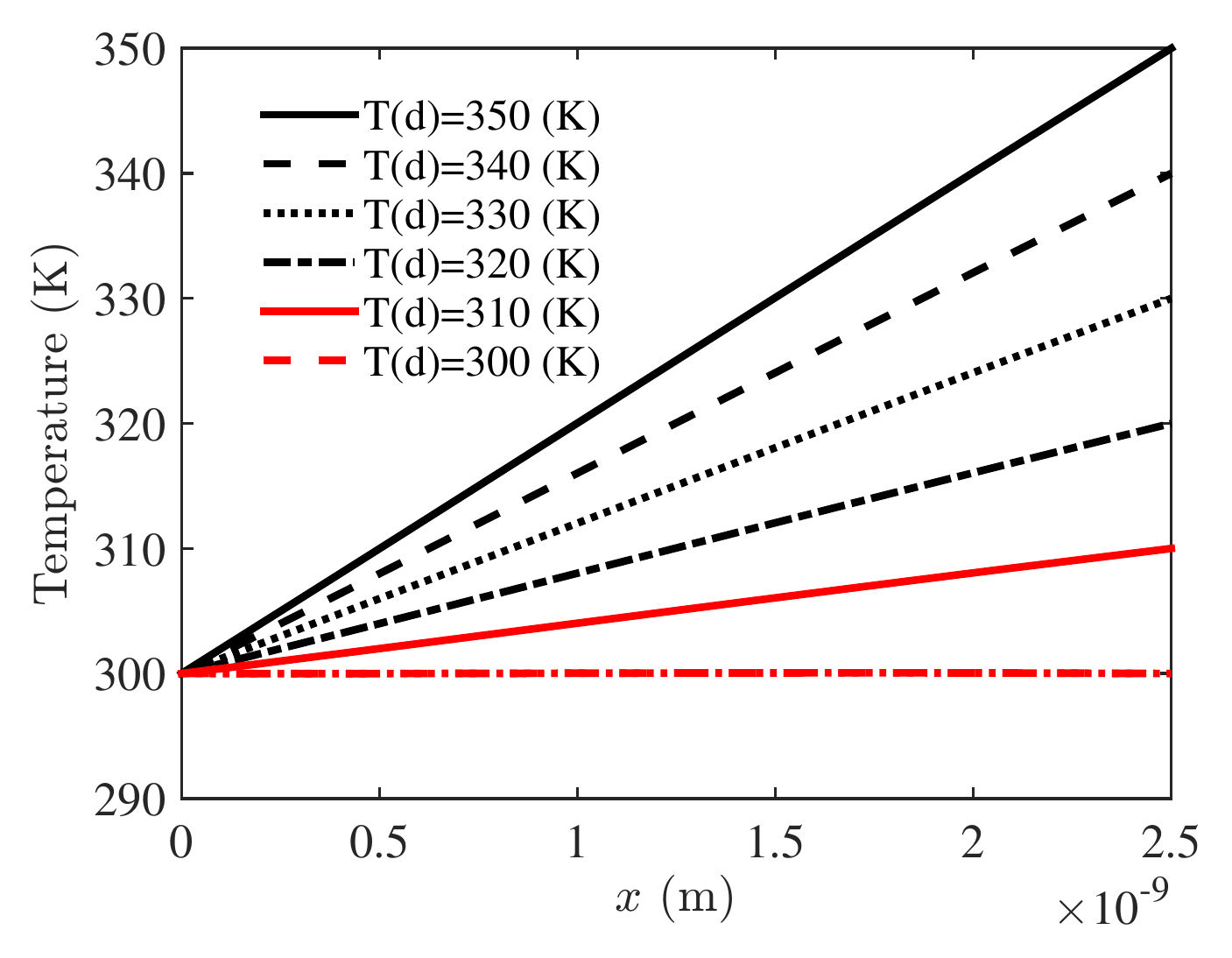}
\caption{Gramicidin-A channel. Electrolyte fluid temperature profiles. 
$T(0)>300 \, K$ (left). $T(d)>300 \, K$ (right).}
\label{fig:ballistic-THD-TL-vs-T0}
\end{figure}

Temperature profiles of the electrolytic fluid are almost linear 
for both choices of the applied thermal drop, see Fig.~\ref{fig:ballistic-THD-TL-vs-T0}. This behaviour is to be ascribed 
to the high value of thermal conductivity of the electrolyte fluid
compared to that of the ion fluid. This makes the fluid behave
as a perfect sink so that its heating is only passively driven by an external
temperature gradient according to Fourier's law~\eqref{eq:THD_sL_1D}
($v_e=0$ in this case).

\subsubsection{Ion velocity}
From ion velocity profiles we can inspect how an externally 
applied temperature gradient may affect ion motion in the channel.
\begin{figure}[h!]
\centering
\textbf{Ion velocity when $T(0)$ or $T(d) > 300\,\unit{K}$. THD model}\par\medskip
\includegraphics[width=0.49\textwidth]{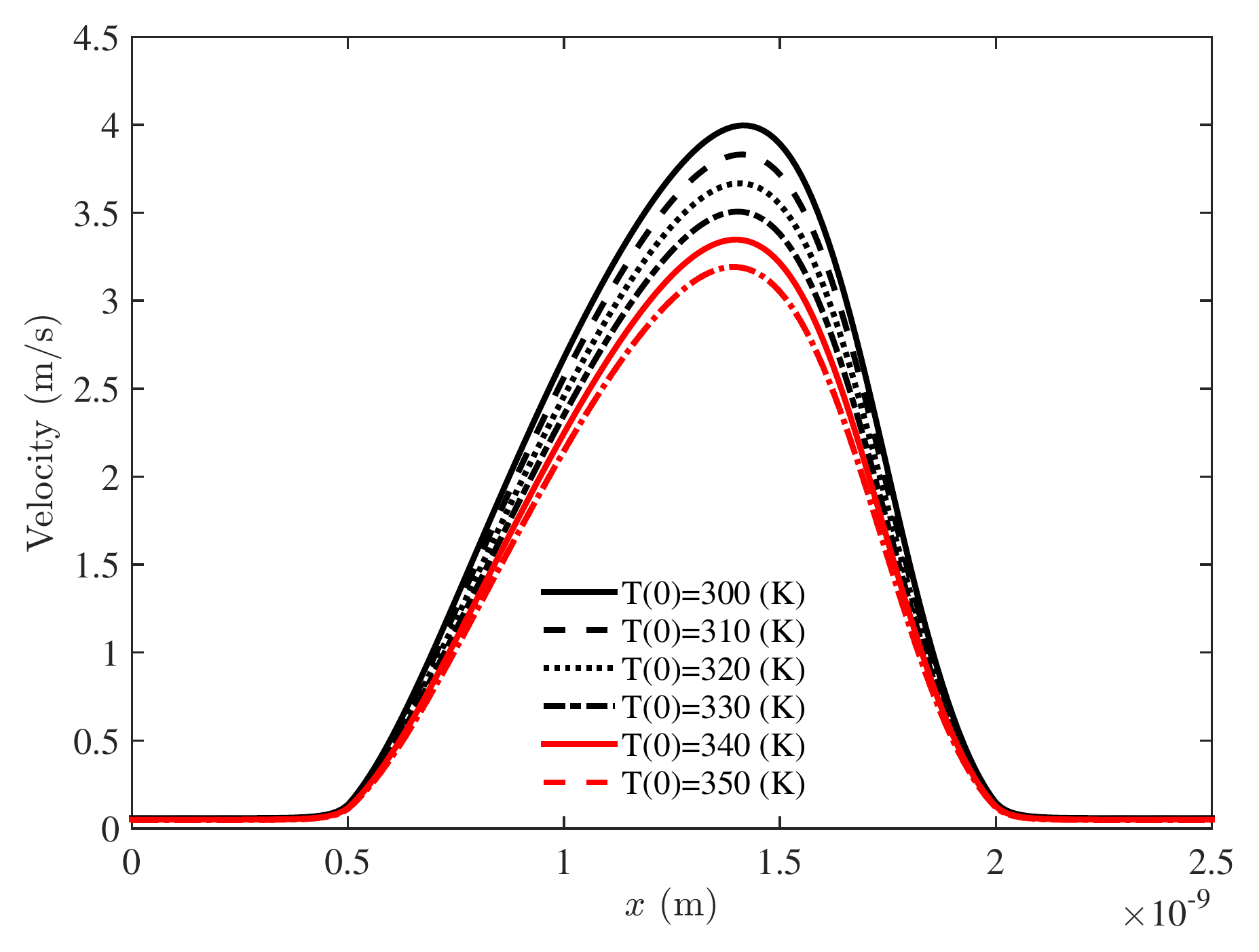}
\includegraphics[width=0.49\textwidth]{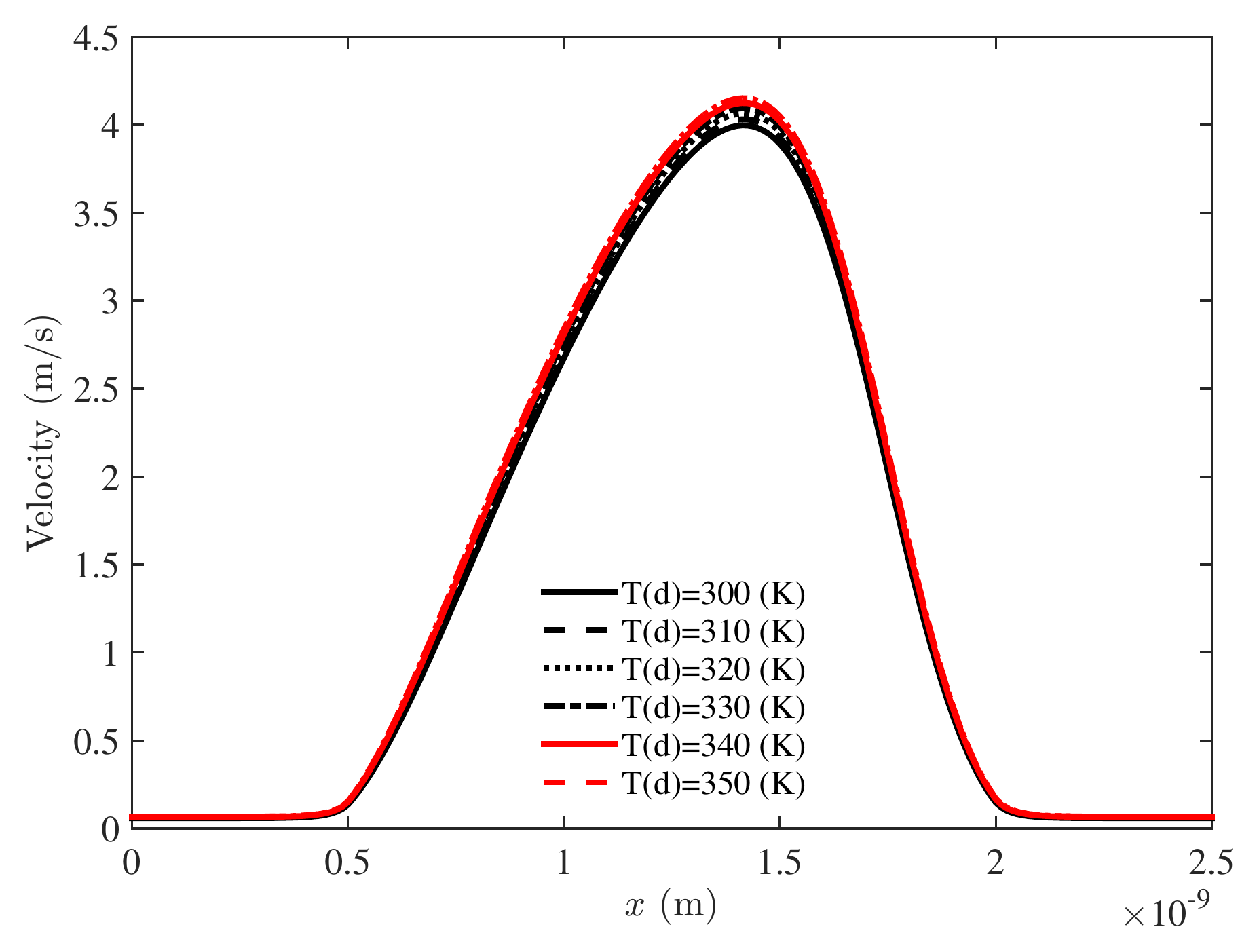}
\caption{Gramicidin-A channel. Cation velocity computed by the THD model.
$T(0)>300 \, K$ (left), $T(d)>300 \, K$ (right).}
\label{fig:ballistic-vel-THD-vs-T0}
\end{figure}
\begin{figure}[h!]
\centering
\textbf{Ion velocity when $T(0)$ or $T(d) > 300\,\unit{K}$. ET model}\par\medskip
\includegraphics[width=0.49\textwidth]{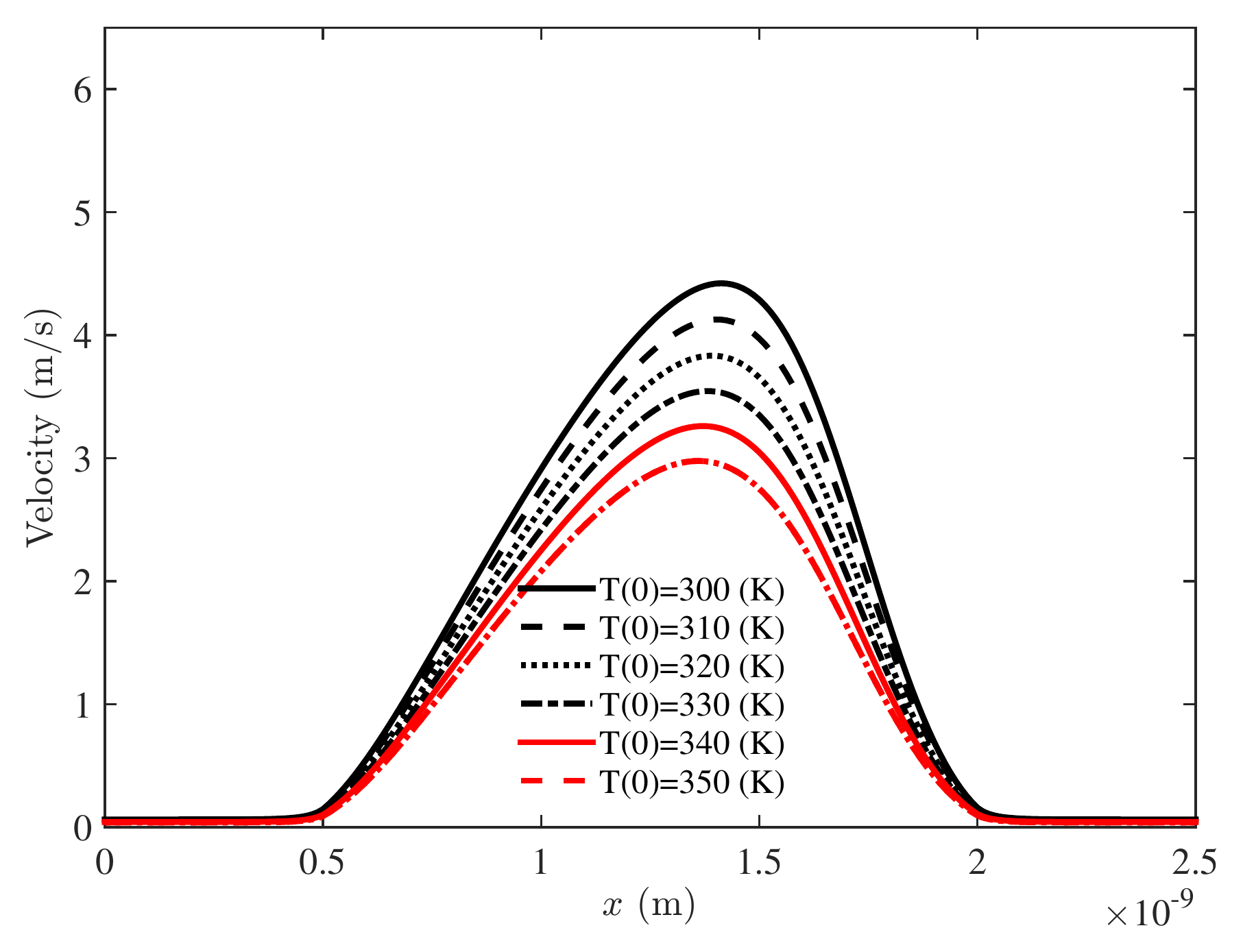}
\includegraphics[width=0.49\textwidth]{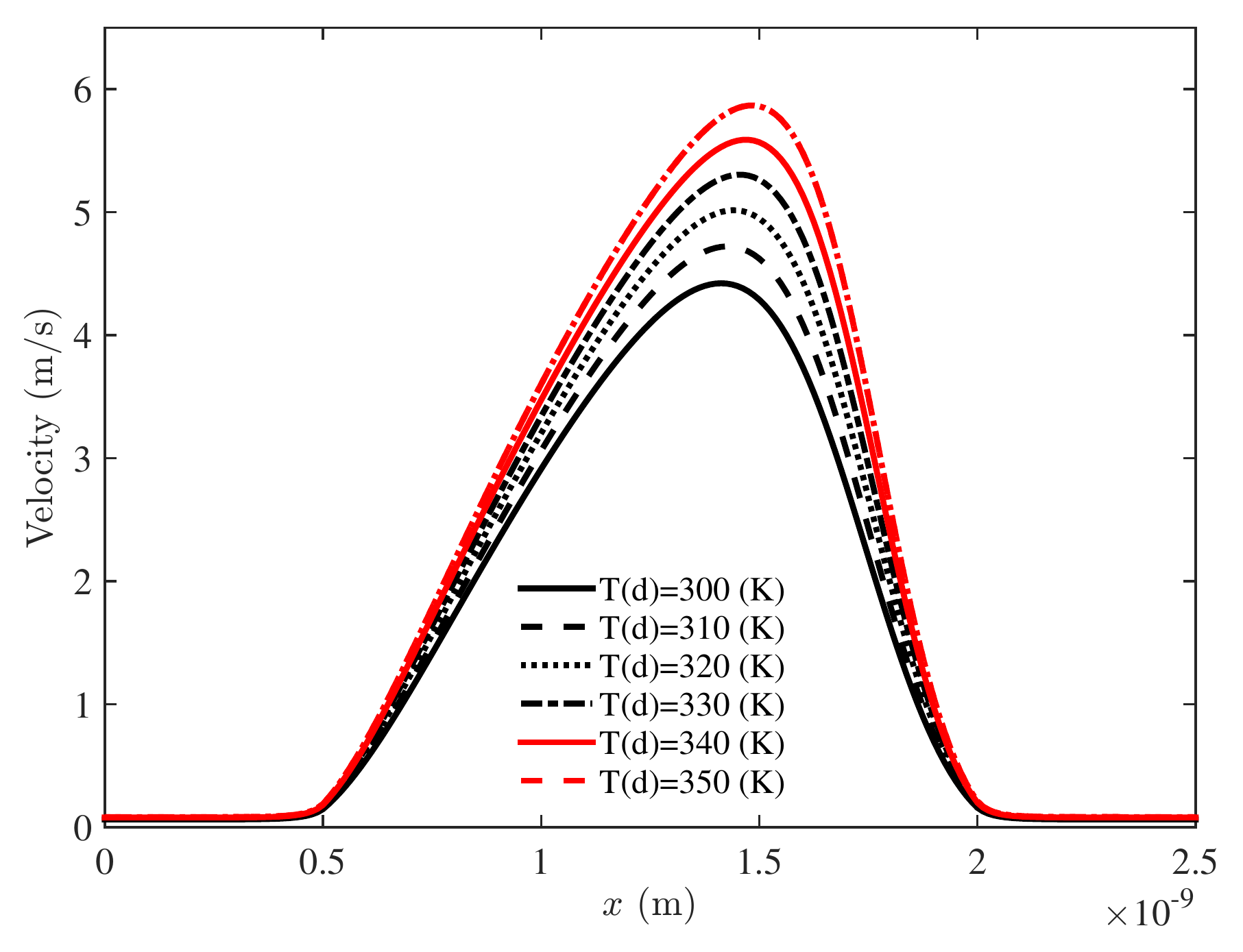}
\caption{Gramicidin-A channel. Cation velocity computed by the ET model.
$T(0)>300 \, K$ (left), $T(d)>300 \, K$ (right).}
\label{fig:ballistic-vel-ET-vs-T0}
\end{figure} 

In Fig.~\ref{fig:ballistic-vel-THD-vs-T0} (left)
we see ions moving slower when $T(0) > T(d)=300\,\unit{K}$.
The opposite happens in Fig.~\ref{fig:ballistic-vel-THD-vs-T0} (right) 
where $T(d) > T(0)=300\,\unit{K}$.
However, in this latter case, velocity variations are much smaller.
Again, the different behavior of the ion velocity 
in the two sets of thermal conditions is related to the interplay
between electric and thermal forces already discussed in 
Sect.~\ref{sec:thermal_variables}. If $T(0) > T(d)$ thermal diffusion
is larger than in the case where $T(d) > T(0)$, and increases 
as $T(0)$ increases. This explains the maximum value of $v_n$ at 
$T(0) = 300 \,\unit{K}$. If $T(d) > T(0)$ ion motion is mainly driven
by the electric field, which explains why the various curves 
are substantially insensitive to the increase of $T(d)$.
The ion velocity in the ET model exhibits a trend similar to that 
of the THD model, see Fig.~\ref{fig:ballistic-vel-ET-vs-T0},
even if the distribution of the peak value shows a much 
larger spread than in the case of the velocity computed by the 
THD model. This is to be ascribed to the larger thermal conductivity
of the (thermally) unified system composed by anions, cations and fluid.

\begin{figure}[h!] 
\centering
\textbf{IV curve when $T(0)$ or $T(d) > 300\,\unit{K}$. THD model}\par\medskip
\includegraphics[width=0.49\textwidth]{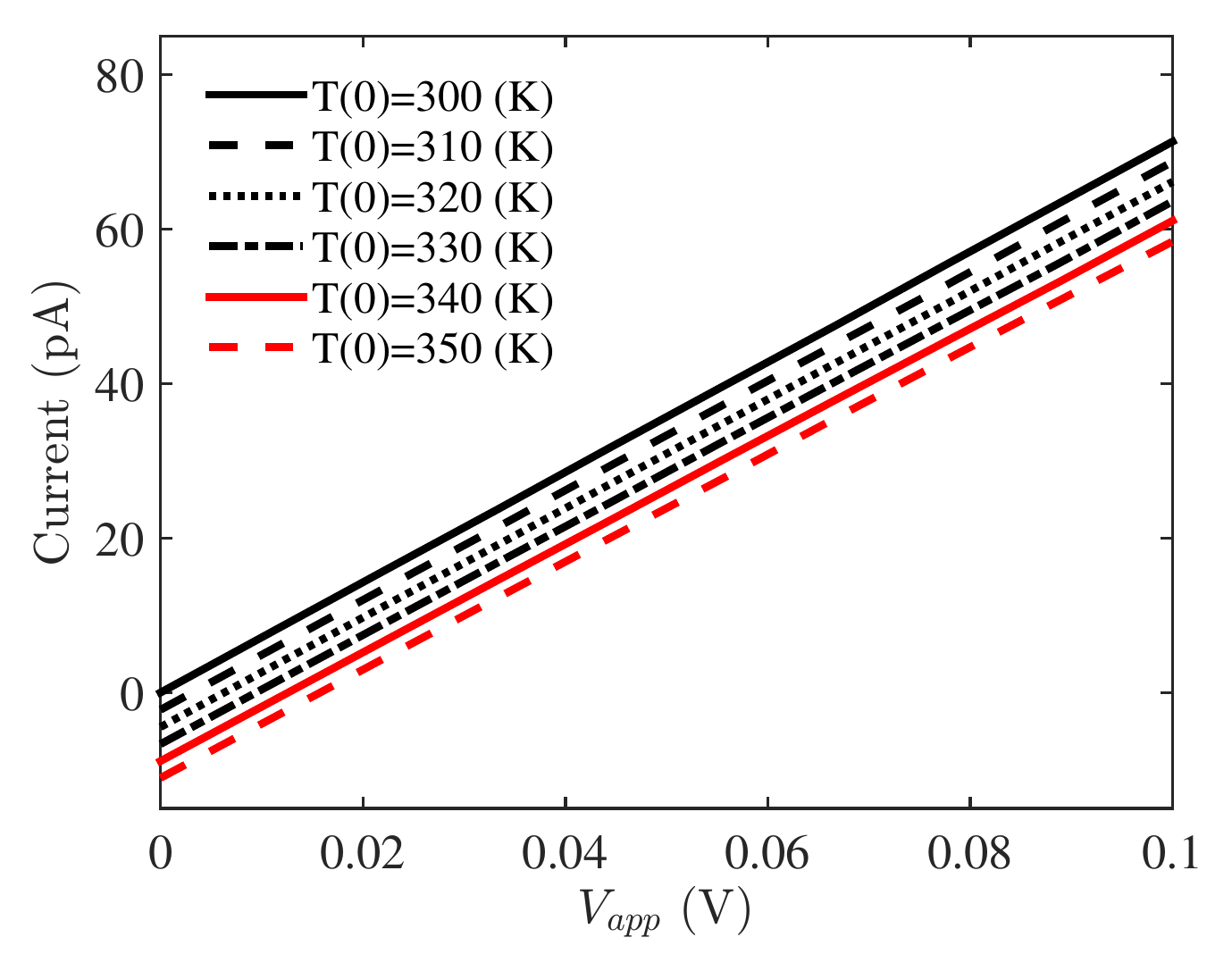}
\includegraphics[width=0.49\textwidth]{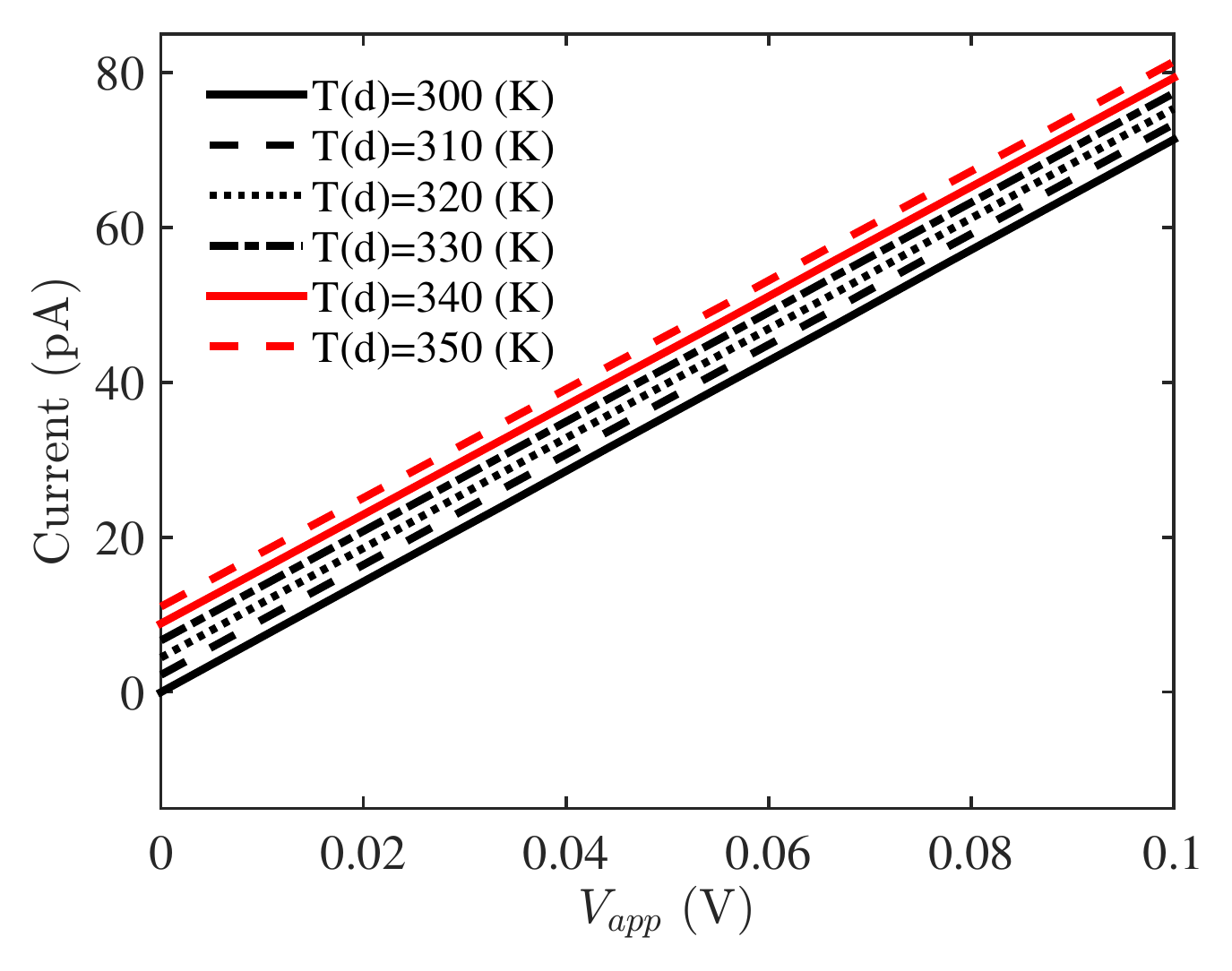}
\caption{Gramicidin-A channel. I-V curves computed by the THD model.
$T(0)>300 \, K$ (left), $T(d)>300 \, K$ (right).}
\label{fig:ballistic-IV-THD-vs-T0}
\end{figure}
\begin{figure}[h!]
\centering
\textbf{IV curve when $T(0)$ or $T(d) > 300\,\unit{K}$. ET model}\par\medskip
\includegraphics[width=0.49\textwidth]{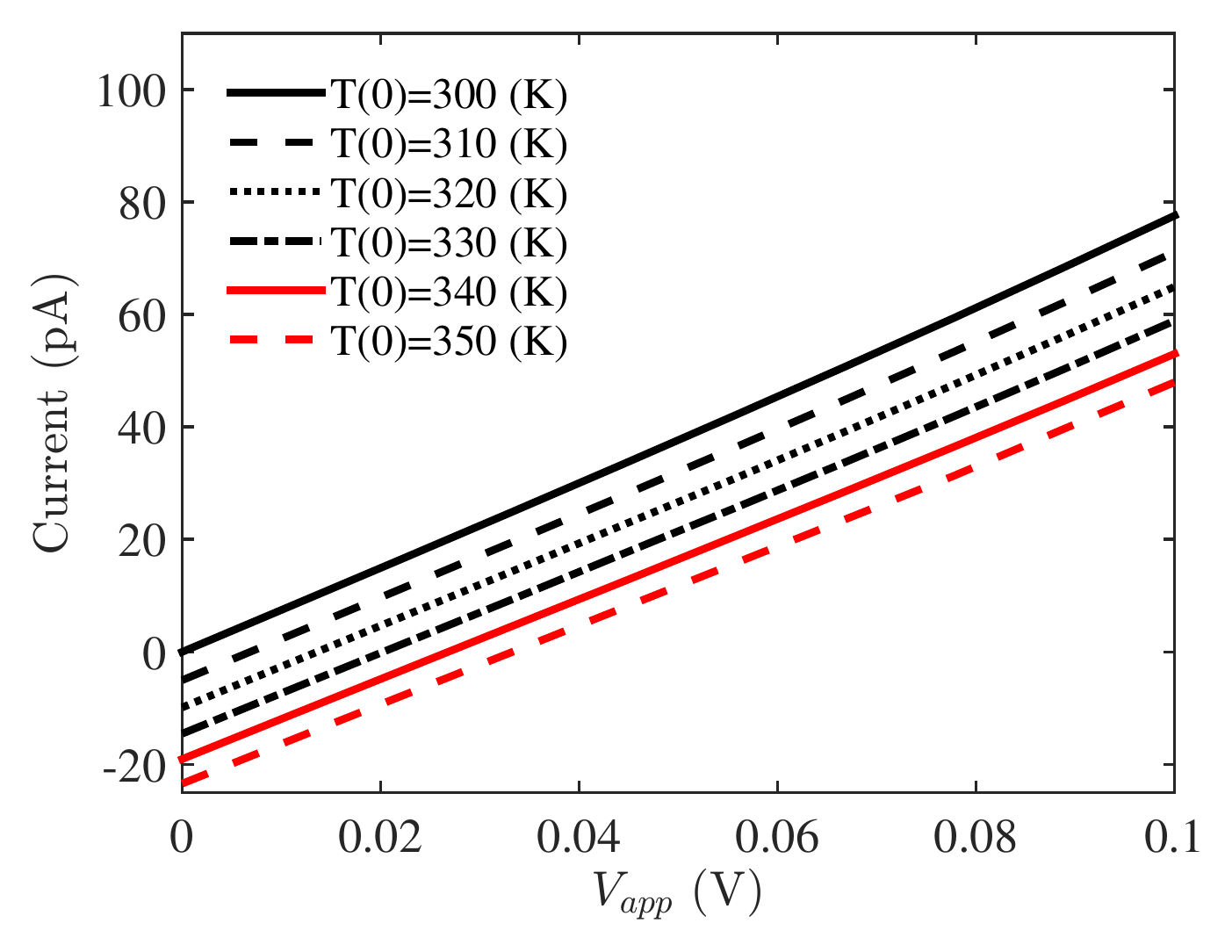}
\includegraphics[width=0.49\textwidth]{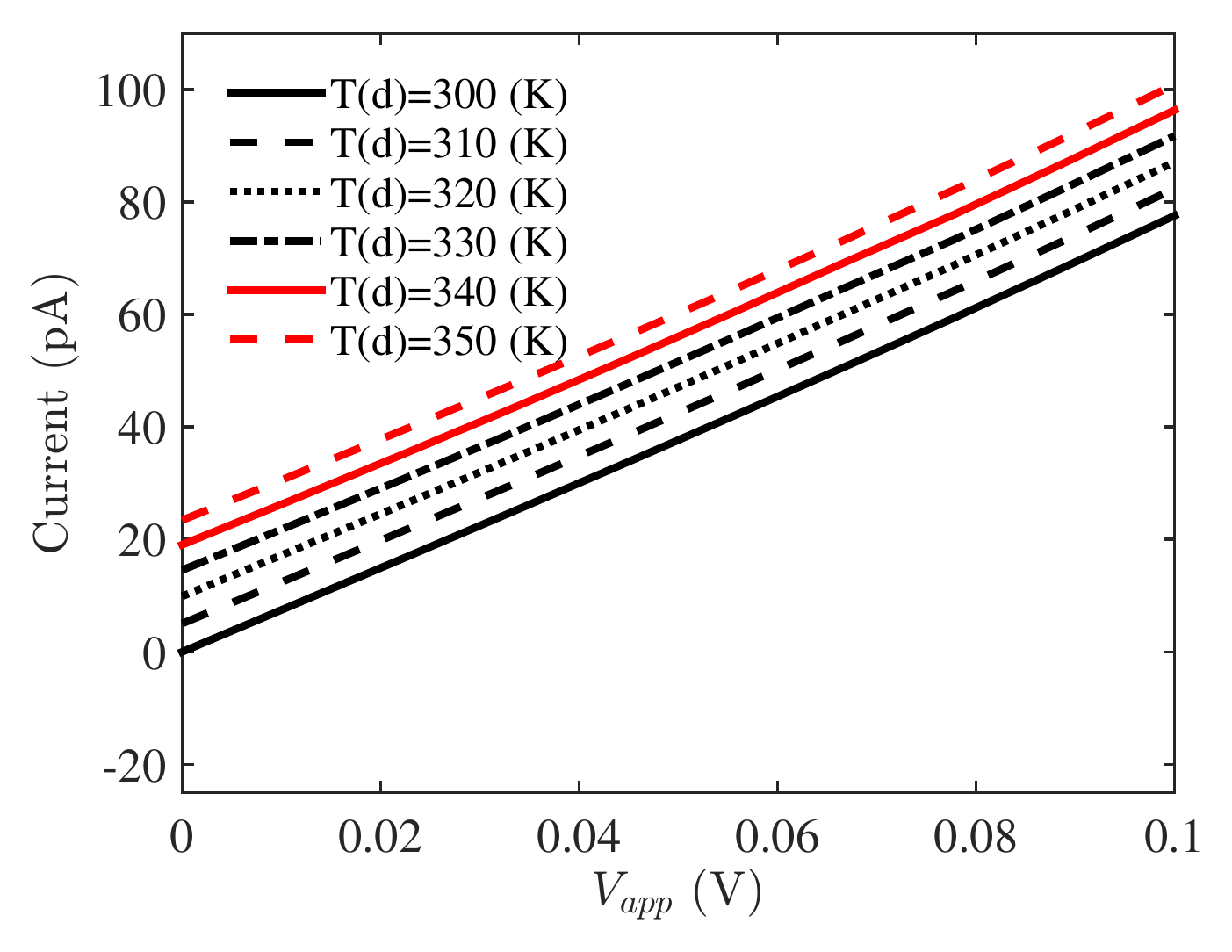}
\caption{Gramicidin-A channel. I-V curves computed by the ET model.
$T(0)>300 \, K$ (left), $T(d)>300 \, K$ (right).}
\label{fig:ballistic-IV-ET-vs-T0}
\end{figure}

The IV curves computed by the two models in the two distinct 
thermal sets of boundary conditions reflect what reported 
about ion velocities:
a smaller ion velocity means a smaller current flowing in the channel.
The externally applied voltage $\Vapp$ ranges from 0 to 0.1 V. 
Fig.~\ref{fig:ballistic-IV-THD-vs-T0} (left) shows 
the IV curves for the THD model when $T(0)>T(d)=300\,\unit{K}$.
Notice that the current assumes negative values when $\Vapp$ is close to zero and the temperature gradient due to thermal BCs is strong enough to 
counterbalance the affect of the applied voltage.
Consistently with the ion velocity profiles of Fig.~\ref{fig:ballistic-vel-THD-vs-T0} (left), 
the current flowing in the channel is lower than in 
the case of homogeneous thermal BCs. 
As $T(0)$ increases, the I-V profile is shifted down 
with respect to the case of homogeneous BCs 
($T(0)=300\,\unit{K}$ solid black line).
The case when $T(d)>T(0)=300\,\unit{K}$ is shown in Fig.~\ref{fig:ballistic-IV-THD-vs-T0} (right).
In this situation the IV relationships are shifted up, meaning that as $T(d)$ increases, a higher current flows in the channel, just as predicted from the ion velocity profiles of Fig.~\ref{fig:ballistic-vel-THD-vs-T0} (right).
Similar results and observations hold for the IV curves 
predicted by the ET model, see Figs.~\ref{fig:ballistic-IV-ET-vs-T0}.

\subsection{IV curves with velocity extended THD and ET models}\label{sec:gramicidin-uNS}
In the simulations conducted so far, we set $v_e=0$, 
which corresponds to neglecting the electroosmosis effect.
We have also conducted simulations with $v_e$ varying in the range 
$0 \div 0.01\, \unit{m\,s^{-1}}$, and we found that it has
no appreciable influence on the current flowing in the channel, so that no results are shown in this case.
However, different conclusions are drawn in the channel 
presented in the following section.

\subsection{Bipolar nanofluidic diode}\label{sec:BP}
The nanofluidic channel investigated in the present section (whose detailed representation can be found in~\cite{Siwy7,manganini2013}) is called 
Bipolar (BP) nanofluidic diode because of its similarity with a $p$-$n$ 
electronic diode: 
as a matter of fact, the BP channel can operate into two different states, depending on the sign of the applied voltage: (1) an 'open' state  (also called 'forward' bias), characterized by high current flowing through the device, and (2) a 'closed' state (also called 'reverse' bias), with very little current flowing in the channel.


\subsubsection{Electrochemical variables}
\begin{figure}[h!]
\centering
\textbf{Permanent charge and concentrations}\par\medskip
\includegraphics[width=0.32\textwidth]{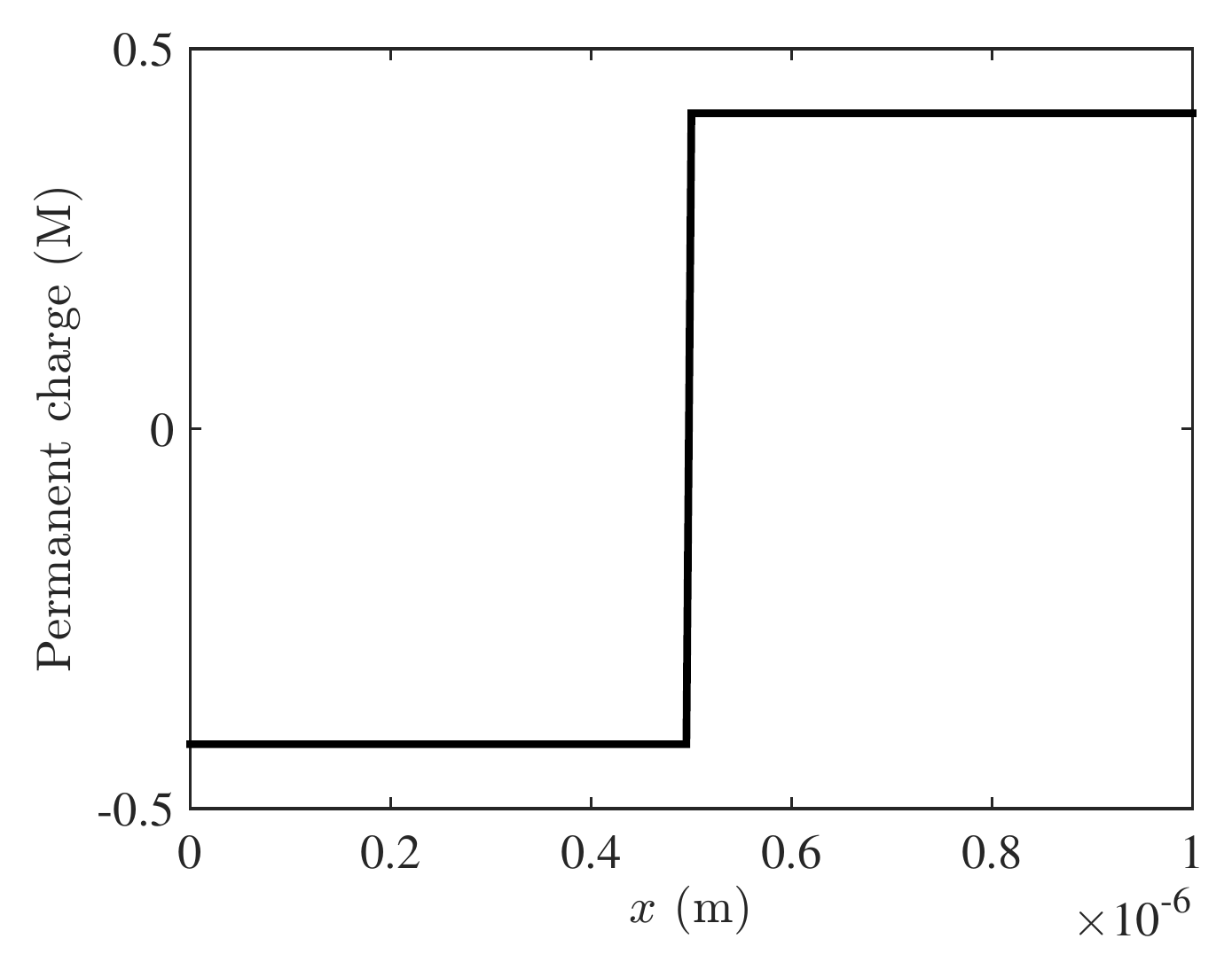}
\includegraphics[width=0.32\textwidth]{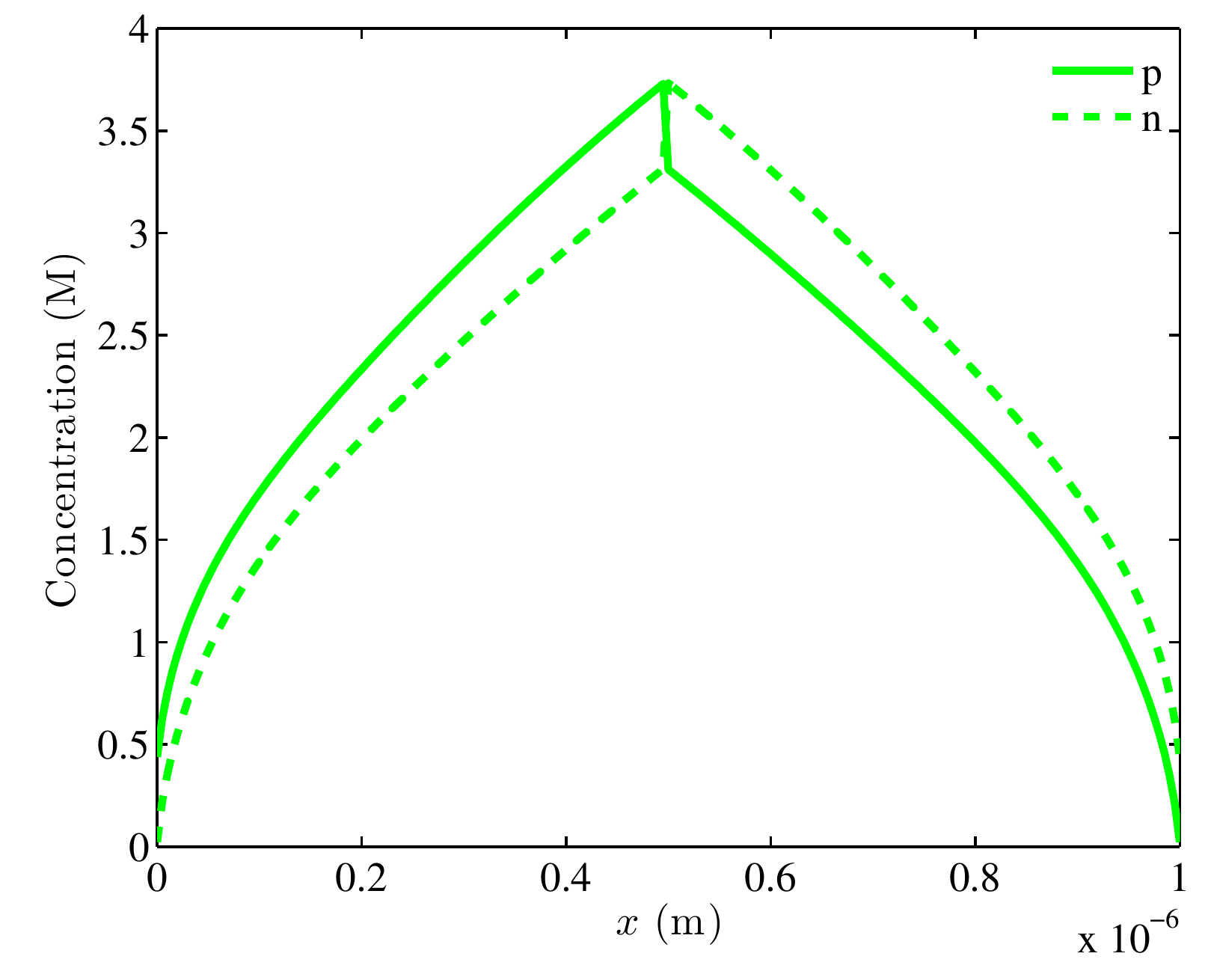}
\includegraphics[width=0.32\textwidth]{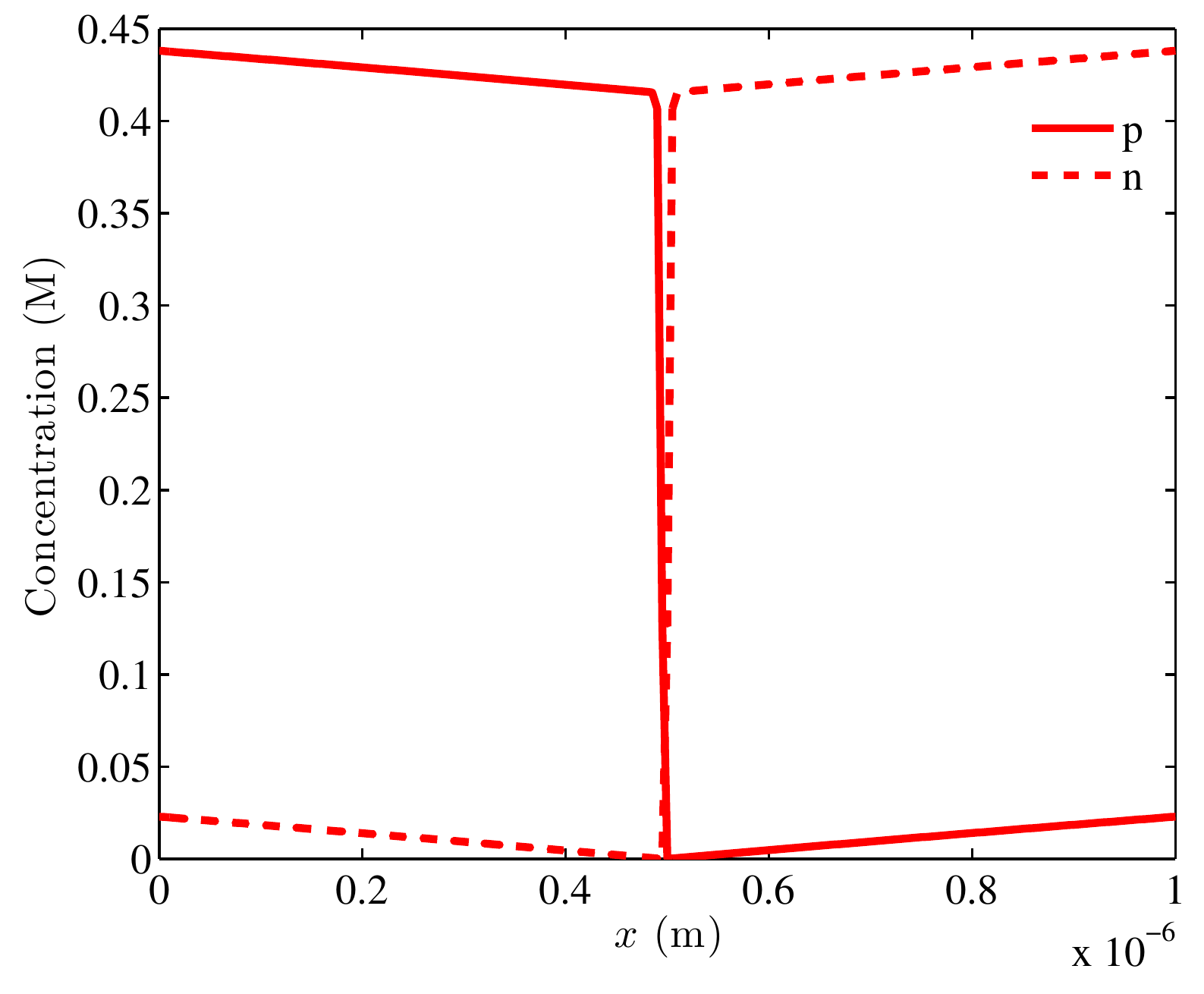}
\caption{BP channel. Permanent charge profile (left). Ion concentrations: forward bias $\Vapp = +1 \, \unit{V}$ (middle),
 reverse bias $\Vapp = -1 \, \unit{V}$ (right).}
\label{fig:BP_Ndop_and_ion_conc}
\end{figure}

\begin{figure}[h!]
\centering
\textbf{Potential and electric field}\par\medskip
\includegraphics[width=0.49\textwidth]{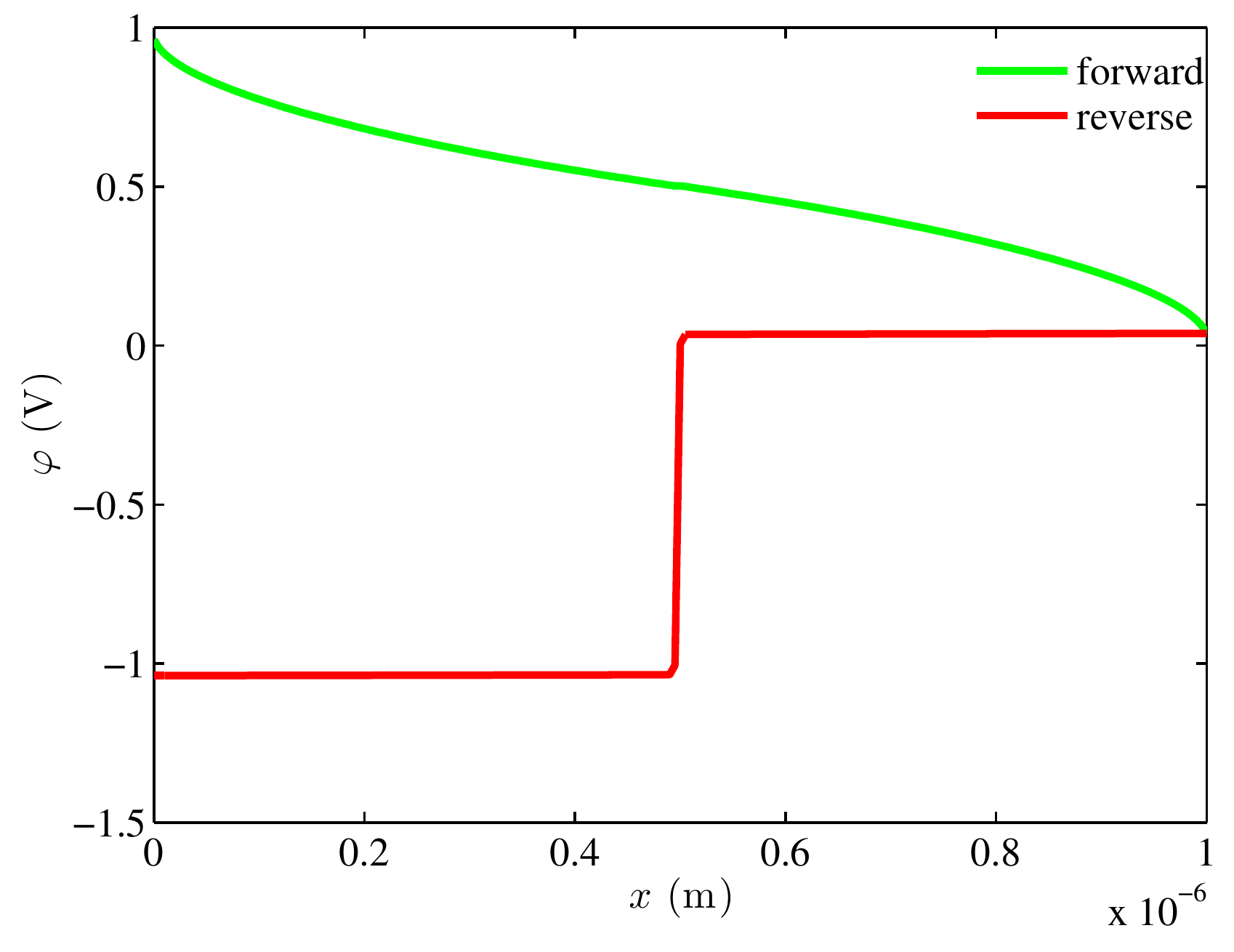}
\includegraphics[width=0.49\textwidth]{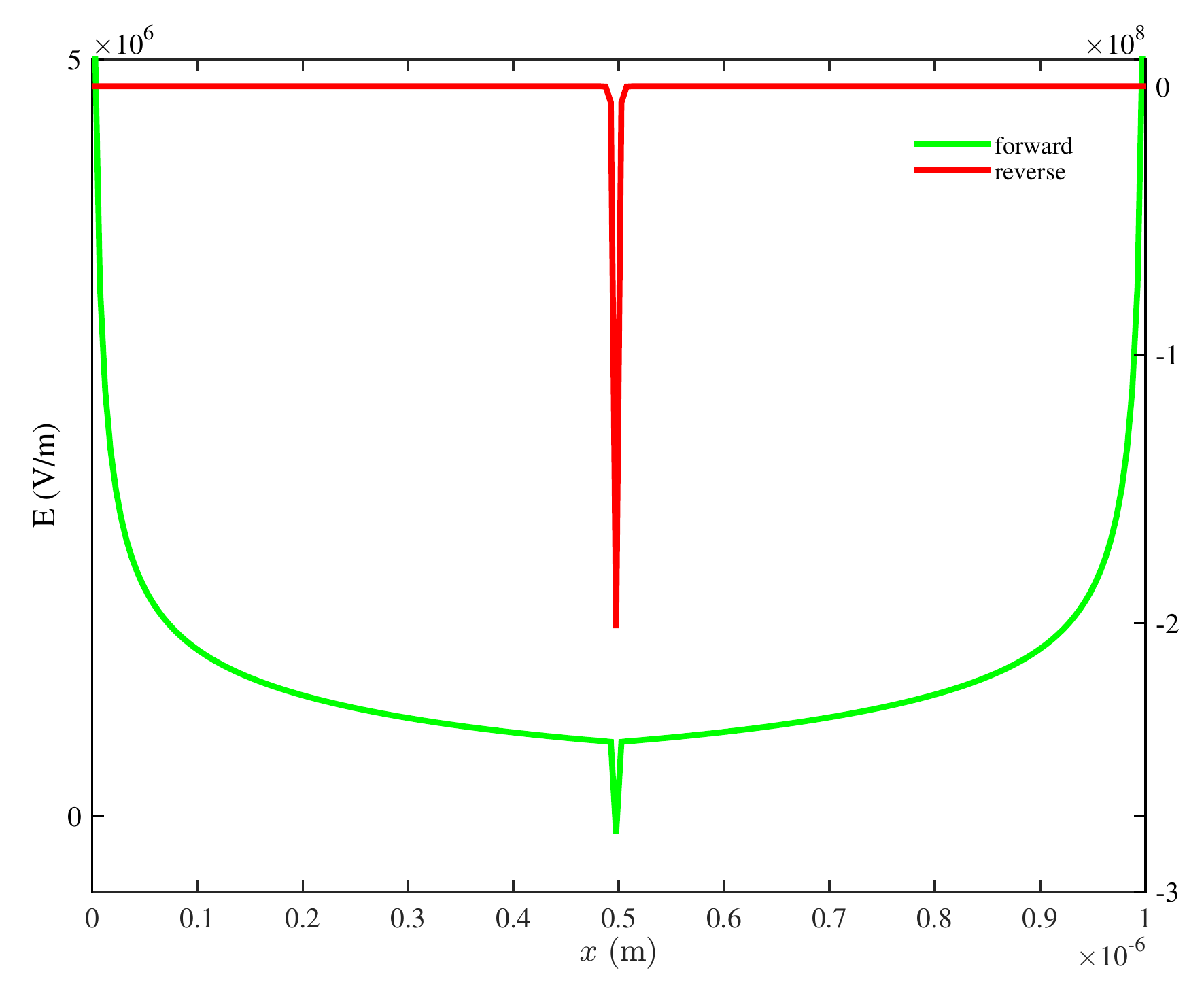}
\caption{BP channel. Electric potential profile: forward (green) and reverse (red) bias. Electric field profiles: forward bias (green, left axis scale) and reverse (red, right axis scale).}
\label{fig:BP_efield}
\end{figure}

The input-output behavior of the BP channel 
is related to the surface permanent charge profile
(cf. Fig.~\ref{fig:channel}).
The channel has a negative surface charge along half of its length, while the remaining part has positive surface charge (of the same magnitude), see Fig.~\ref{fig:BP_Ndop_and_ion_conc} (left).
Thus, anions and cations carry equal weight to channel behavior.
This is confirmed by the symmetric spatial distribution of ion concentrations,
both in the case of forward and reverse bias, see Fig.~\ref{fig:BP_Ndop_and_ion_conc} (middle) and Fig.~\ref{fig:BP_Ndop_and_ion_conc} (right), respectively.
Notice that ion concentration is much higher in forward bias than reverse bias.
Electric potential and electric field profile are shown in Fig.~\ref{fig:BP_efield}.
One can see that, in the closed-state, carrier flow is inhibited by 
the potential barrier at the middle of the channel.
Conversely, in the open-state, the potential drop enhances ion flow.
The computed IV curves show the on-off trend just described.
In particular, we see that if $\Vapp$ is negative the current is almost equal to zero, while it increases in a nonlinear manner if $\Vapp$ is positive 
(cf. Fig.~\ref{fig:BP-IV-THD-T0-Td}). The predicted marked on-off function
mode of the BP channel and the shape of the I-V curves are also in very 
good qualitative agreement with data reported in the study of heat biosensors
in~\cite{caterina1997,voets2004,vay2012}. This fact is an encouraging 
motivation to a further use and calibration of the computational
models and methodologies proposed in the present article on biophysical 
novel applications.

\subsubsection{Temperature}\label{sec:BP-temp}
\begin{figure}[h!]
\centering
\textbf{IV curves when $T(0)=T(d) \neq 300\,\unit{K}$}\par\medskip
\includegraphics[width=0.49\textwidth]
{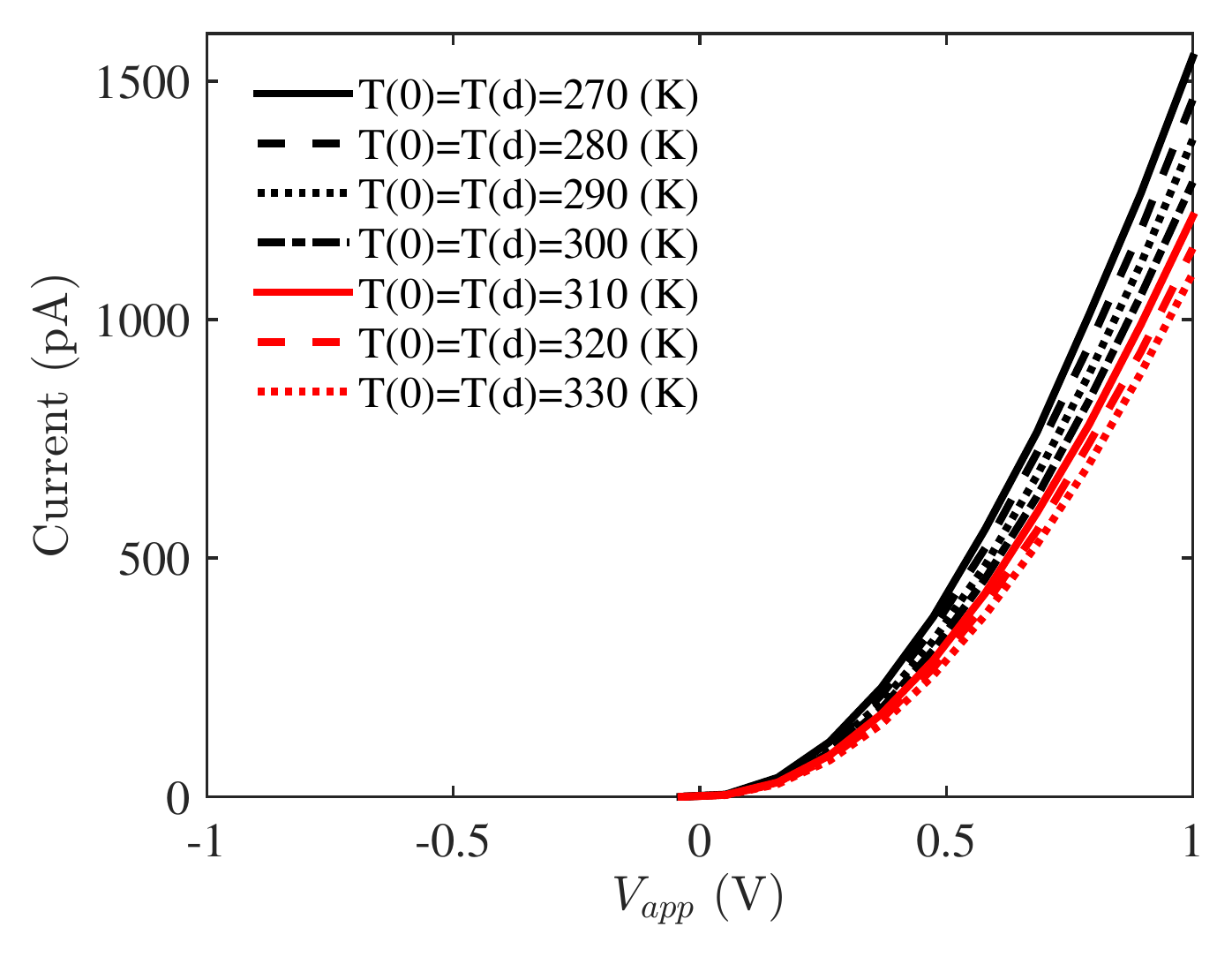}
\includegraphics[width=0.49\textwidth]
{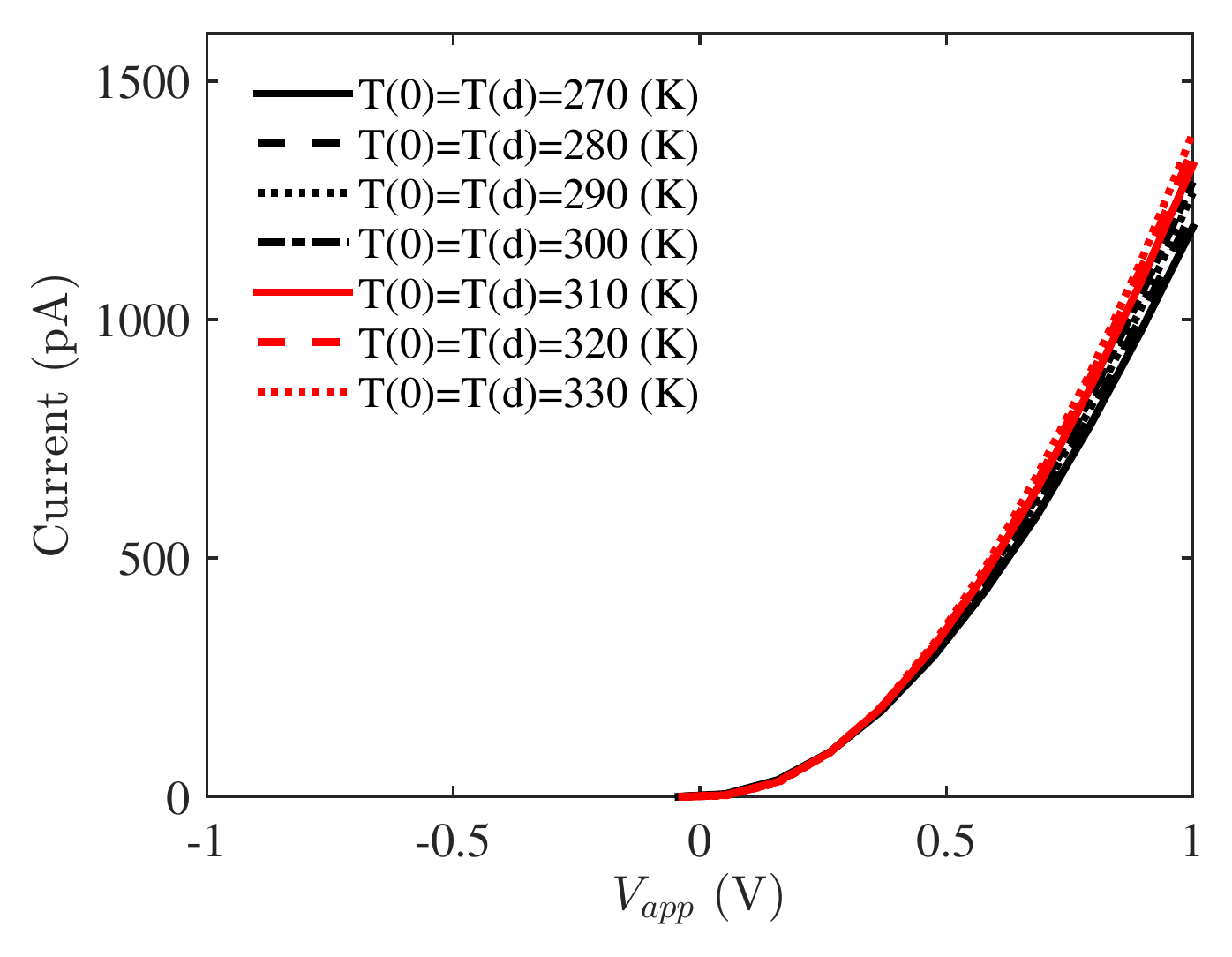}
\caption{BP diode. IV curves when $T(0)=T(d) \neq 300 \, K$: THD 
model (left), ET model (right).}
\label{fig:BP-IV-THD-T0-Td}
\end{figure}

The IV curves, 
when $T(0)=T(d)$ assumes lower or higher values than the reference value of $300\,\unit{K}$, are shown in Fig.~\ref{fig:BP-IV-THD-T0-Td} 
for the THD model (left) and ET model (right), respectively.
The externally applied voltage $\Vapp$ ranges from -1 to 1 V. 
The IV curves computed by the THD model follow the same trend, 
as a function of the bath temperature, as in 
Sect.~\ref{sec:gramicidin} (cf. Fig.~\ref{fig:ballistic_IV_THD_T0_Td}): 
the higher the bath temperature the lower
the current flowing in the channel because of increased frictional effects,
\rs{see Fig~\ref{fig:BP-IV-THD-T0-Td} (left).
The IV curves computed by the ET model in this type of channel 
differ remarkably from those computed by the THD model 
with respect to boundary temperature: the higher
the bath temperature, the higher the current, 
see Fig~\ref{fig:BP-IV-THD-T0-Td} (right). 
We also observe that the spread in the value of the maximum 
channel current as a function of bath temperature is much wider for the THD model than for the ET model. This trend was exactly the opposite in the simulation of the Gramicidin-A channel of Sect.~\ref{sec:gramicidin}, as 
demonstrated by Fig.~\ref{fig:ballistic_IV_THD_T0_Td}. 
We point out that the application to the BP channel of 
the theoretical prediction for channel current given by the ideal diode
model (see, e.g.,~\cite{muller2002device}) would yield
\begin{equation}\label{eq:idealdiode}
I = I_0 \left[ \exp\left( \Frac{q \Vapp}{K_B T_{sys}}  \right) 
- 1 \right]
\end{equation}
where $I$ is the total current flowing in the diode and 
$I_0$ is the saturation current. Thus, according 
to~\eqref{eq:idealdiode}, an increase of channel system temperature 
turns out into a decrease in channel current, in accordance with
the results computed by the THD model. Despite this preliminary indication
in favor of the predictions of the THD model, we feel
that the strongly different response of the two models 
when applied to different channel configurations certainly warrants further investigation and will be the object of a further step of our 
research activity.}

\begin{figure}
\centering
\textbf{$T_e(x)$ profile when $T(0) \neq 300\,\unit{K}$}\par\medskip
\includegraphics[width=0.49\textwidth]{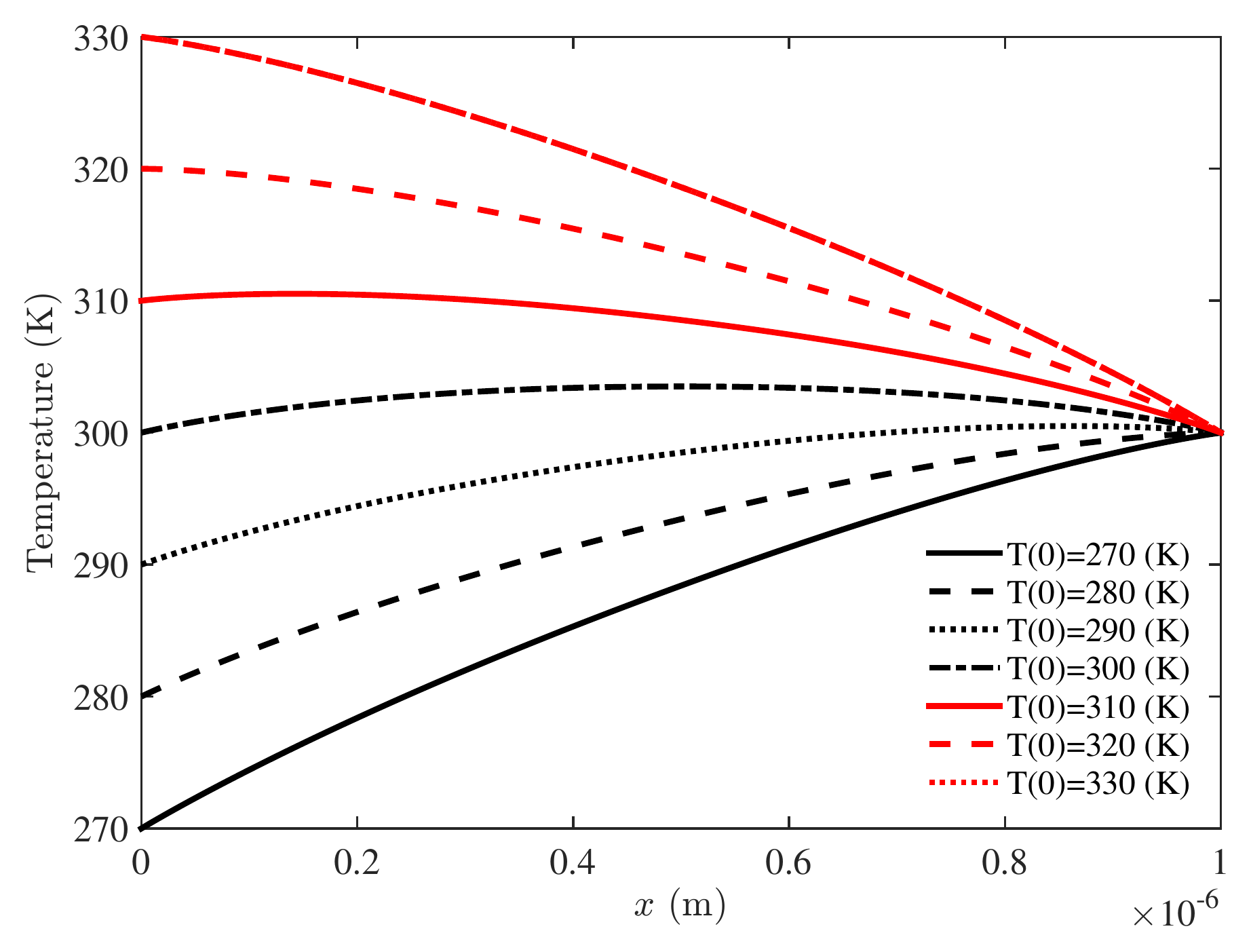}
\includegraphics[width=0.49\textwidth]{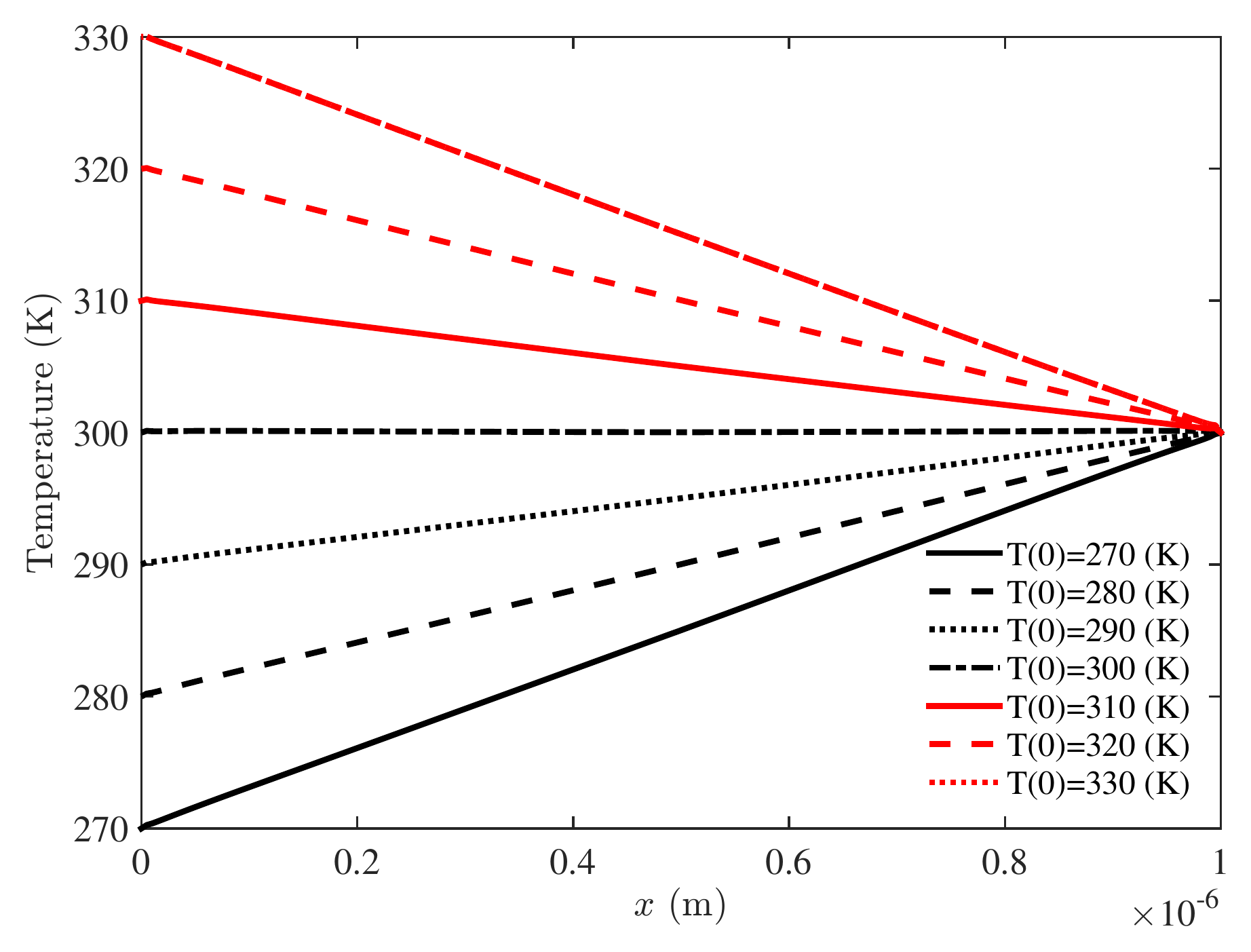}
\caption{BP diode. Water temperature in the ET model (left) and THD model (right).}
\label{fig:BP-Te-ET-T0}
\end{figure}

Also the temperature profiles of electrolytic fluid differ between the ET model and the THD model.
The profiles from ET model have a non linear profile, while those from THD model are linear, see Fig.~\ref{fig:BP-Te-ET-T0} (left) and Fig.~\ref{fig:BP-Te-ET-T0} (right), 
respectively, when $T(0)$ ranges from 270 to 330 K.

%

\subsubsection{Electroosmosis}\label{sec:BP-osmosis}
\begin{figure}
\centering
\textbf{IV curves when $v_e > 0\,\unit{m/s}$}\par\medskip
\includegraphics[width=0.49\textwidth]
{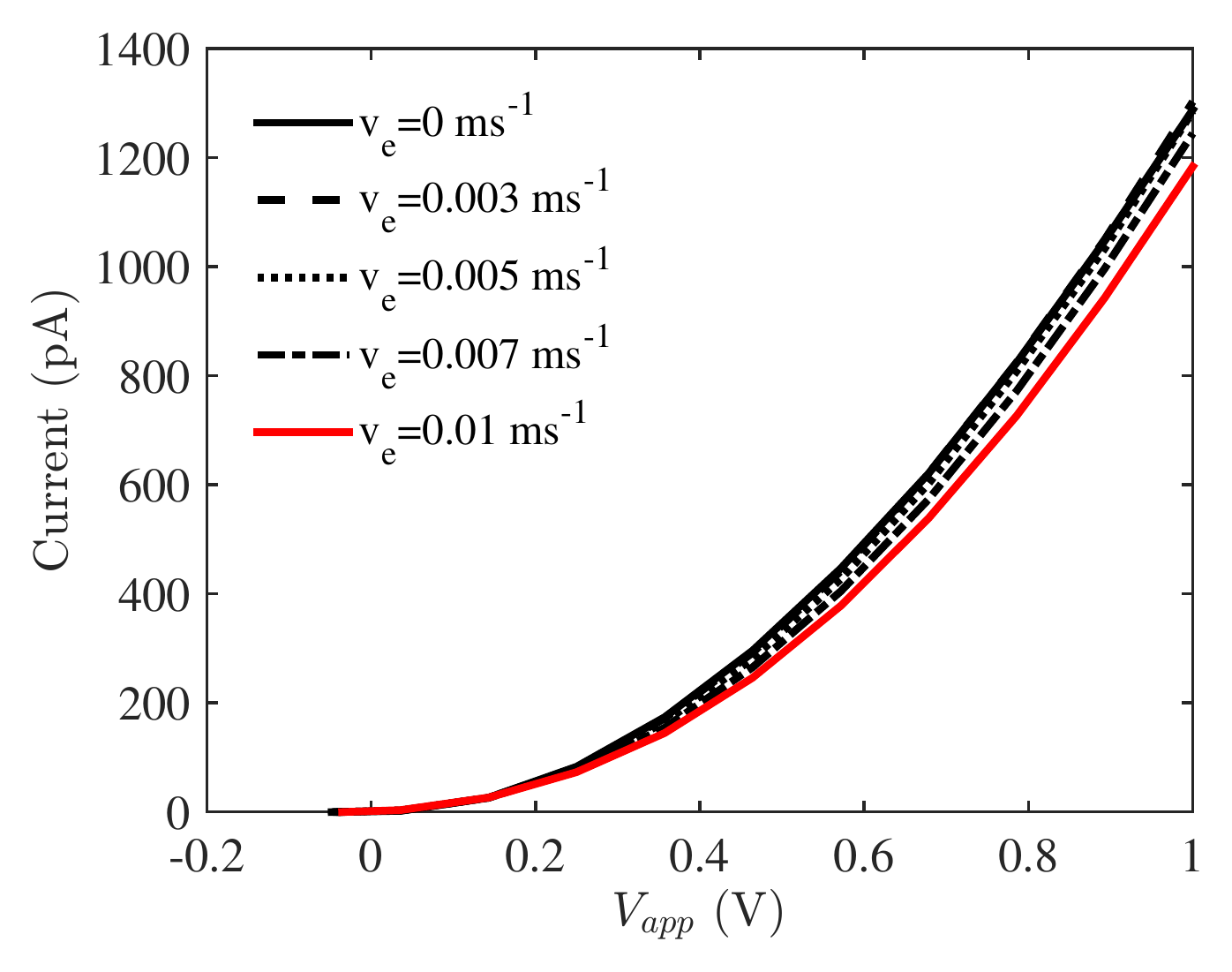}
\includegraphics[width=0.49\textwidth]
{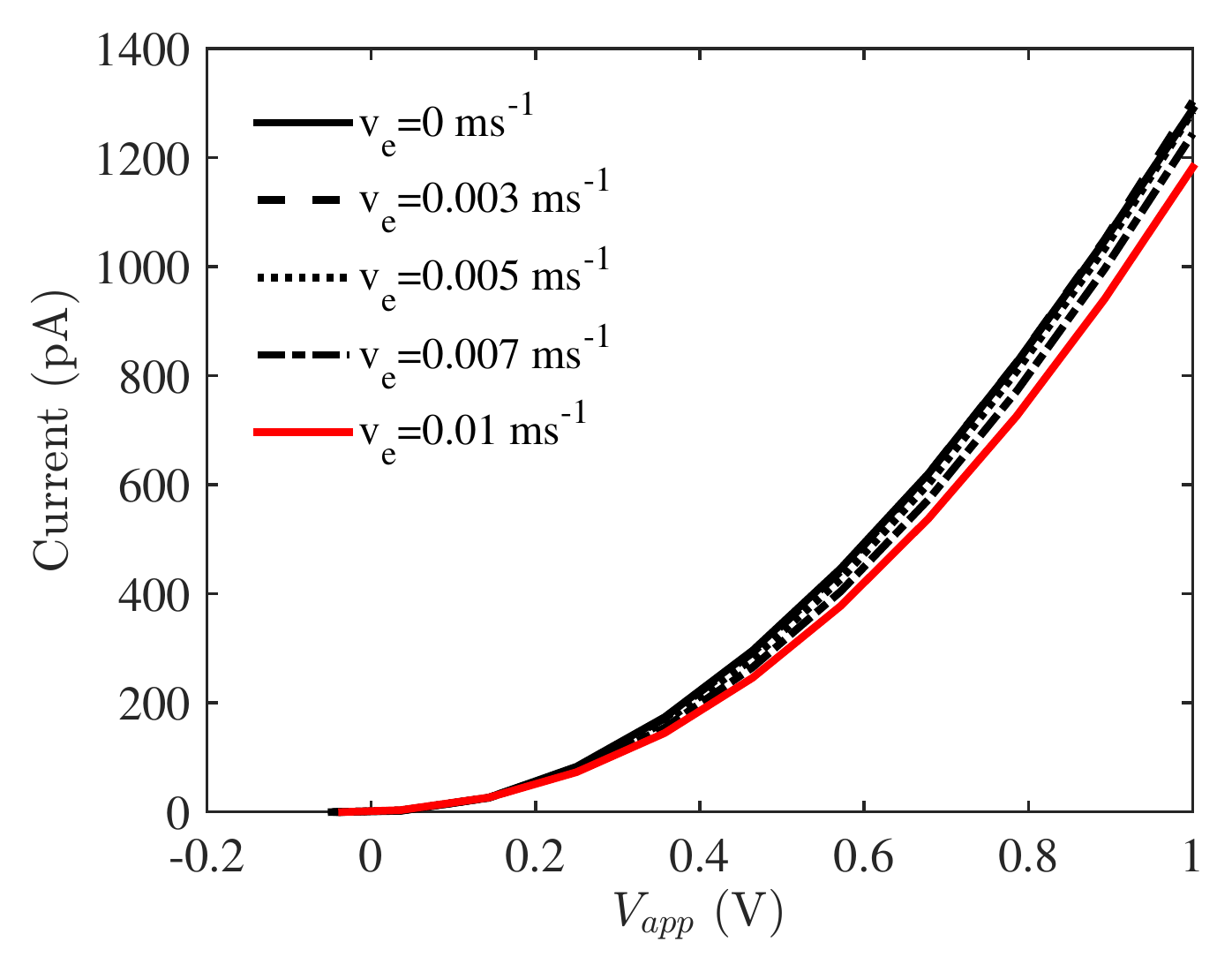}
\caption{BP diode. IV curve when $v_e$ ranges in $0\div 0.01 \, \unit{m/s}$. vTHD model (left). vET model (right).}
\label{fig:BP-IV-THD-uNS}
\end{figure}

In this type of channel the influence of the electrolytic fluid velocity is shown in Fig~\ref{fig:BP-IV-THD-uNS} (left) for the vTHD model 
and in Fig.~\ref{fig:BP-IV-THD-uNS} (right) for the vET model.
The selected range for the electrolyte fluid velocity ($v_e$ from $0$ to $0.01\,\unit{m/s}$) is rather arbitrary since our model does not account,
at the moment, for a self-consistent coupling between the 
vTHD system and the NS equations for the fluid.
As the fluid velocity has positive sign, electrolytic particles flow from left to right.
Accordingly, this additional translational force should enhance the movement of positive ions, while it should reduce that of negative ions, since the two ions move in opposite directions.
From the IV curves shown, one can see that the additional driving force due to $v_e$ lowers the total current flowing in the channel in both models.
Analogous results hold for negative values of $v_e$ (not shown).

\begin{figure}
\centering
\textbf{$T_e(x)$ profile for different $v_e$ (vET vs. vTHD)}\par\medskip
\includegraphics[width=0.32\textwidth]{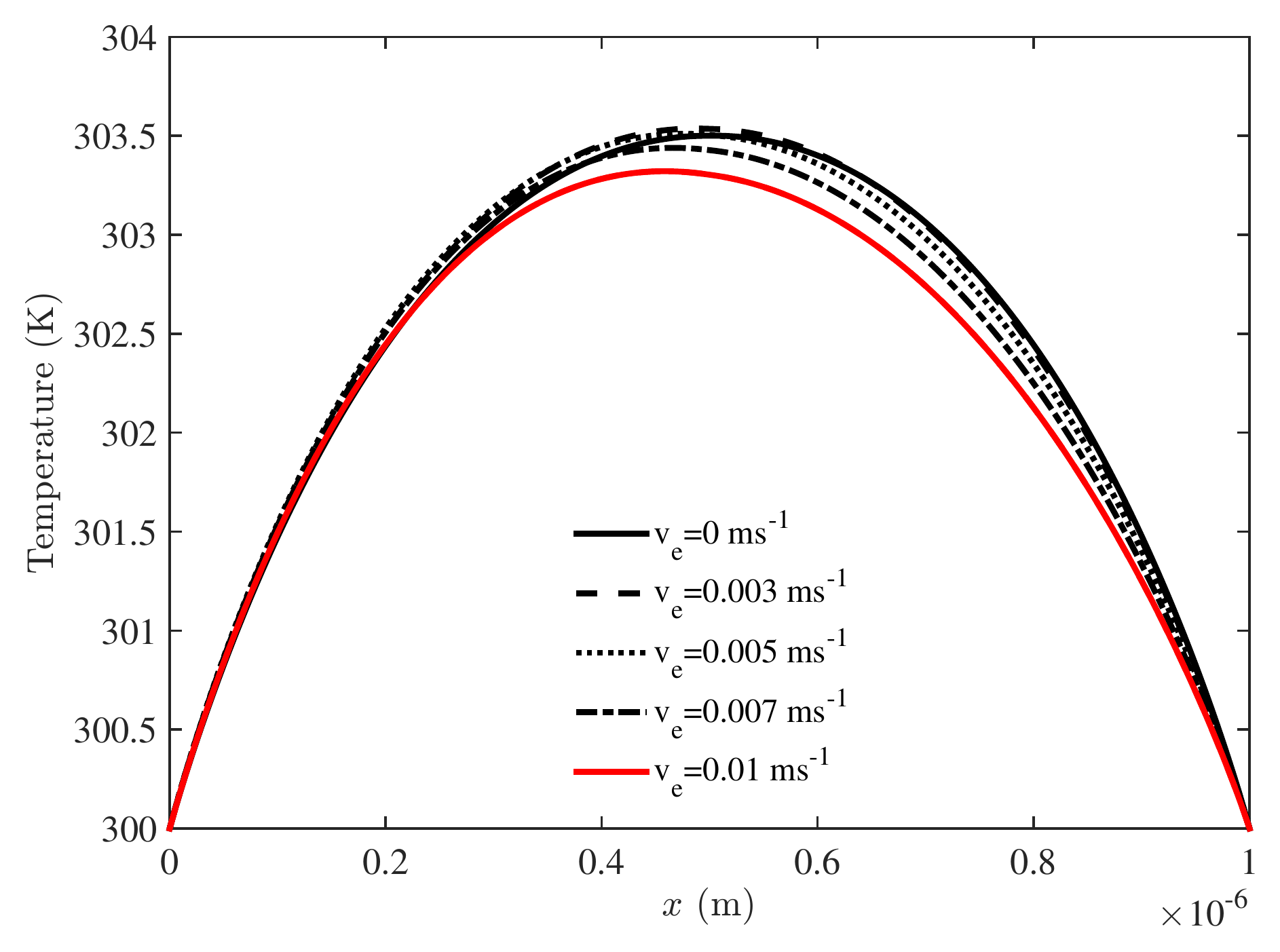}
\includegraphics[width=0.32\textwidth]{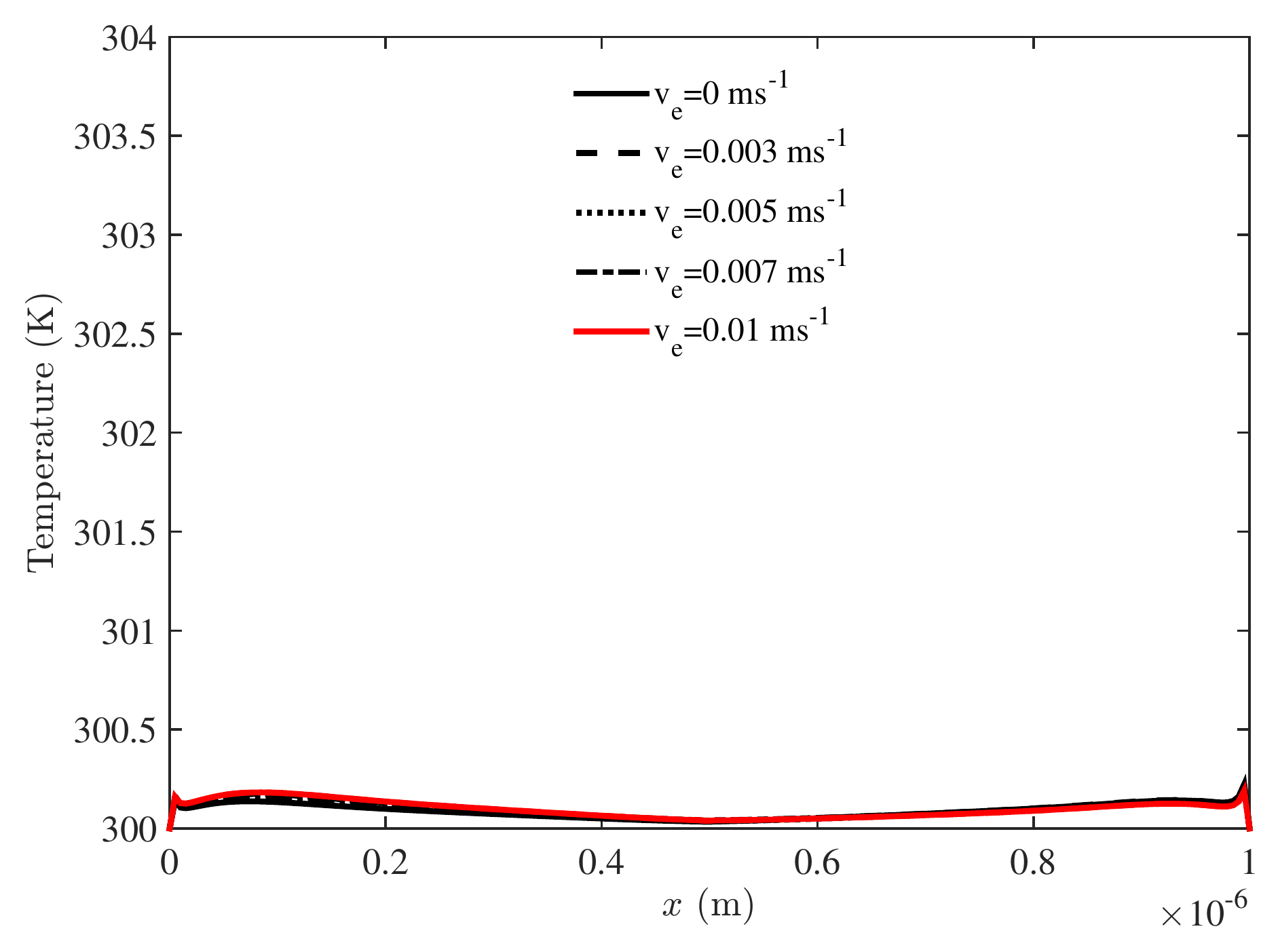}
\includegraphics[width=0.32\textwidth]{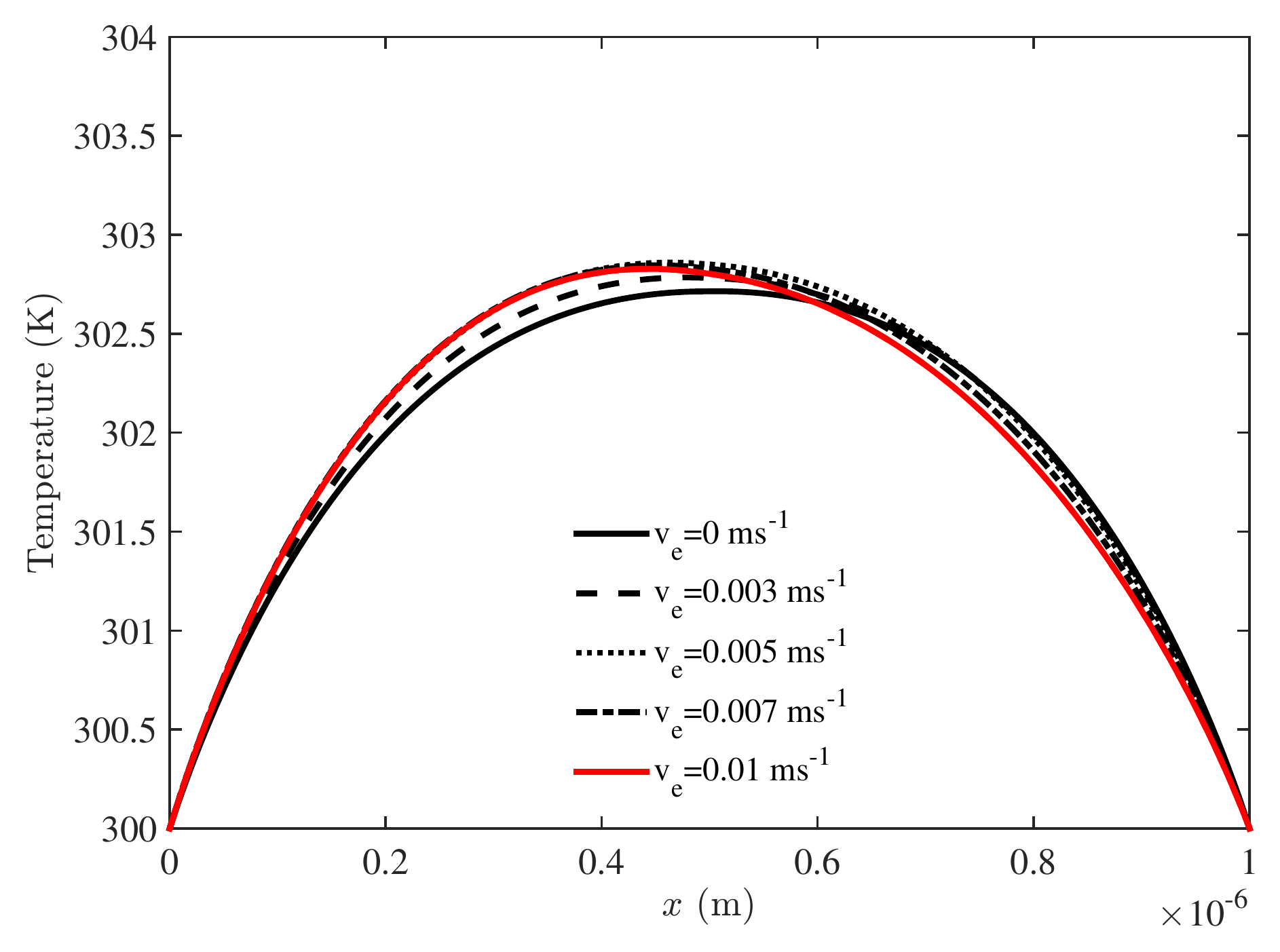}
\caption{BP diode. $v_e$ ranges in $0 \div 0.01\,ms^{-1}$. vET model (left) differs greatly from vTHD at $\vsat=10\,\unit{m/s}$ (middle), while is quite close to vTHD at $\vsat=1\,\unit{m/s}$.}
\label{fig:BP_T_ET_uNS}
\end{figure}

The temperature of the electrolytic fluid velocity for different 
values of $v_e$ is shown in Fig.~\ref{fig:BP_T_ET_uNS}. Temperature is only partially affected by $v_e$.
The increase in temperature predicted by the vET model (left) is about 
3 $\unit{K}$, which is quite remarkable, especially if compared to the 
results obtained with the vTHD model (middle), 
where temperature changes are much less significant, instead.
The very different prediction of the vTHD formulation is
related to the parameter $\vsat$ that strongly affects
the value of the relaxation time parameter $\tau_{w_{\nu}}$ 
(cf.~\eqref{eq:tau_w_vsat}) 
in the energy balance equation~\eqref{eq:THD1D_en}: the larger
$\vsat$, the smaller $\tau_{w_{\nu}}$, which corresponds to almost
instantaneous restoration of equilibrium conditions in the electrolyte.
Indeed, taking a smaller value of the saturation velocity, for example $\vsat=1\,\unit{m/s}$, leads to a result that is much closer 
to the thermal profile predicted by the vET model, see 
Fig.~\ref{fig:BP_T_ET_uNS} (right). 


%

\section{Conclusions and Research Perspectives}\label{sec:conclusions}
In the present article we have proposed and numerically investigated 
a hierarchy of mathematical models for the simulation of thermal,
fluid and electrochemical phenomena in biological transmembrane channels.
The hierarchy is an extension of the classic Poisson-Nernst-Planck
model for ion electrodiffusion and is conducted along the same lines
of thought that have guided the development of the so-called
hydrodynamic transport model in the analysis of semiconductor devices.
To discretize the proposed models we have devised in the 1D case 
a robust finite element dual-mixed hybridized method that ensures flux conservation, self-equilibrium
and satisfaction of a positivity principle for ion concentrations and temperatures. The numerical scheme has been thoroughly studied in several
benchmark problems that demonstrate its accuracy and stability.
An appropriate solution map is used to successively solve the nonlinear
system of equations arising from model hierarchy and the resulting 
computational tool has been successfully calibrated and validated in the simulation
of two realistic biological channels. Future developments of this study
include:
\begin{itemize}
\item time dependent simulations to describe the response of the
channel to externally applied stimuli;
\item self-consistent coupling of the hierarchy with the solution of 
the Navier-Stokes equations for the electrolyte fluid;
\item deeper investigation of the dependence of model predictions
on biophysical parameters, for instance, saturation velocity that 
seems to play a critical role in determining the self-heating
effect in the electrolyte fluid;
\item extension of the numerical scheme to 2D and 3D channel simulation.
\end{itemize}

\bibliographystyle{plain} 
\bibliography{biblio}

\end{document}